%% file: main.tex
\documentclass{article}

\usepackage[final,nonatbib]{neurips_2021}

\usepackage[utf8]{inputenc} 
\usepackage[T1]{fontenc}    
\usepackage{hyperref}       
\usepackage{url}            
\usepackage{booktabs}       
\usepackage{amsfonts}       
\usepackage{nicefrac}       
\usepackage{microtype}      
\usepackage{xcolor}         

\usepackage{amsmath, amsthm, amssymb}
\usepackage{algorithm, algorithmicx}
\usepackage[noend]{algpseudocode}
\usepackage{subfiles}
\usepackage{graphicx}
\usepackage{subfigure}
\usepackage{enumitem}
\setlist[itemize]{topsep = 0pt, itemsep = 0pt, partopsep = 0pt, itemsep = 0pt, leftmargin = *}
\setlist[enumerate]{topsep = 0pt, itemsep = 0pt, partopsep = 0pt, itemsep = 0pt, leftmargin = *}
\setlist[description]{topsep = 0pt, itemsep = 0pt, partopsep = 0pt, itemsep = 0pt, leftmargin = *}
\usepackage[margin = 0pt]{caption}
\usepackage{mathtools}
\mathtoolsset{showonlyrefs}
\allowdisplaybreaks

\newtheorem{theorem}{Theorem}
\newtheorem{corollary}[theorem]{Corollary}
\newtheorem{lemma}[theorem]{Lemma}

\theoremstyle{definition}
\newtheorem{assumption}[theorem]{Assumption}
\newtheorem{definition}[theorem]{Definition}
\theoremstyle{remark}

\DeclareMathOperator*{\argmin}{argmin}
\DeclareMathOperator*{\argmax}{argmax}
\newcommand{\st}{\ \mathop{\mathrm{s.t.}}\ }

\newcommand{\bbN}{\mathbb{N}}
\newcommand{\bbR}{\mathbb{R}}

\DeclareMathOperator{\bigO}{O}
\DeclareMathOperator{\dom}{dom}
\DeclareMathOperator{\intmath}{int}

\title{A Gradient Method for Multilevel Optimization}

\author{%
	Ryo Sato\\
	The University of Tokyo\\
	\And
	Mirai Tanaka\\
	The Institute of Statistical Mathematics\\
	RIKEN\\
	\And
	Akiko Takeda\\
	The University of Tokyo\\
	RIKEN\\
}

\begin{document}
\maketitle

\begin{abstract}
Although application examples of multilevel optimization have already been discussed since the 1990s,
the development of solution methods was almost limited to bilevel cases due to the difficulty of the problem.
In recent years, in machine learning, Franceschi et al.\ have proposed a method for solving bilevel optimization problems
by replacing their lower-level problems with the $T$ steepest descent update equations with some prechosen iteration number $T$.
In this paper, we have developed a gradient-based algorithm for multilevel optimization with $n$ levels based on their idea
and proved that our reformulation asymptotically converges to the original multilevel problem.
As far as we know, this is one of the first algorithms with some theoretical guarantee for multilevel optimization.
Numerical experiments show that
a trilevel hyperparameter learning model considering data poisoning
produces more stable prediction results
than an existing bilevel hyperparameter learning model in noisy data settings.
\end{abstract}

\subfile{chapters/Introduction}

\subfile{chapters/Preliminaries}

\subfile{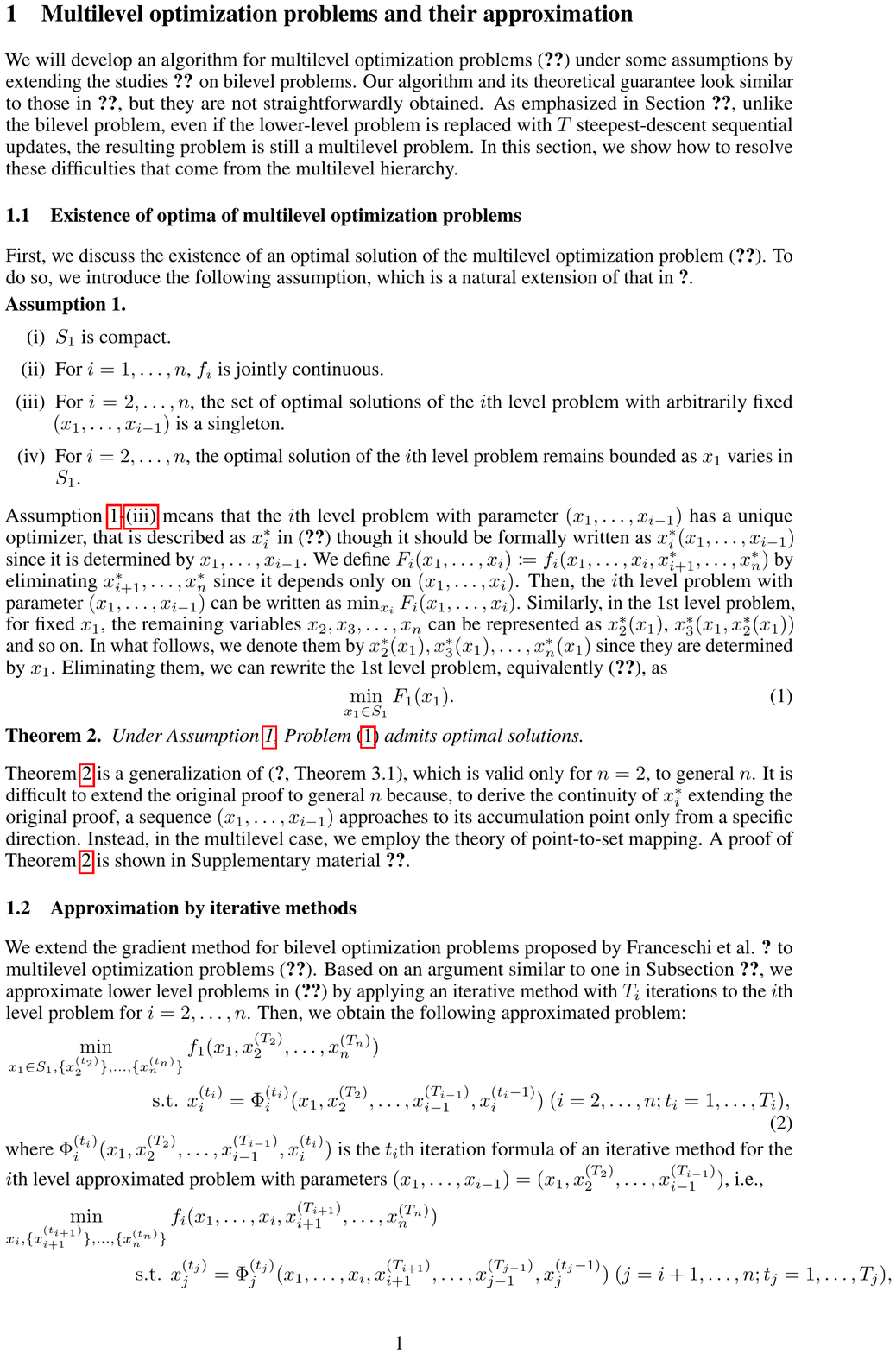}

\subfile{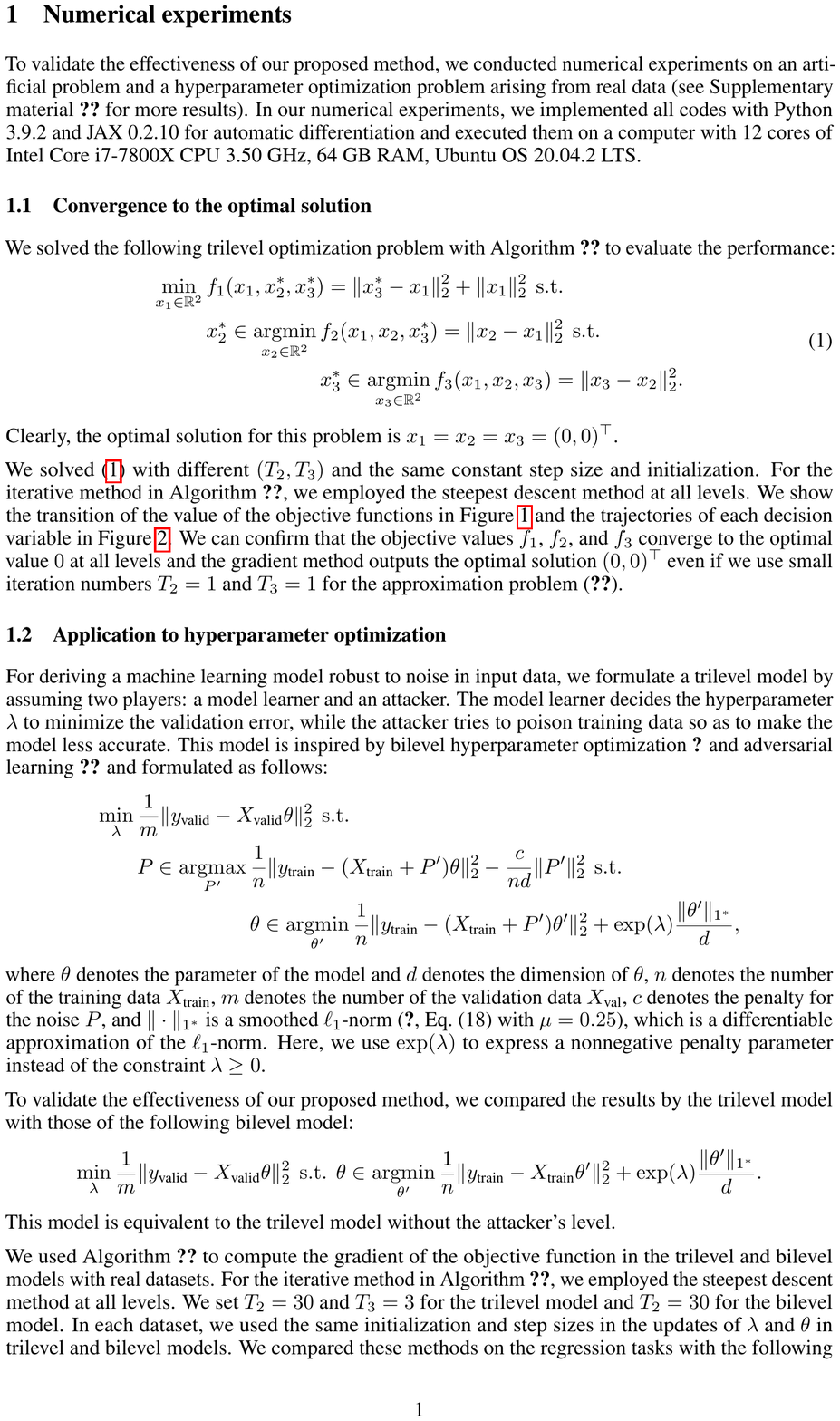}

\subfile{chapters/Conclusion}

\begin{ack}
This work was partially supported by JSPS KAKENHI (JP19K15247, JP17H01699, and JP19H04069).
\end{ack}

\bibliographystyle{abbrv}
\bibliography{reference}

\newpage
\appendix

\subfile{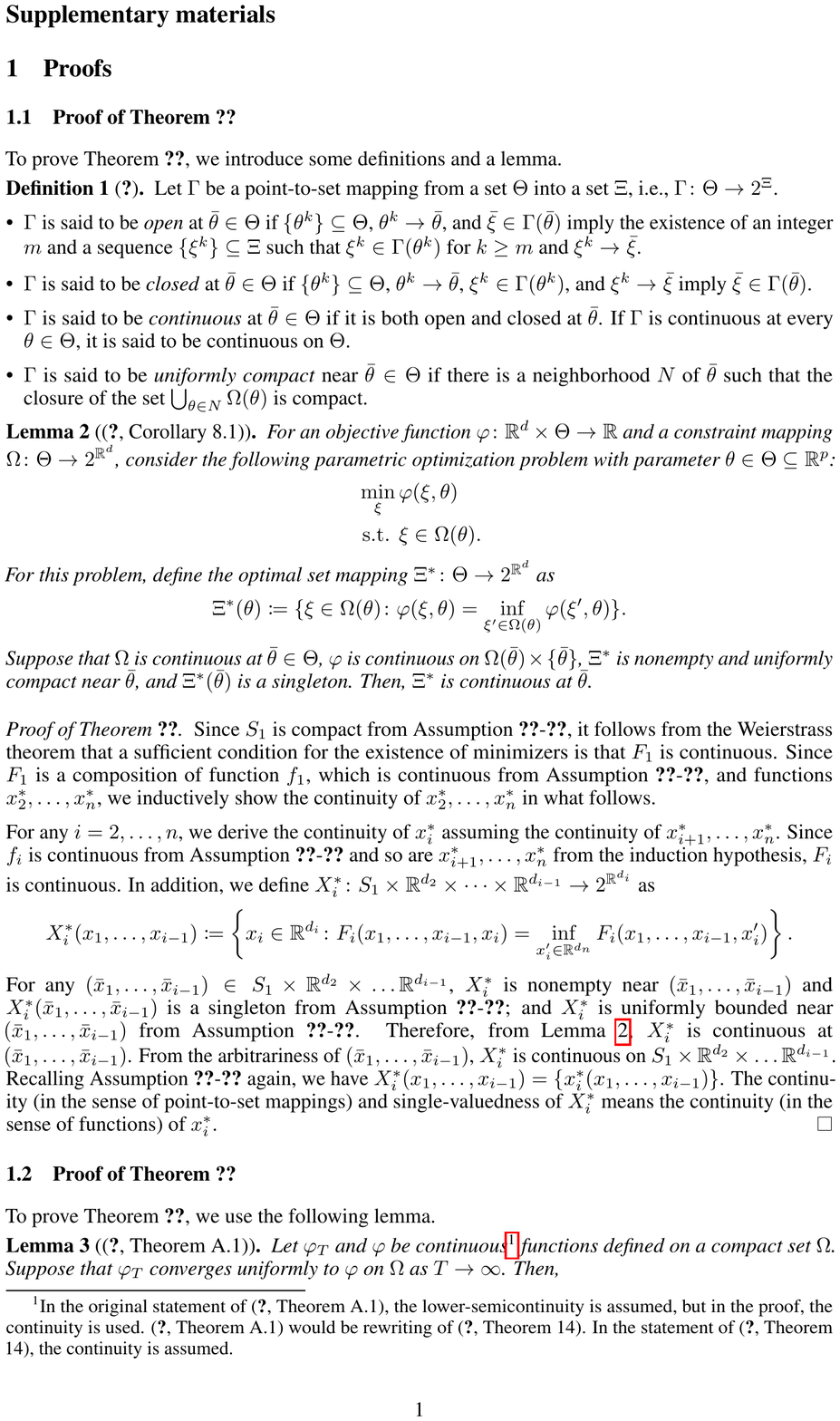}
\end{document}

%% file: chapters/Introduction.tex
\section{Introduction}
\label{sect: intro}

\paragraph{Multilevel optimization}
When modeling real-world problems,
there is often a hierarchy of decision-makers, and decisions are taken at different levels in the hierarchy.
In this work, we consider multilevel optimization problems whose simplest form%
\footnote{This assumes that lower-level problems have unique optimal solutions
$x_{2}^{\ast}, \dots, x_{n}^{\ast}$ for the simplicity.
}
is the following:
\begin{equation}
\begin{aligned}
\min_{x_{1} \in S_{1} \subseteq \bbR^{d_{1}}, x_{2}^{\ast}, \dots, x_{n}^{\ast}}{}&
	f_{1}(x_{1}, x_{2}^{\ast},\dots, x_{n}^{\ast}) \st\\
&	\begin{aligned}
	x_{2}^{\ast} = \argmin_{x_{2} \in \bbR^{d_{2}}, x_{3}^{\ast}, \dots, x_{n}^{\ast}}{}&
		f_{2}(x_{1}, x_{2}, x_{3}^{\ast},\dots, x_{n}^{\ast}) \st\\
		\ddots\\
		&	\begin{aligned}
			&	x_{n}^{\ast} = \argmin_{x_{n}\in \bbR^{d_{n}}} f_{n}(x_{1}, x_{2}, \dots, x_{n}),
			\end{aligned}
	\end{aligned}
\end{aligned}
\label{eq: multilevel}
\end{equation}
where $d_{i}$ is positive integers,
$x_{i} \in \bbR^{d_{i}}$ is the decision variable at the $i$th level,
and $f_{i}\colon \bbR^{d_{i}} \to \bbR$ is the objective function of the $i$th level for $i = 1, \dots, n$.
In the formulation, an optimization problem contains another optimization problem as a constraint.
The framework can be used to formulate problems in which decisions are made in sequence and earlier decisions influence later decisions.
An optimization problem in which this structure is $n$-fold
is called an $n$-level optimization problem or a multilevel optimization problem especially when $n \ge 3$.

\paragraph{Existing study on bilevel optimization}
Especially, when $n = 2$ in \eqref{eq: multilevel},
it is called a bilevel optimization problem.
Bilevel optimization or multilevel optimization has been a well-known problem in the field of optimization
since the 1980s or 1990s respectively (see the survey paper \cite{VC94} for the research at that time),
and their solution methods have been studied mainly on bilevel optimization until now.
Recently, studies on bilevel optimization have received significant attention
from the machine learning community
due to practical applications and the development of efficient algorithms.
For example, bilevel optimization was used to formulate decision-making problems for players in conflict \cite{WC17, SMD18},
and to formulate hyperparameter learning problems in machine learning \cite{BHJKP06, FDFP17, FFSGP18}.
Franceschi et al.\ \cite{FDFP17} constructed a single-level problem that approximates a bilevel problem
by replacing the lower-level problem with $T$ steepest descent update equations using a prechosen iteration number $T$.
It is shown that the optimization problem with $T$ steepest descent update equations asymptotically converges to
the original bilevel problem as $T \rightarrow \infty$.
It is quite different from the ordinary approach (e.g., \cite{AS13})
that transforms bilevel optimization problems into single-level problems using the optimality condition of the lower-level problem.
The application of the steepest descent method to bilevel optimization problems has recently attracted attention,
resulting in a stream of papers; e.g., \cite{FNPH19, SCHB19, LVD20}.


\paragraph{Existing study on multilevel optimization}
Little research has been done on decision-making models \eqref{eq: multilevel} that assume the multilevel hierarchy,
though various applications of multilevel optimization have been already discussed in the 1990s \cite{VC94}.
As far as we know, few solution methods with a theoretical guarantee have been proposed.
A metaheuristic algorithm based on a perturbation, i.e.,
a solution method in which the neighborhood is randomly searched for a better solution
in the order of $x_{1}$ to $x_{n}$ has been proposed in \cite{TKMO12} for general $n$-level optimization problems.
The effect of updating $x_{i}$ on the decision variables $x_{i + 1}, \dots, x_{n}$ of the lower-level problems is not considered,
and hence, the hierarchical structure of the multilevel optimization problem cannot be fully utilized in updating the solution.
Because of this, the algorithm does not have a theoretical guarantee for the obtained solution.
Very recently (two days before the NeurIPS 2021 abstract deadline),
a proximal gradient method based on the fixed-point theory for trilevel optimization problem \cite{SKM21} has been proposed.
This paper has proved convergence to an optimal solution assuming the convexity of the objective functions.
However, it does not include numerical results and hence its practical efficiency is unclear.

\paragraph{Contribution of this paper}
In this paper, by extending the gradient method for bilevel optimization problems \cite{FDFP17} to multilevel optimization problems,
we propose an algorithm with a convergence guarantee for multilevel optimization other than one by \cite{SKM21}.
Increasing the problem hierarchy from two to three drastically makes developing algorithms and showing theoretical guarantees difficult.
In fact, as discussed above, in the case of bilevel,
replacing the second-level optimization problem with its optimality condition
or replacing it with $T$ steepest-descent sequential updates immediately results in a one-level optimization problem.
However, when it comes to $n \ge 3$ levels,
it may be necessary to perform this replacement $n$ times,
and it seems not easy to construct a solution method with a convergence guarantee.
It is not straightforward at all to give our method a theoretical guarantee similar to
the one \cite{FDFP17} developed for bilevel optimization.
We have confirmed the effectiveness of the proposed method by applying it to hyperparameter learning with real data.
Experimental verifications of the effectiveness of trilevel optimization using real-world problems are the first ones as far as we know.

%% file: chapters/Preliminaries.tex
\section{Related existing methods}

\subsection{Two-stage robust optimization problems}

When making long-term decisions, it is also necessary to make many decisions according to the situation.
If the same objective function is acceptable throughout the period, i.e., $f_{1} = \dots = f_{n}$,
and there are hostile players, such a problem is often formulated as
a multistage robust optimization problem.
Especially, the concept of two-stage robust optimization (so-called  adjustable robust optimization)
$\min_{x_{1} \in S_{1}} \max_{x_{2} \in S_{2}} \min_{x_{3} \in S_{3}} f_{1}(x_{1}, x_{2}, x_{3})$
was introduced with feasible sets $S_{i} \subseteq \bbR^{d_{i}}$ for $i = 1, 2, 3$,
and methodology based on affine decision rules was developed by \cite{BGGN04}.
Researches based on affine policies are still being actively conducted; see, e.g., \cite{BR16, YGH19}.
The two-stage robust optimization can be considered as a special case of \eqref{eq: multilevel} with $n = 3$
when constraints $x_{2} \in S_{2} \subseteq \bbR^{d_{2}}$ and $x_{3} \in S_{3} \subseteq \bbR^{d_{3}}$ are added.
This approach is unlikely to be applicable to our problem \eqref{eq: multilevel}
without strong assumptions such as affine decision rules.

There is another research stream on min-max-min problems including integer constraints (see e.g., \cite{BCP14, WC17}).
They transform the problem into the min-max problem by taking the dual for the inner ``min''
and apply cut-generating approaches or Benders decomposition techniques with the property of integer variables.
While the approach is popular in the application studies of electric grid planning,
it is restricted to the specific min-max-min problems and no more applicable to multilevel optimization problems \eqref{eq: multilevel}.


\subsection{Existing methods for bilevel optimization problems}
\label{sect: bilevel}

Various applications are known for bilevel optimization problems, i.e., the case of $n = 2$ in \eqref{eq: multilevel}:
\begin{equation}
\min_{x_{1} \in S_{1}} f_{1}(x_{1}, x_{2}) \st x_{2} = \argmin_{x_{2}'} f_{2}(x_{1}, x_{2}').
\label{eq: bilevel 1st}
\end{equation}
One of the most well-known applications in machine learning is hyperparameter learning.
For example, the problem that finds the best hyperparameter value in the ridge regression is formulated as
\begin{equation}
\min_{\lambda \ge 0}{} \|y_{\text{valid}} - X_{\text{valid}} \theta\|_{2}^{2}
\st \theta = \argmin_{\theta'} \|y_{\text{train}} - X_{\text{train}} \theta'\|_{2}^{2} + \lambda \|\theta'\|_{2}^{2},
\end{equation}
where $(X_{\text{train}}, y_{\text{train}})$ are training samples,
$(X_{\text{valid}}, y_{\text{valid}})$ are validation samples,
and $\lambda$ is a hyperparameter that is optimized in this problem.
Under Assumption \ref{asp: existence} restricted to $n = 2$ shown later,
the existence of an optimal solution to the bilevel problem \eqref{eq: bilevel 1st} is ensured by \cite[Theorem 3.1]{FFSGP18}.

There are mainly two approaches to solve bilevel optimization problems.
The old practice is to replace the lower-level problem with its optimality condition and to solve the resulting single-level problem.
It seems difficult to use this approach to multilevel problems because we need to apply the replacement $n$-times.
The other one is to use gradient-based methods developed by Franceschi et al. \cite{FDFP17, FFSGP18} for bilevel optimization problems.
Their approach reduces \eqref{eq: bilevel 1st} to a single-level problem by
replacing the lower-level problem with $T$ equations using a prechosen number $T$.

Hereinafter, we briefly summarize the results of Franceschi et al. \cite{FDFP17, FFSGP18} for bilevel optimization.
Under the continuous differentiability assumption for $f_{2}$,
$x_{2}$ is iteratively updated by $x_{2}^{(t)} = \Phi^{(t)}(x_{1}, x_{2}^{(t - 1)})$ at the $t$th iteration
using an iterative method, e.g., the gradient descent method:
\begin{equation}
\Phi^{(t)}(x_{1}, x_{2}^{(t - 1)}) = x_{2}^{(t - 1)} - \alpha^{(t)} \nabla_{x_{2}} f_{2}(x_{1}, x_{2}^{(t - 1)}).
\end{equation}
Then, the bilevel optimization problem is approximated by the single-optimization problem with $T$ equality constraints
and new variables $\{x_{2}^{(t)}\}_{t = 1}^{T}$ instead of $x_{2}$:
\begin{equation}
\min_{x_{1} \in S_{1}, \{x_{2}^{(t)}\}}
	f_{1}(x_{1}, x_{2}^{(T)})
\st
	x_{2}^{(t)} = \Phi^{(t)}(x_{1}, x_{2}^{(t - 1)})\ (t = 1, \dots, T),
\label{eq: bilevel approx}
\end{equation}
where $x_{2}^{(0)}$ is a given constant.
Eliminating $x_{2}^{(1)}, \dots, x_{2}^{(T)}$ using constraints,
we can equivalently recast the problem above into the unconstrained problem $\min_{x_{1} \in S_{1}} \tilde{F}_{1}(x_{1})$.

Under assumption on differentiability,
the gradient of $\tilde{F}_{1}(x_{1})$ is given in \cite[Section 3]{FDFP17} by
\begin{align}
\nabla_{x_{1}} \tilde{F}_{1}(x_{1})
&= \nabla_{x_{1}} f_{1}(x_{1}, x_{2}^{(T)})
+ \sum_{t = 1}^{T} B^{(t)} \left(\prod_{s = t + 1}^{T} A^{(s)}\right) \nabla_{x_{2}} f_{1}(x_{1}, x_{2}^{(T)}),\\
A^{(t)}
&= \nabla_{x_{2}} \Phi^{(t)}(x_{1}, x_{2}^{(t - 1)})\ (t = 1, \dots, T),\\
B^{(t)}
&= \nabla_{x_{1}} \Phi^{(t)}(x_{1}, x_{2}^{(t - 1)})\ (t = 1, \dots, T).
\end{align}
Roughly speaking, \cite[Theorem 3.2]{FFSGP18} proved that
the optimal value and the solution set of \eqref{eq: bilevel approx} converge to those of \eqref{eq: bilevel 1st} as $T \to \infty$
under Assumptions \ref{asp: existence} and \ref{asp: convergence} restricted to $n = 2$.

%% file: chapters/Algorithms.tex
\section{Multilevel optimization problems and their approximation}
\label{sect: approx}

We will develop an algorithm for multilevel optimization problems \eqref{eq: multilevel} under some assumptions
by extending the studies \cite{FDFP17, FFSGP18} on bilevel problems.
Our algorithm and its theoretical guarantee look similar to those in \cite{FDFP17, FFSGP18},
but they are not straightforwardly obtained.
As emphasized in Section \ref{sect: intro},
unlike the bilevel problem,
even if the lower-level problem is replaced with $T$ steepest-descent sequential updates,
the resulting problem is still a multilevel problem.
In this section, we show how to resolve these difficulties that come from the multilevel hierarchy.

\subsection{Existence of optima of multilevel optimization problems}

First, we discuss the existence of an optimal solution of the multilevel optimization problem \eqref{eq: multilevel}.
To do so, we introduce the following assumption,
which is a natural extension of that in \cite{FFSGP18}.

\begin{assumption}
\label{asp: existence}
\mbox{}
\begin{enumerate}[label = (\roman*)]
	\item \label{enumi: compact}
	$S_{1}$ is compact.
	\item \label{enumi: continuous}
	For $i = 1, \dots, n$, $f_{i}$ is jointly continuous.
	\item \label{enumi: singleton}
	For $i = 2, \dots, n$,
	the set of optimal solutions of the $i$th level problem with arbitrarily fixed $(x_{1}, \dots, x_{i - 1})$ is a singleton.
	\item \label{enumi: bounded}
	For $i = 2, \dots, n$, the optimal solution of the $i$th level problem remains bounded as $x_{1}$ varies in $S_{1}$.
\end{enumerate}
\end{assumption}

Assumption \ref{asp: existence}-\ref{enumi: singleton} means that
the $i$th level problem with parameter $(x_{1}, \dots, x_{i - 1})$ has a unique optimizer,
that is described as $x_{i}^{\ast}$ in \eqref{eq: multilevel}
though it should be formally written as $x_{i}^{\ast}(x_{1}, \dots, x_{i - 1})$
since it is determined by $x_{1}, \dots, x_{i - 1}$.
We define $F_{i}(x_{1}, \dots, x_{i})
\coloneqq f_{i}(x_{1}, \dots, x_{i}, x_{i + 1}^{\ast}, \dots, x_{n}^{\ast})$ by eliminating $x_{i + 1}^{\ast}, \dots, x_{n}^{\ast}$
since it depends only on $(x_{1}, \dots, x_{i})$.
Then, the $i$th level problem with parameter $(x_{1}, \dots, x_{i - 1})$ can be written as
$\min_{x_{i}} F_{i}(x_{1}, \dots, x_{i})$.
Similarly, in the $1$st level problem, for fixed $x_{1}$,
the remaining variables $x_{2}, x_{3}, \dots, x_{n}$
can be represented as $x_{2}^{\ast}(x_{1})$, $x_{3}^{\ast}(x_{1}, x_{2}^{\ast}(x_{1}))$ and so on.
In what follows, we denote them by $x_{2}^{\ast}(x_{1}), x_{3}^{\ast}(x_{1}), \dots, x_{n}^{\ast}(x_{1})$
since they are determined by $x_{1}$.
Eliminating them, we can rewrite the $1$st level problem, equivalently \eqref{eq: multilevel}, as
\begin{equation}
\min_{x_{1} \in S_{1}} F_{1}(x_{1}).
\label{eq: 1st level asp}
\end{equation}

\begin{theorem}
\label{thm: existence}
Under Assumption \ref{asp: existence},
Problem \eqref{eq: 1st level asp} admits optimal solutions.
\end{theorem}

Theorem \ref{thm: existence} is a generalization of \cite[Theorem 3.1]{FFSGP18}, which is valid only for $n = 2$, to general $n$.
It is difficult to extend the original proof to general $n$
because, to derive the continuity of $x_{i}^{\ast}$ extending the original proof,
a sequence $(x_{1}, \dots, x_{i - 1})$ approaches to its accumulation point only from a specific direction.
Instead, in the multilevel case, we employ the theory of point-to-set mapping.
A proof of Theorem \ref{thm: existence} is shown in Supplementary material \ref{sect: existence proof}.

\subsection{Approximation by iterative methods}
We extend the gradient method for bilevel optimization problems proposed by Franceschi et al.\ \cite{FDFP17}
to multilevel optimization problems \eqref{eq: multilevel}.
Based on an argument similar to one in Subsection \ref{sect: bilevel},
we approximate lower level problems in \eqref{eq: multilevel}
by applying an iterative method with $T_{i}$ iterations to the $i$th level problem for $i = 2, \dots, n$.
Then, we obtain the following approximated problem:
\begin{equation}
\begin{alignedat}{2}
\min_{x_{1} \in S_{1}, \{x_{2}^{(t_{2})}\}, \dots, \{x_{n}^{(t_{n})}\}}{}&
	f_{1}(x_{1}, x_{2}^{(T_{2})}, \dots, x_{n}^{(T_{n})})\\
\st{}&
	x_{i}^{(t_{i})} = \Phi_{i}^{(t_{i})}(x_{1}, x_{2}^{(T_{2})}, \dots, x_{i - 1}^{(T_{i - 1})}, x_{i}^{(t_{i} - 1)})&\
	&(i = 2, \dots, n; t_{i} = 1, \dots, T_{i}),
\end{alignedat}
\label{eq: approx}
\end{equation}
where $\Phi_{i}^{(t_{i})}(x_{1}, x_{2}^{(T_{2})}, \dots, x_{i - 1}^{(T_{i - 1})}, x_{i}^{(t_{i})})$ is
the $t_{i}$th iteration formula of an iterative method for the $i$th level approximated problem
with parameters $(x_{1}, \dots, x_{i - 1}) = (x_{1}, x_{2}^{(T_{2})}, \dots, x_{i - 1}^{(T_{i - 1})})$, i.e.,
\begin{equation}
\begin{alignedat}{2}
\min_{\substack{x_{i}, \{x_{i + 1}^{(t_{i + 1})}\}, \dots, \{x_{n}^{(t_{n})}\}}}{}&
	f_{i}(x_{1}, \dots, x_{i}, x_{i + 1}^{(T_{i + 1})}, \dots, x_{n}^{(T_{n})})\\
\st{}&
	x_{j}^{(t_{j})}
	= \Phi_{j}^{(t_{j})}(x_{1}, \dots, x_{i}, x_{i + 1}^{(T_{i + 1})}, \dots, x_{j - 1}^{(T_{j - 1})}, x_{j}^{(t_{j} - 1)})&\
	&	(j = i + 1, \dots, n; t_{j} = 1, \dots, T_{j}),
\end{alignedat}
\end{equation}
and its initial point $x_{i}^{(0)}$ for each $i$ are regarded as a given constant.
For $i = 1, \dots, n$, if we fix $x_{1}, \dots, x_{i}$,
each remaining variable $x_{j}^{(t)}$ is determined by constraints
and thus, we denote it by $x_{j}^{(t_{j})}(x_{1}, \dots, x_{i})$ for $j = i + 1, \dots, n$.
Using this notation, the $i$th level objective function can be written as
\begin{multline}
\tilde{F}_{i}(x_{1}, \dots, x_{i})\\
\coloneqq
f_{i}(x_{1}, \dots, x_{i - 1}, x_{i}, x_{i + 1}^{(T_{i + 1})}(x_{1}, \dots, x_{i}), \dots,
x_{n}^{(T_{n})}(x_{1}, \dots, x_{n - 1}^{(T_{n - 1})}(\cdots(x_{i + 1}^{(T_{i + 1})}(x_{1}, \dots, x_{i}))\cdots))),
\end{multline}
where $(x_{1}, \dots, x_{i - 1})$ acts as a fixed parameter in the $i$th level problem.
When we approximate the $i$th level problem with the steepest descent method for example, we use
$\Phi_{i}^{(t_{i})}(x_{1}, \dots, x_{i - 1}, x_{i}^{(t_{i})})
= x_{i}^{(t_{i} - 1)} - \alpha_{i}^{(t_{i} - 1)} \nabla_{x_{i}} \tilde{F}_{i}(x_{1}, \dots, x_{i - 1}, x_{i}^{(t_{i} - 1)})$.

Especially, all variables other than $x_{1}$ in \eqref{eq: approx} can be expressed as
$\{x_{2}^{(t_{2})}(x_{1})\}, \dots, \{x_{n}^{(t_{n})}(x_{1})\}$
since it is determined by $x_{1}$.
By plugging the constraints into the objective function and eliminating them,
we can reformulate \eqref{eq: approx} as
\begin{equation}
\min_{x_{1} \in S_{1}} \tilde{F}_{1}(x_{1}).
\label{eq: 1st level approx}
\end{equation}

In the remainder of this subsection, we assume $T_{2} = \dots = T_{n} = T$ for simplicity.
Then, the optimal value and solutions of Problem \eqref{eq: 1st level approx} converge to
those of Problem \eqref{eq: 1st level asp} as $T \to \infty$ in some sense.
To derive it, we introduce the following assumption,
which is also a natural extension of that in \cite{FFSGP18}.

\begin{assumption}
\label{asp: convergence}
\mbox{}
\begin{enumerate}[label = (\roman*)]
	\item \label{enumi: upper unif conv}
	$f_{1}(x_{1}, \cdot, \dots, \cdot)$ is uniformly Lipschitz continuous on $S_{1}$.
	\item \label{enumi: lower unif conv}
	For all $i = 2, \dots, n$, sequence $\{x_{i}^{(T)}(x_{1})\}$ converges uniformly to $x_{i}^{\ast}(x_{1})$ on $S_{1}$ as $T \to \infty$.
\end{enumerate}
\end{assumption}

\begin{theorem}
\label{thm: convergence}
Under Assumptions \ref{asp: existence} and \ref{asp: convergence}, the followings hold:
\begin{enumerate}[label = (\alph*)]
	\item The optimal value of Problem \eqref{eq: 1st level approx} converges to that of Problem \eqref{eq: 1st level asp}
	as $T \to \infty$.
	\item The set of the optimal solutions of Problem \eqref{eq: 1st level approx} converges to that of Problem \eqref{eq: 1st level asp};
	more precisely, denoting an optimal solution of Problem \eqref{eq: 1st level approx}
	by $x_{1, T}^{\ast}$,
	we have
	\begin{itemize}
		\item $\{x_{1, T}^{\ast}\}_{T = 1}^{\infty}$ admits a convergent subsequence;
		\item for every subsequence $\{x_{1, T_{k}}^{\ast}\}_{k = 1}^{\infty}$ such that
		$x_{1, T_{k}}^{\ast} \to x_{1}^{\ast}$ as $k \to \infty$,
		the accumulation point $x_{1}^{\ast}$ is an optimal solution of Problem \eqref{eq: 1st level asp}.
	\end{itemize}
\end{enumerate}
\end{theorem}

See Supplementary material \ref{sect: convergence proof} for a proof of this theorem.

\section{Proposed method: Gradient computation in the approximated problem}
\label{sect: algo}

Now we propose to apply a projected gradient method to the approximated problem \eqref{eq: 1st level approx}
for multilevel optimization problems \eqref{eq: multilevel}
because \eqref{eq: 1st level approx} asymptotically converges to \eqref{eq: multilevel} as shown in Theorem \ref{thm: convergence}.
In this section, we derive the formula of $\nabla_{x_{1}} \tilde{F}_{1}(x_{1})$
and confirm the global convergence of the resulting projected gradient method.

\subsection{Gradient of the objective function in the approximated problem}

The following theorem provides a computation formula of $\nabla_{x_{1}} \tilde{F}_{1}(x_{1})$.

\begin{theorem}[Gradient formula for the $n$-level optimization problems]
\label{thm: grad formula}
The gradient $\nabla_{x_{1}} \tilde{F}_{1}(x_{1})$ can be expressed as follows:
\begin{align}
\nabla_{x_{1}} \tilde{F}_{1}(x_{1})
&= \nabla_{x_{1}} f_{1}(x_{1}, x_{2}^{(T_{2})}, \dots, x_{n}^{(T_{n})})
+ \sum_{i = 2}^{n} Z_{i} \nabla_{x_{i}} f_{1}(x_{1}, x_{2}^{(T_{2})}, \dots, x_{n}^{(T_{n})}),\\
Z_{i}
&= \sum_{t = 1}^{T_{i}} \left(\sum_{j = 2}^{i - 1} Z_{j} C_{i j}^{(t)} + B_{i}^{(t)}\right) \prod_{s = t + 1}^{T_{i}} A_{i}^{(s)},\\
A_{i}^{(t)}
&= \nabla_{x_{i}} \Phi_{i}^{(t)}(x_{1}, x_{2}^{(T_{2})}, \dots, x_{i - 1}^{(T_{i - 1})}, x_{i}^{(t - 1)}),\\
B_{i}^{(t)}
&= \nabla_{x_{1}} \Phi_{i}^{(t)}(x_{1}, x_{2}^{(T_{2})}, \dots, x_{i - 1}^{(T_{i - 1})}, x_{i}^{(t - 1)}),\\
C_{i j}^{(t)}
&= \nabla_{x_{j}} \Phi_{i}^{(t)}(x_{1}, x_{2}^{(T_{2})}, \dots, x_{i - 1}^{(T_{i - 1})}, x_{i}^{(t - 1)})
\end{align}
for any $i = 2, \dots, n$; $t = 1, \dots, T_{i}$; and $j = 2, \dots, i - 1$,
where we define $\prod_{s = t + 1}^{T_{i}} A_{i}^{(s)} \coloneqq A_{i}^{(t + 1)} A_{i}^{(t + 2)} \dots A_{i}^{(T_{i})}$ for $t < T_{i}$
and $\prod_{s = T_{i} + 1}^{T_{i}} A_{i}^{(s)} = I$.
\end{theorem}

See supplementary material \ref{sect: grad formula proof} for a proof of this theorem.


We consider computing $\nabla_{x_{1}} \tilde{F}_{1}(x_{1})$ using Theorem \ref{thm: grad formula}.
Notice that we can easily compute $Z_{2} = \sum_{t = 1}^{T_{2}} B_{2}^{(t)} \prod_{s = t + 1}^{T_{2}} A_{2}^{(s)}$.
For $i = 3, \dots, n$, when we have $Z_{2}, \dots, Z_{i - 1}$, we can compute $Z_{i}$.
We show an algorithm
that computes $\nabla_{x_{1}} \tilde{F}_{1}(x_{1})$ by computing $Z_{2}, \dots, Z_{n}$ in this order in Algorithm \ref{algo: grad comp}.

\begin{algorithm}
\caption{Computation of $\nabla_{x_{1}} \tilde{F}_{1}(x_{1})$}
\label{algo: grad comp}
\begin{algorithmic}[1]
\Require{$x_{1}$: current value of the $1$st level variable.
	$\{x_{i}^{(0)}\}_{i = 2}^{n}$: initial values of the lower level iteration.}
\Ensure{The exact value of $\nabla_{x_{1}} \tilde{F}_{1}(x_{1})$.}
\State{$g \coloneqq (0, \dots, 0)^{\top}$.}
\For{$i \coloneqq 2, \dots, n$}
	\State{$Z_{i} \coloneqq O$.}
	\For{$t \coloneqq 1, \dots, T_{i}$}
		\State{$x_{i}^{(t)} \coloneqq \Phi_{i}^{(t)}(x_{1}, x_{2}^{(T_{2})}, \dots, x_{i - 1}^{(T_{i - 1})}, x_{i}^{(t - 1)})$.}
		\State{$\bar{B}_{i}^{(t)} \coloneqq \sum_{l = 2}^{i - 1} Z_{l} C_{i l}^{(t)} + B_{i}^{(t)}$.}
		\State{$Z_{i} \coloneqq Z_{i} A_{i}^{(t)} + \bar{B}_{i}^{(t)}$.}
	\EndFor
\EndFor
\For{$i = 2, \dots, n$}
	\State{$g \coloneqq g + Z_{i} \nabla_{x_{i}} f_{1}$.}
\EndFor
\State{$g \coloneqq g + \nabla_{x_{1}} f_{1}$.}
\State{\textbf{return} $g$}
\end{algorithmic}
\end{algorithm}

For $i = 2, \dots, n$ and $t = 1, \dots, T_{i}$,
$\Phi_{i}^{(t)}$, which appears in the 5th line of Algorithm \ref{algo: grad comp},
is the update formula based on the gradient $\nabla_{x_{i}} \tilde{F}_{i}(x_{1}, \dots, x_{i})$ of the $i$th level objective function.
$\nabla_{x_{i}} \tilde{F}_{i}(x_{1}, \dots, x_{i})$  can be computed by 
applying Algorithm \ref{algo: grad comp} to the $(n - i + 1)$-level optimization problem
with objective functions $\tilde{F}_{i}, \dots, \tilde{F}_{n}$.
Therefore, recursively calling Algorithm \ref{algo: grad comp} in the computation of $\Phi_{i}^{(t)}$,
we can compute $\nabla_{x_{1}} \tilde{F}_{1}(x_{1})$.
For an example of applying Algorithm \ref{algo: grad comp} to Problem \eqref{eq: 1st level approx}
arising from a trilevel optimization problem, i.e., Problem \eqref{eq: multilevel} with $n = 3$,
see Supplementary material \ref{sect: trilevel grad comp}.


\subsection{Complexity of the gradient computation}

We analyze the complexity for computing $\nabla_{x_{1}} \tilde{F}_{1}(x_{1})$ 
by recursively calling Algorithm \ref{algo: grad comp}.
In the following theorem, the asymptotic big O notation is denoted by $\bigO(\cdot)$.

\begin{theorem}
\label{thm: complexity recursive}
Let the time and space complexity for computing $\nabla_{x_{i}} \tilde{F}_{i}(x_{i})$ be $c_{i}$ and $s_{i}$, respectively.
We use $\Phi_{i}^{(t_{i})}$ based on $\nabla_{x_{i}} \tilde{F}_{i}(x_{i})$
and recursively call Algorithm \ref{algo: grad comp} for computing $\nabla_{x_{i}} \tilde{F}_{i}(x_{1}, \dots, x_{i})$.
In addition, we assume the followings:
\begin{itemize}
	\item The time and space complexity for evaluating $\Phi_{i}^{(t_{i})}$ are $\bigO(c_{i})$ and $\bigO(s_{i})$.
	\item The time and space complexity of $\nabla_{x_{1}} f_{1}$ and $\nabla_{x_{i}} f_{1}$ for $i = 1, \dots, n$
	are smaller in the sense of the order than those of \texttt{for} loops in lines 2--9 in Algorithm \ref{algo: grad comp}.
\end{itemize}
Then, the overall time complexity $c_{1}$ and space complexity $s_{1}$
for computing $\nabla_{x_{i}} \tilde{F}_{i}(x_{1}, \dots, x_{i})$ can be written as
\begin{alignat}{2}
c_{1} &= \bigO\left(p^{n} n! c_{n} \prod_{i = 1}^{n - 1} (T_{i + 1} d_{i})\right),&\quad
s_{1} &= \bigO(q^{n} s_{n}),
\label{eq: complexity recursive}
\end{alignat}
respectively, for some constant $p, q > 1$.
\end{theorem}

For a proof, see Supplementary material \ref{sec:complexity}.
Note that, if $n$ is a fixed parameter, 
those complexity reduces to a polynomial of $T_{i}$'s, $d_{i}$'s, $c_{n}$, and $s_{n}$.
Hence, Algorithm \ref{algo: grad comp} can be regarded as a fixed-parameter tractable algorithm.

\subsection{Global convergence of the projected gradient method}

Here, we consider solving Problem \eqref{eq: 1st level approx} by the projected gradient method,
which calculates the gradient vector by Algorithm \ref{algo: grad comp}
and projects the updated point 
on $S_{1}$ in each iteration.
When all lower-level updates are based on the steepest descent method,
we can derive the Lipschitz continuity of the gradient of the objective function of \eqref{eq: 1st level approx}.
Hence, we can guarantee the global convergence of the projected gradient method for \eqref{eq: 1st level approx}
by taking a sufficiently small step size.

\begin{theorem}
\label{thm: Lip cont grad}
Suppose $\Phi_{i}^{(t)}(x_{1}, \dots, x_{i - 1}, x_{i}^{(t - 1)})
= x_{i}^{(t - 1)} - \alpha_{i}^{(t - 1)} \nabla_{x_{i}} \tilde{F}_{i} (x_{1}, \dots, x_{i - 1}, x_{i}^{(t - 1)})$
for all $i = 2, \dots, n$ and $t_{i} = 1, \dots, T_{i}$,
where $\alpha_{i}^{(t - 1)}$ and $x_{i}^{(0)}$ are given parameters for all $i$ and $t$.
Assume that $\nabla_{x_{j}} f_{i}$ is Lipschitz continuous and bounded for all $i = 1, \dots, n$ and $j = 1, \dots, n$ and also
$\nabla_{x_{i}} \Phi_{i}^{(t)}$, $\nabla_{x_{1}} \Phi_{i}^{(t)}$, and $\nabla_{x_{j}} \Phi_{i}^{(t)}$ are Lipschitz continuous and bounded
for all $i = 2, \dots, n$; $j = 2, \dots, i - 1$; $t = 1, \dots, T_{i}$.
Then, $\nabla_{x_{1}} \tilde{F}_{1}$ is Lipschitz continuous.
\end{theorem}

See Supplementary material \ref{sect: Lip cont grad proof} for a proof of this theorem.

\begin{corollary}
Suppose the same assumption as Theorem \ref{thm: Lip cont grad} and assume $S_{1}$ is a compact convex set.
Let $L$ be the Lipschitz constant of $\nabla_{x_{1}} \tilde{F}_{1}$.
Then, a sequence $\{x_{1}^{(t)}\}$ generated by the projected gradient method with sufficiently small constant step size,
e.g., smaller than $2 / L$, for Problem \eqref{eq: 1st level approx} from any initial point
has a convergent subsequence that converges to a stationary point with convergence rate $\bigO(1 / \sqrt{t})$.
\end{corollary}

\begin{proof}
From Theorem \ref{thm: Lip cont grad},
the gradient of the objective function of Problem \eqref{eq: 1st level approx} is $L$-Lipschitz continuous.
Let $G\colon \intmath(\dom(\tilde{F}_{1})) \to \bbR^{d_{1}}$ be the gradient mapping \cite[Definition 10.5]{Bec17}
corresponding to $\tilde{F}_{1}$, the indicator function of $S_{1}$,
and the constant step size $\alpha_{1}^{(t)}$ with satisfying $0 < \alpha_{1}^{(t)} < 2 / L$ for all $t$.
Note that $\|G(x_{1})\| = 0$ if and only if
$x_{1}$ is a stationary point of Problem \eqref{eq: 1st level approx} \cite[Theorem 10.7]{Bec17}.
By applying \cite[Theorem 10.15]{Bec17},
we obtain $\min_{s = 0}^{t} \|G(x_{1}^{(s)})\| \le \bigO(1 / \sqrt{t})$ and $\|G(\bar{x}_{1})\| = 0$,
where $\bar{x}_{1}$ is a limit point of $\{x_{1}^{(t)}\}$.
\end{proof}

%% file: chapters/Experiments.tex
\section{Numerical experiments}
\label{sect: experiments}

To validate the effectiveness of our proposed method,
we conducted numerical experiments on an artificial problem and a hyperparameter optimization problem arising from real data
(see Supplementary material \ref{sect: experiments detail} for complete results).
In our numerical experiments,
we implemented all codes with Python 3.9.2 and JAX 0.2.10 for automatic differentiation
and executed them on a computer with 12 cores of Intel Core i7-7800X CPU 3.50 GHz, 64 GB RAM, Ubuntu OS 20.04.2 LTS.

In this section, we used Algorithm 1 to calculate the gradient of $\tilde{F}_{i}$ in problem \ref{eq: approx},
and used automatic differentiation to calculate the gradient of $\Phi_{i}^{(t)}$ in problem \ref{eq: approx}.

\subsection{Convergence to the optimal solution}
\label{exp: convergence}

We solved the following trilevel optimization problem with Algorithm \ref{algo: grad comp} to evaluate the performance:
\begin{equation}
\begin{aligned}
\min_{x_{1} \in \bbR^{2}}{}&
	f_{1}(x_{1}, x_{2}^{\ast}, x_{3}^{\ast}) = \|x_{3}^{\ast} - x_{1}\|_{2}^{2} + \|x_{1}\|_{2}^{2} \st\\
 &	\begin{aligned}
	x_{2}^{\ast} \in \argmin_{x_{2} \in \bbR^{2}}{}&
		f_{2}(x_{1}, x_{2}, x_{3}^{\ast}) = \|x_{2} - x_{1}\|_{2}^{2} \st\\
&	x_{3}^{\ast} \in \argmin_{x_{3} \in \bbR^{2}} f_{3}(x_{1}, x_{2}, x_{3}) = \|x_{3} - x_{2}\|_{2}^{2}.
	\end{aligned}
\end{aligned}
\label{eq: artificial instance}
\end{equation}
Clearly, the optimal solution for this problem is $x_{1} = x_{2} = x_{3} = (0, 0)^{\top}$.

\begin{figure}
\centering
\subfigure[$T_{2} = 10$, $T_{3} = 10$]{
\includegraphics[width = 0.49\textwidth]{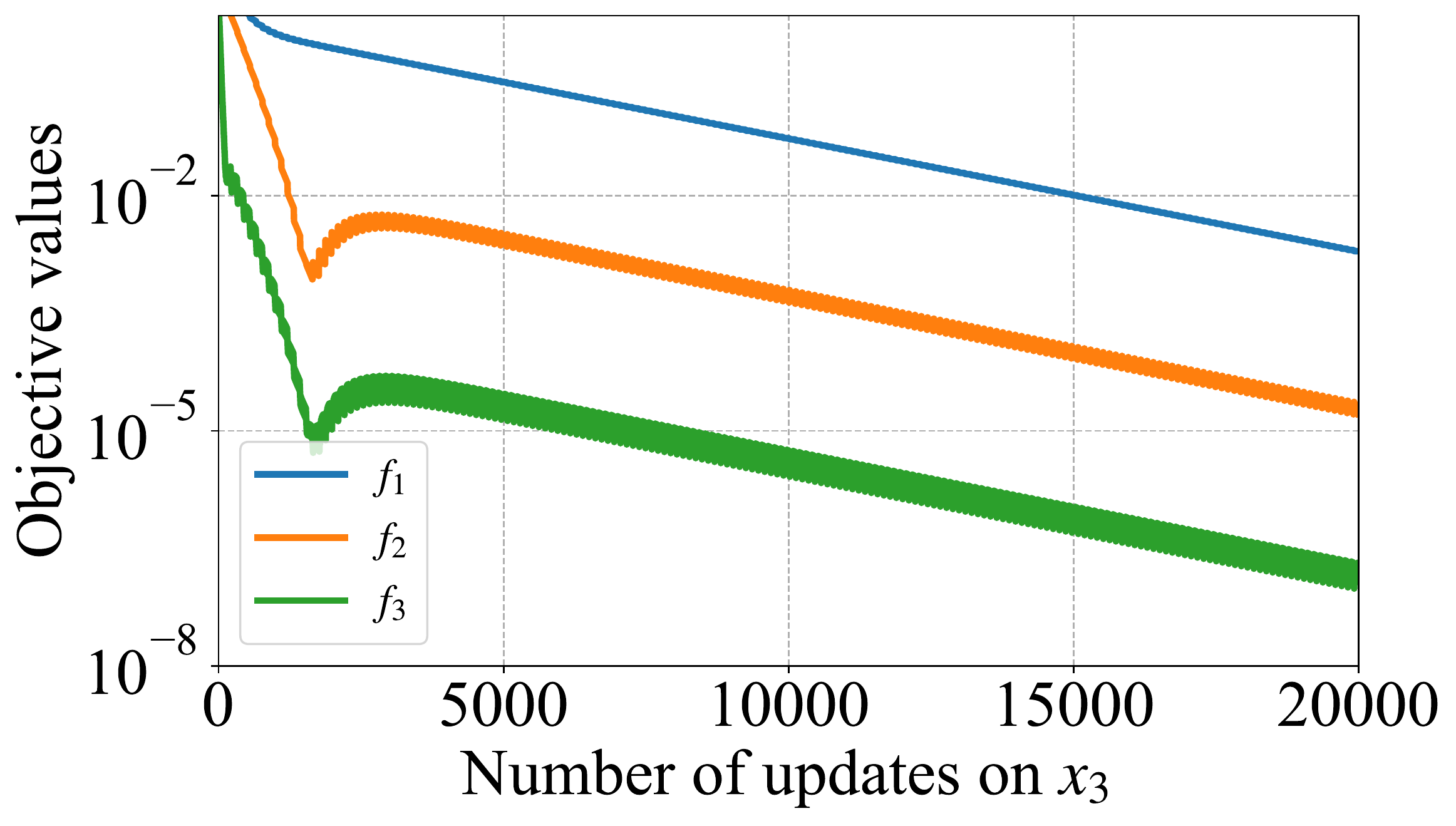}}
\subfigure[$T_{2} = 1$, $T_{3} = 1$]{
\includegraphics[width = 0.49\textwidth]{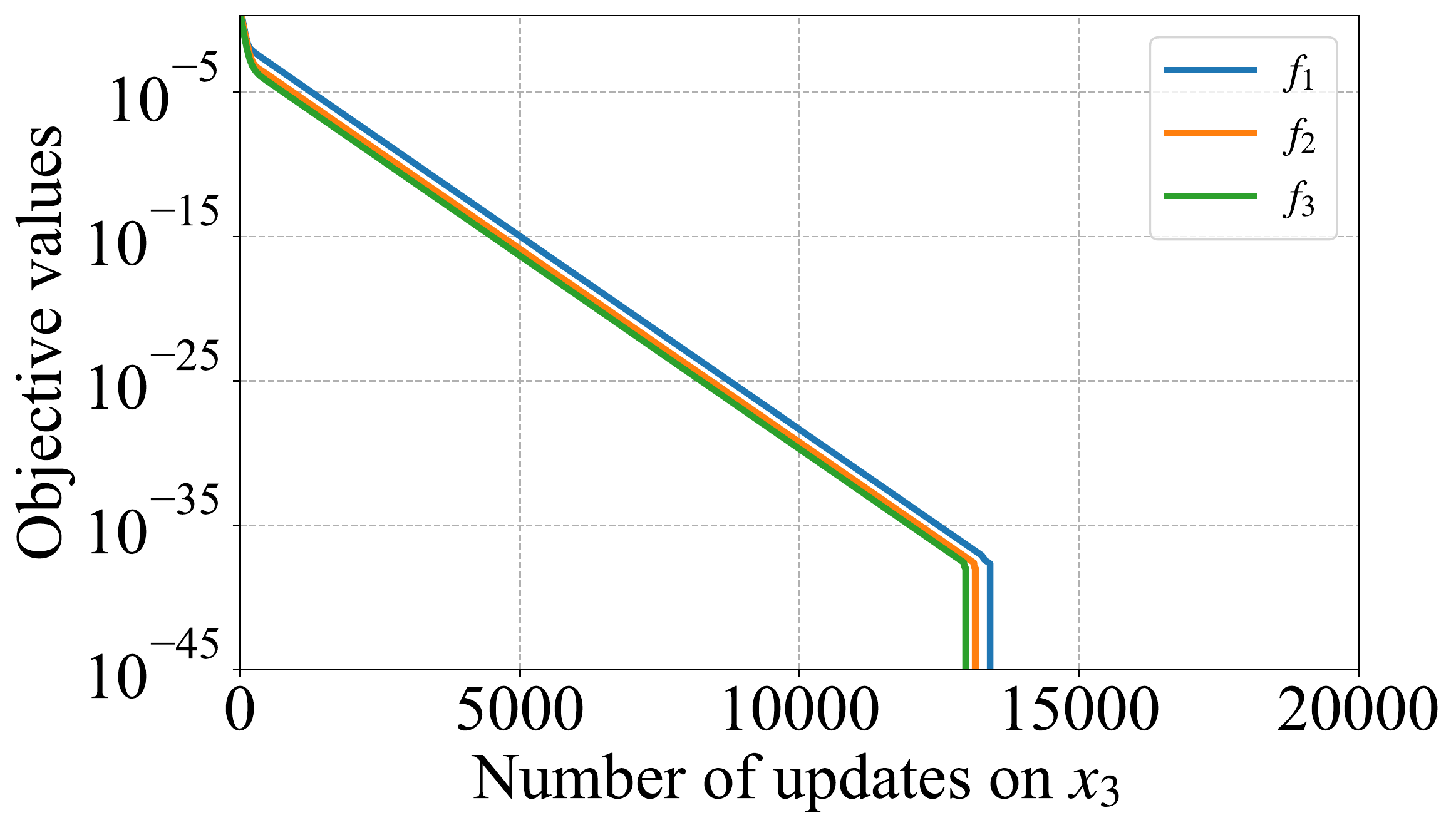}}
\caption{Performance of Algorithm \ref{algo: grad comp} for Problem \eqref{eq: artificial instance}.
	The objective values linearly decreased.}
\label{fig: objective values}
\centering
\subfigure[$T_{2} = 10$, $T_{3} = 10$]{%
\includegraphics[width = 0.49\textwidth]{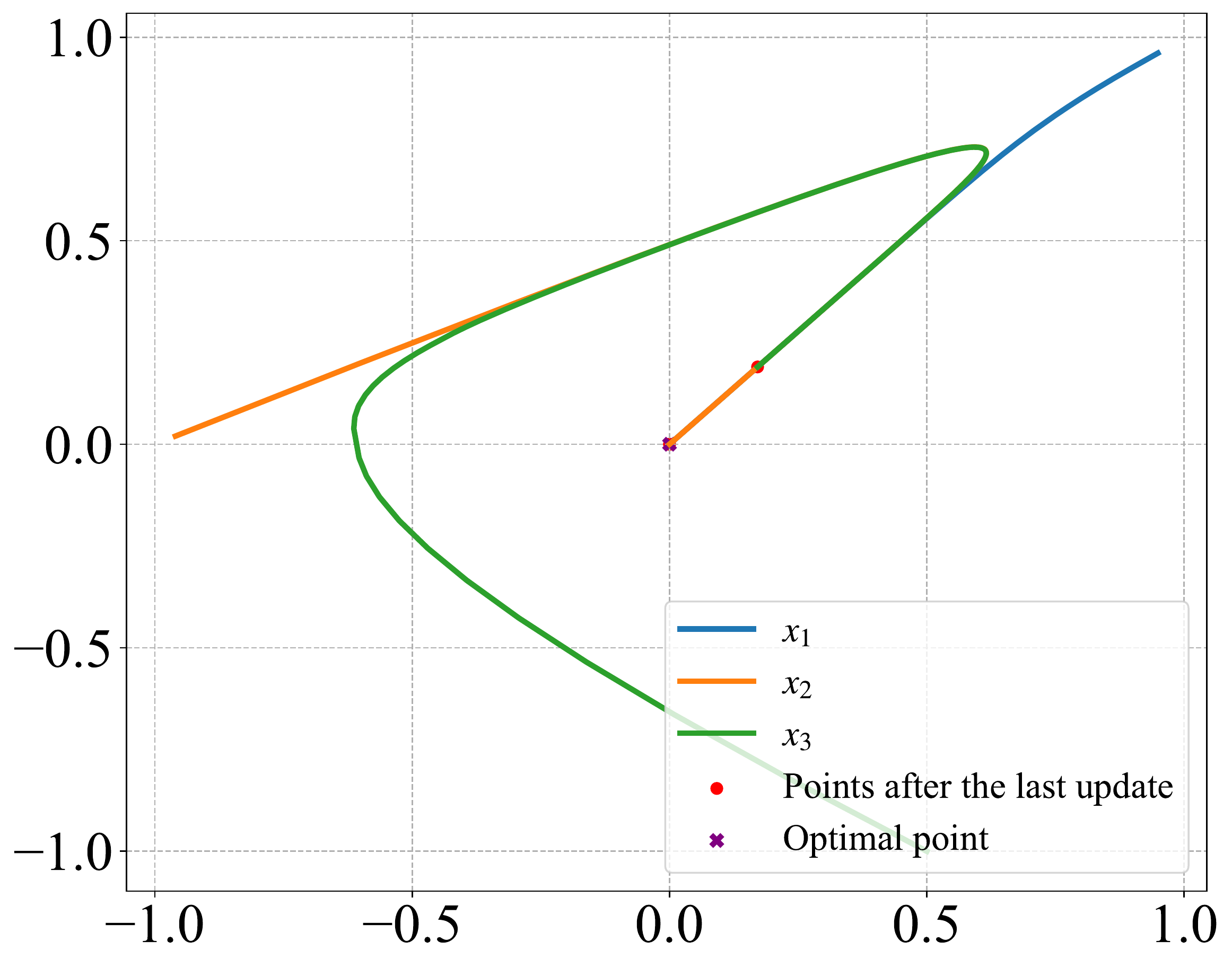}}
%
%
\subfigure[$T_{2} = 1$, $T_{3} = 1$]{%
\includegraphics[width = 0.49\textwidth]{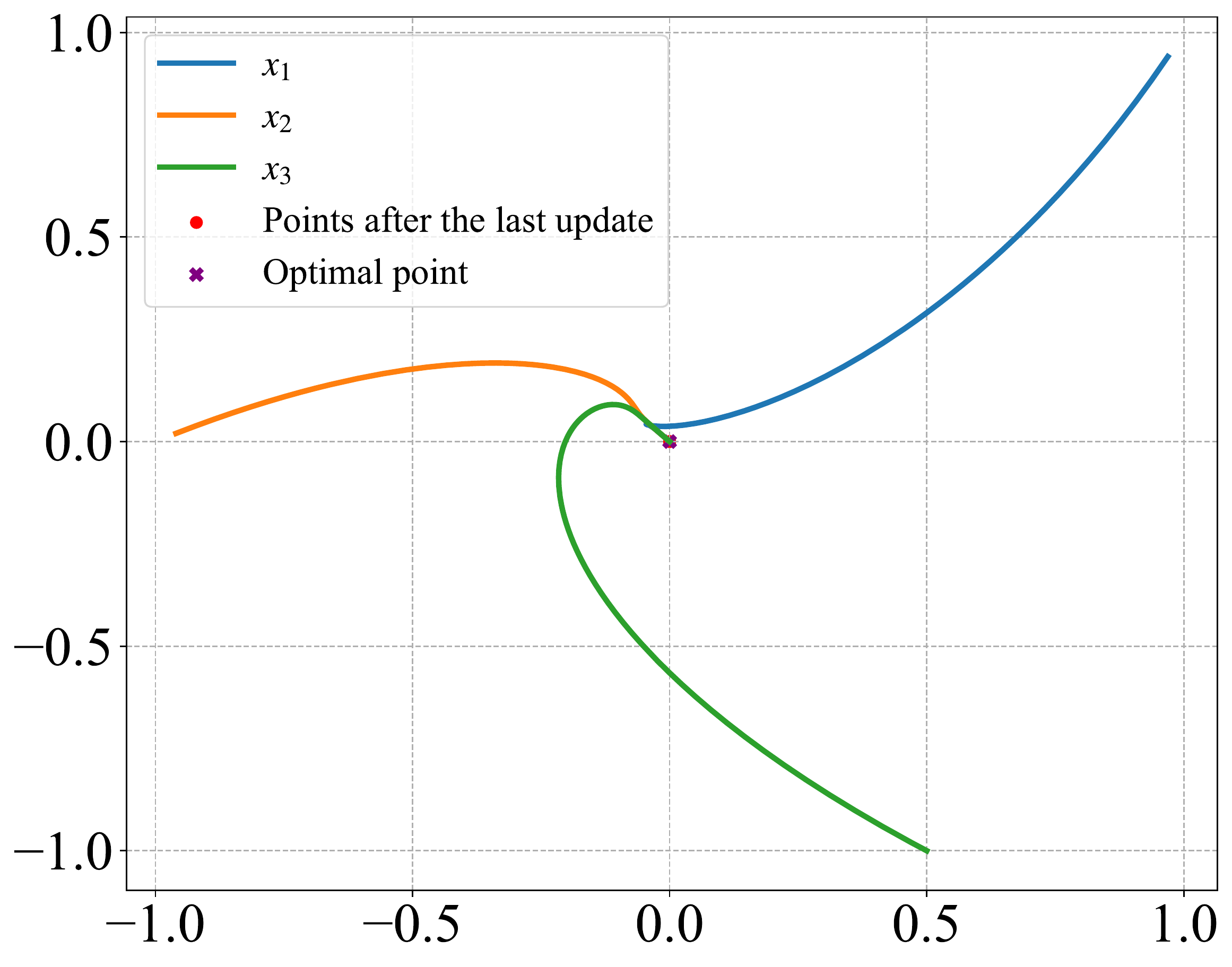}}
\caption{Trajectories of each variable.
	$(0, 0)^{\top}$ is the optimal solution of each level.
	Our algorithm with few iterations for lower-level problems $T_{2} = T_{3} = 1$ performed well.}
\label{fig: trajectory}
\end{figure}

We solved \eqref{eq: artificial instance} with fixed constant step size and initialization but different $(T_{2}, T_{3})$.
For the iterative method in Algorithm \ref{algo: grad comp},
we employed the steepest descent method at all levels.
We show the transition of the value of the objective functions in
Figure \ref{fig: objective values}
and the trajectories of each decision variable in
Figure \ref{fig: trajectory}.
Since updates of $x_{3}$ is the most inner iteration,
the time required to update $x_{3}$ does not change when $T_{2}$ or $T_{3}$ changes.
Hence the number of updates of $x_{3}$ is proportional to the total computational time.
Therefore, we can compare the time efficiency of the optimization algorithm by focusing on the number of updates of $x_{3}$.
We confirmed that the values of $f_{1}$, $f_{2}$, and $f_{3}$ converge to the optimal value $0$ at all levels and
the gradient method outputs the optimal solution $(0, 0)^{\top}$
even if we use small iteration numbers $T_{2} = 1$ and $T_{3} = 1$ for the approximation problem \eqref{eq: approx}.

We also made comparison Algorithm \ref{algo: grad comp} and an existing algorithm \cite{TKMO12} based on evolutionary strategy
by solving \eqref{eq: artificial instance} using both algorithms.
Algorithm \ref{algo: grad comp} outperformed the existing algorithm.
For detail, see Supplementary material \ref{sect: evolutionary}.


\subsection{Application to hyperparameter optimization}\label{exp: HO}
For deriving a machine learning model robust to noise in input data,
we formulate a trilevel model by assuming two players: a model learner and an attacker.
The model learner decides the hyperparameter $\lambda$ to minimize the validation error,
while the attacker tries to poison training data so as to make the model less accurate.
This model is inspired by bilevel hyperparameter optimization \cite{FFSGP18} and
adversarial learning \cite{JOBLNL18, LLCFXAHO20} and formulated as follows:
\begin{equation}
\begin{aligned}
\min_{\lambda}{}&
	\frac{1}{m} \|y_{\text{valid}} - X_{\text{valid}} \theta\|_{2}^{2} \st\\
&	\begin{aligned}
	P \in \argmax_{P'}{}&
		\frac{1}{n} \|y_{\text{train}} - (X_{\text{train}} + P') \theta\|_{2}^{2} - \frac{c}{n d}\|P'\|_{2}^{2} \st\\
	&	\theta \in \argmin_{\theta'} \frac{1}{n} \|y_{\text{train}} - (X_{\text{train}} + P') \theta'\|_{2}^{2}
		+ \exp(\lambda) \frac{\|\theta'\|_{1^{\ast}}}{d},
	\end{aligned}
\end{aligned}
\end{equation}
where $\theta$ denotes the parameter of the model,
$d$ denotes the dimension of $\theta$,
$n$ denotes the number of the training data $X_{\text{train}}$,
$m$ denotes the number of the validation data $X_{\text{val}}$,
$c$ denotes the penalty for the noise $P$,
and $\|\cdot\|_{1^{\ast}}$ is a smoothed $\ell_{1}$-norm \cite[Eq.~(18) with $\mu = 0.25$]{SNC19},
which is a differentiable approximation of the $\ell_{1}$-norm.
Here, we use $\exp(\lambda)$ to express a nonnegative penalty parameter
instead of the constraint $\lambda \ge 0$.

To validate the effectiveness of our proposed method,
we compared the results by the trilevel model with those of the following bilevel model:
\begin{equation}
\min_{\lambda} \frac{1}{m} \|y_{\text{valid}} - X_{\text{valid}} \theta\|_{2}^{2}
\st \theta \in \argmin_{\theta'}
\frac{1}{n} \|y_{\text{train}} - X_{\text{train}} \theta'\|_{2}^{2} + \exp(\lambda) \frac{\|\theta'\|_{1^{\ast}}}{d}.
\label{ho-2level}
\end{equation}
This model is equivalent to the trilevel model without the attacker's level.

We used Algorithm \ref{algo: grad comp} to compute the gradient of the objective function in
the trilevel and bilevel models with real datasets.
For the iterative method in Algorithm \ref{algo: grad comp},
we employed the steepest descent method at all levels.
We set $T_{2} = 30$ and $T_{3} = 3$ for the trilevel model and $T_{2} = 30$ for the bilevel model.
In each dataset, we used the same initialization and step sizes in the updates of $\lambda$ and $\theta$ in trilevel and bilevel models.
We compared these methods on the regression tasks with the following datasets:
the diabetes dataset \cite{data-diabetes},
the (red and white) wine quality datasets \cite{data-wine},
the Boston dataset \cite{data-Boston}.
For each dataset, we standardized each feature and the objective variable;
 randomly chose 40 rows for training data $(X_{\text{train}}, y_{\text{train}})$,
 chose other 100 rows for validation data $(X_{\text{valid}}, y_{\text{valid}})$,
and used the rest of the rows for test data.

\begin{figure}
	\centering
	\subfigure[$\sigma = 0.01$]{
	\includegraphics[width = 0.49\textwidth]{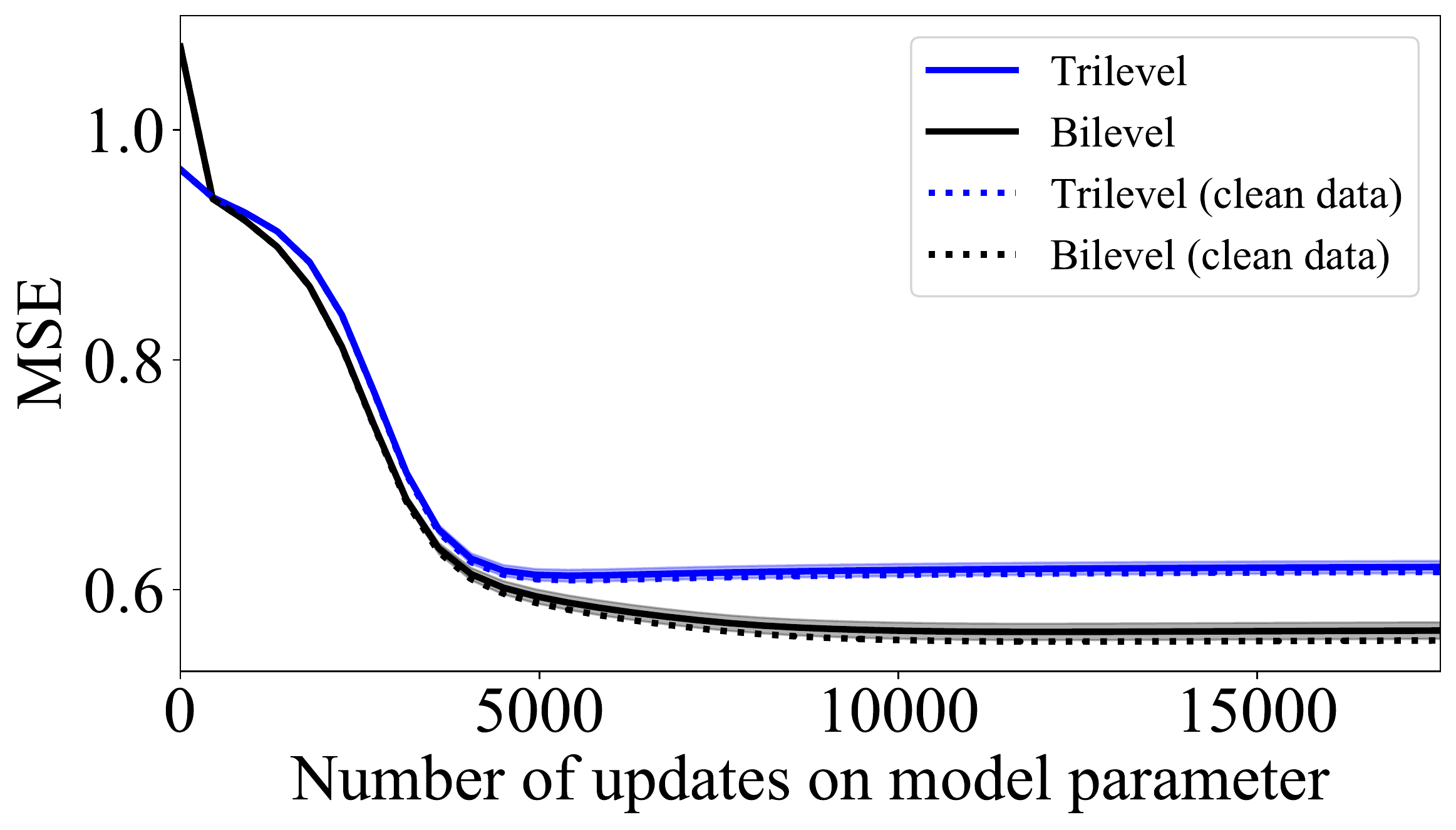}}
	\subfigure[$\sigma = 0.08$]{
	\includegraphics[width = 0.49\textwidth]{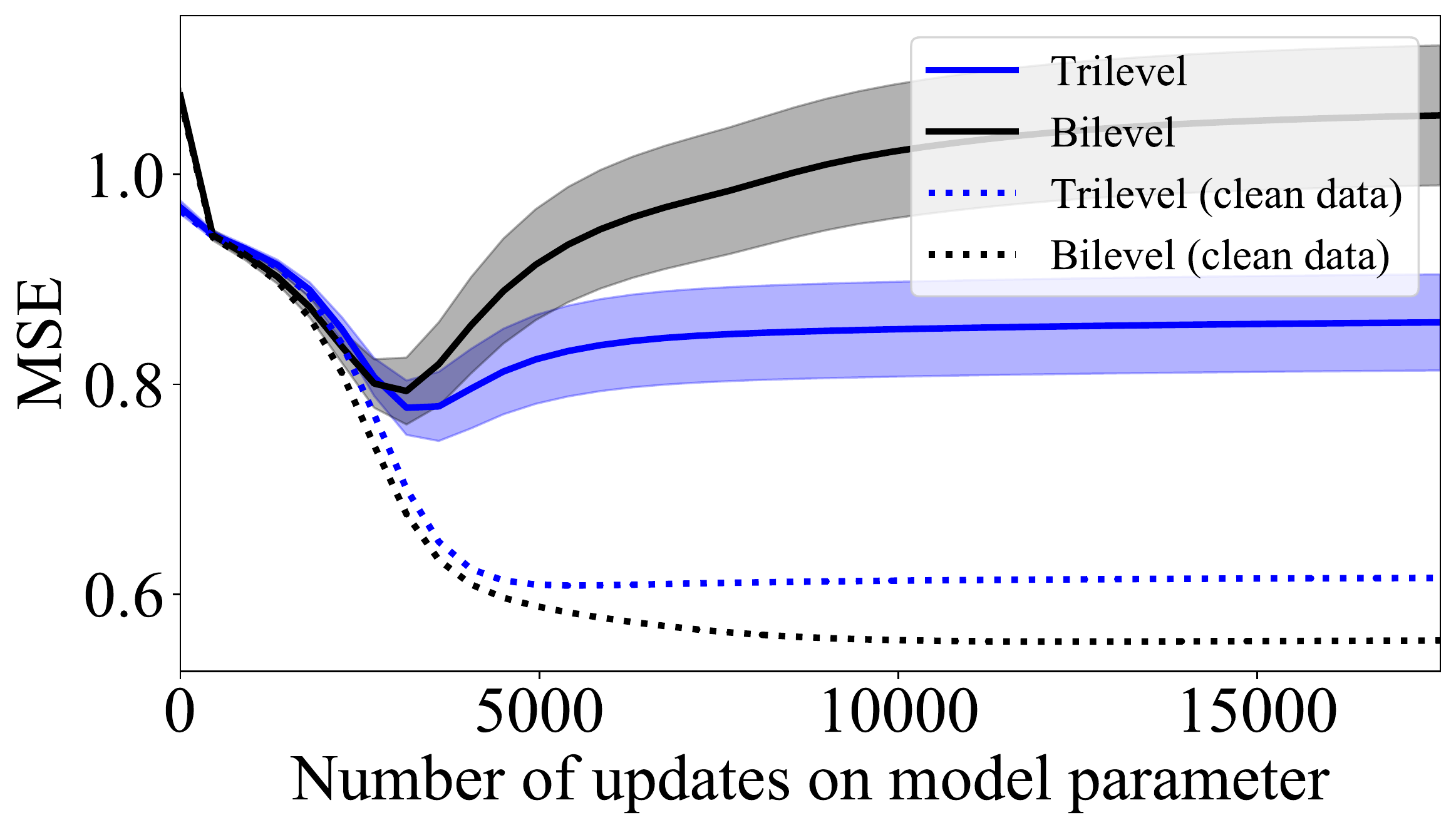}}
	\caption{MSE of test data with Gaussian noise with the standard deviation of $\sigma$ on diabetes dataset.}
	\label{fig: iter-diabetes}
	\centering
	\subfigure[Diabetes dataset]{
	\includegraphics[width = 0.49\textwidth]{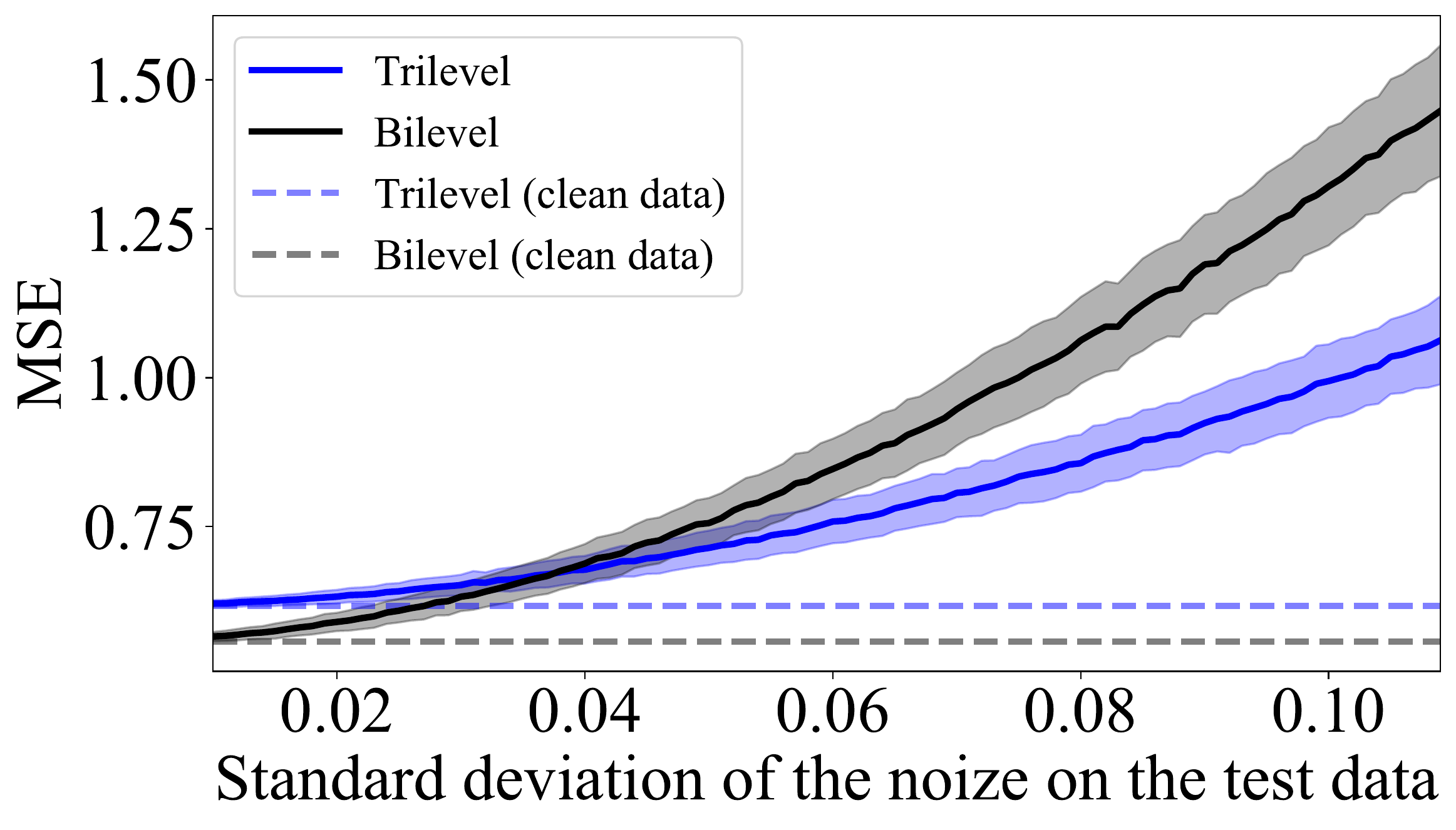}}
	\subfigure[Red wine quality dataset]{
	\includegraphics[width = 0.49\textwidth]{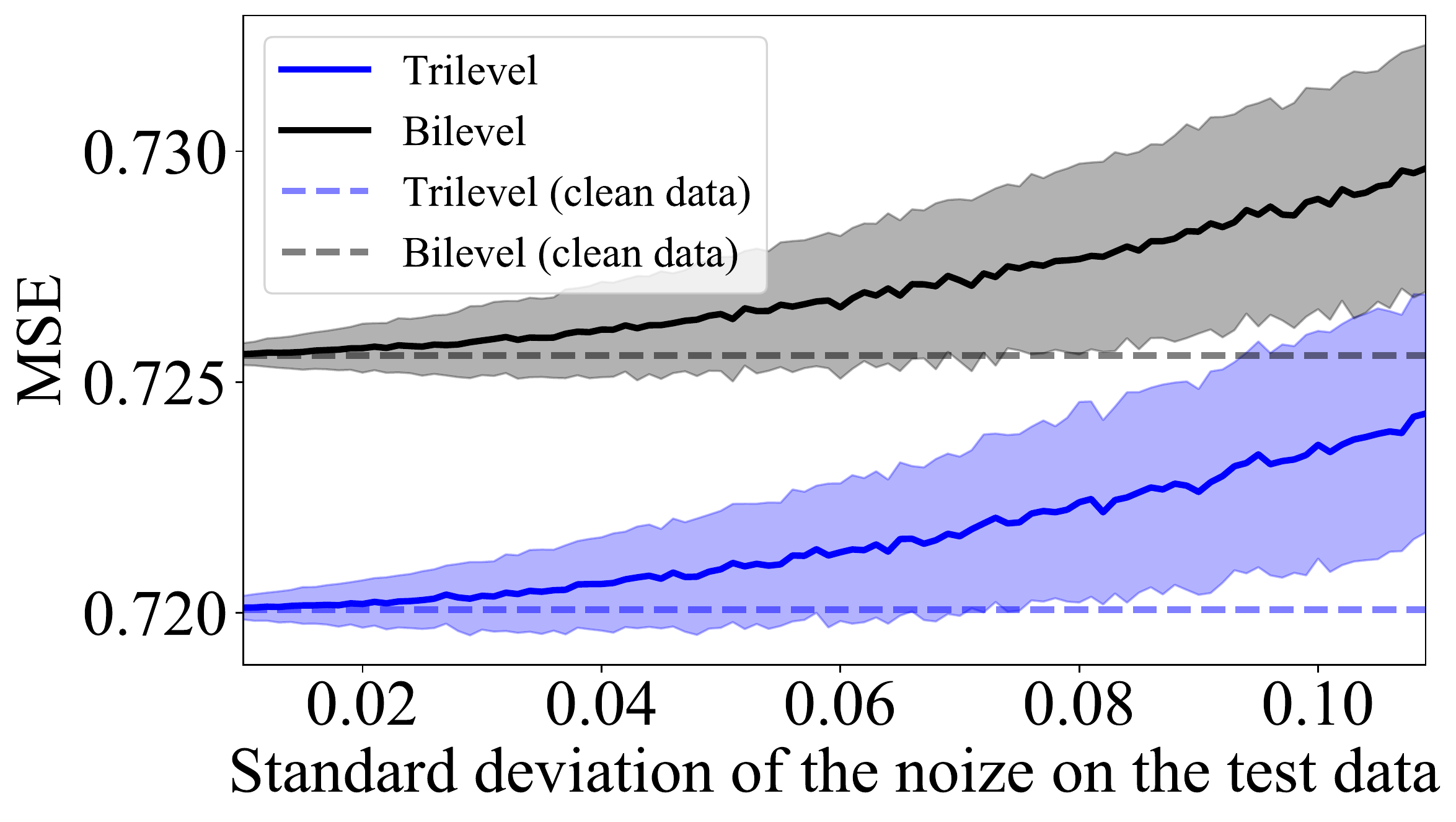}}
	\caption{MSE of test data with Gaussian noise, using early-stopped parameters for prediction.}
	\label{fig: MSE early stop}
\end{figure}

We show the transition of the mean squared error (MSE) by test data with Gaussian noise
in Figure \ref{fig: iter-diabetes}.
The solid line and colored belt respectively indicate the mean and the standard deviation
over 500 times of generation of Gaussian noise.
The dashed line indicates the MSE without noise as a baseline.
In the results of the diabetes dataset (Figure \ref{fig: iter-diabetes}),
the trilevel model provided a more robust parameter than the bilevel model,
because the MSE of the trilevel model rises less in the large noise setting for test data.

Next, we compared the quality of the resulting model parameters by the trilevel and bilevel models.
We set an early-stopping condition on learning parameters:
after 1000 times of updates on the model parameter,
if one time of update on hyperparameter $\lambda$ did not improve test error, 
terminate the iteration and return the parameters at that time.
By using the early-stopped parameters, 
we show the relationship of test error and standard deviation of the noise on the test data
in Figure \ref{fig: MSE early stop} and Table \ref{tbl: MSE}.
For the diabetes dataset, the growth of MSE of the trilevel model was slower than that of the bilevel model.
For wine quality and Boston house-prices datasets,
the MSE of the trilevel model was consistently lower than that of the bilevel model.
Therefore, the trilevel model provides more robust parameters than the bilevel model in these settings of problems.

\begin{table}
\centering
\caption{MSE of test data with Gaussian noise with the standard deviation of 0.08, using early-stopped parameters for prediction.
	The better values are shown in boldface.}
\label{tbl: MSE}
\begin{tabular}{rr@{${} \pm {}$}lr@{${} \pm {}$}lr@{${} \pm {}$}lr@{${} \pm {}$}l}
\toprule
&	\multicolumn{2}{c}{diabetes}&			\multicolumn{2}{c}{Boston}&
	\multicolumn{2}{c}{wine (red)}&			\multicolumn{2}{c}{wine (white)}\\
\midrule
\multicolumn{1}{c}{Trilevel}&
	\textbf{0.8601}&	\textbf{0.0479}&	\textbf{0.4333}&	\textbf{0.0032}&
	\textbf{0.7223}&	\textbf{0.0019}&	\textbf{0.8659}&	\textbf{0.0013}\\
\multicolumn{1}{c}{Bilevel}&
	1.0573&				0.0720&				0.4899&				0.0033&
	0.7277&				0.0019&				0.8750&				0.0014\\
\bottomrule
\end{tabular}
\end{table}

\subsection{Relationship between $(T_{2}, T_{3})$ and the convergence speed}
In the first experiment in Section \ref{exp: convergence},
there is no complex relationship between variables at each level,
and therefore, the objective function value at one level is not influenced significantly
when variables at other levels is changed.
In such a problem setting, by setting $T_{2}$ and $T_{3}$ to small values,
we update $x_{1}$ many times and
the generated sequence $\{x_{1}^{(t)}\}$ quickly converges to a stationary point.
On the other hand, for example, if we set $T_{3}$ to a large value,
$x_{3}$ is well optimized for some fixed $x_{1}$ and $x_{2}$,
and hence, our whole algorithm may need more computation time until convergence.
In the second experiment in Section \ref{exp: HO},
the relationship between variables at each level is more complicated than in the first experiment.
Setting $T_{2}$ or $T_{3}$ smaller in such a problem is not necessarily considered to be efficient
because the optimization algorithm proceeds
without fully approximating the optimality condition of $x_{i}$ at the $i$th level.

%% file: chapters/Conclusion.tex
\section{Conclusion}
\label{sect: concl}

\paragraph{Summary}
In this paper, we have provided an approximated formulation for a multilevel optimization problem by iterative methods
and discussed its asymptotical properties of it.
In addition, we have proposed an algorithm for computing the gradient of the objective function of the approximated problem.
Using the gradient information, we can solve the approximated problem by the projected gradient method.
We have also established the global convergence of the projected gradient method.

\paragraph{Limitation and future work}
Our proposed gradient computation algorithm is fixed-parameter tractable
and hence it works efficiently for small $n$.
For large $n$, however, the exact computation of the gradient is expensive.
Development of heuristics for approximately computing the gradient is left for future research.
Weakening assumptions in the theoretical contribution is also left for future work.
In addition, there is a possibility of another algorithm to solve the approximated
problem \ref{eq: approx}.
In this paper, we propose Algorithm \ref{algo: grad comp}, which corresponds to
forward mode automatic differentiation. 
On the other hand, in the prior research for bilevel optimization \cite{FDFP17},
two algorithms were proposed from the perspective of
forward mode automatic differentiation and reverse mode automatic differentiation, respectively.
Therefore, there is a possibility of another algorithm for problem \ref{eq: approx} which corresponds to the reverse mode automatic differentiation,
and that is left for future work.

%% file: chapters/Supplemental.tex
\section*{Supplementary materials}

\section{Proofs}
\label{sect: proofs}

\subsection{Proof of Theorem \ref{thm: existence}}
\label{sect: existence proof}
To prove Theorem \ref{thm: existence}, we introduce some definitions and a lemma.

\begin{definition}[{\cite{Hog73}}]
Let $\Gamma$ be a point-to-set mapping from a set $\Theta$ into a set $\Xi$, i.e., $\Gamma\colon \Theta \to 2^{\Xi}$.
\begin{itemize}
	\item $\Gamma$ is said to be \emph{open} at $\bar{\theta} \in \Theta$
	if $\{\theta^{k}\} \subseteq \Theta$, $\theta^{k} \to \bar{\theta}$, and $\bar{\xi} \in \Gamma(\bar{\theta})$ imply
	the existence of an integer $m$ and a sequence $\{\xi^{k}\} \subseteq \Xi$ such that
	$\xi^{k} \in \Gamma(\theta^{k})$ for $k \ge m$ and $\xi^{k} \to \bar{\xi}$.
	\item $\Gamma$ is said to be \emph{closed} at $\bar{\theta} \in \Theta$
	if $\{\theta^{k}\} \subseteq \Theta$, $\theta^{k} \to \bar{\theta}$, $\xi^{k} \in \Gamma(\theta^{k})$, and $\xi^{k} \to \bar{\xi}$ imply
	$\bar{\xi} \in \Gamma(\bar{\theta})$.
	\item $\Gamma$ is said to be \emph{continuous} at $\bar{\theta} \in \Theta$
	if it is both open and closed at $\bar{\theta}$.
	If $\Gamma$ is continuous at every $\theta \in \Theta$, it is said to be continuous on $\Theta$.
	\item $\Gamma$ is said to be \emph{uniformly compact} near $\bar{\theta} \in \Theta$
	if there is a neighborhood $N$ of $\bar{\theta}$ such that the closure of the set $\bigcup_{\theta \in N} \Omega(\theta)$ is compact.
\end{itemize}
\end{definition}


\begin{lemma}[{\cite[Corollary 8.1]{Hog73}}]
\label{lem: Hog73}
For an objective function $\varphi\colon \bbR^{d} \times \Theta \to \bbR$ and
a constraint mapping $\Omega\colon \Theta \to 2^{\bbR^{d}}$,
consider the following parametric optimization problem with parameter $\theta \in \Theta \subseteq \bbR^{p}$:
\begin{equation}
\begin{aligned}
\min_{\xi}{}&
	\varphi(\xi, \theta)\\
\st{}&
	\xi \in \Omega(\theta).
\end{aligned}
\label{eq: param optim}
\end{equation}
For this problem, define the optimal set mapping $\Xi^{\ast}\colon \Theta \to 2^{\bbR^{d}}$ as
\begin{equation}
\Xi^{\ast}(\theta) \coloneqq \{\xi \in \Omega(\theta)\colon \varphi(\xi, \theta) = \inf_{\xi' \in \Omega(\theta)} \varphi(\xi', \theta)\}.
\end{equation}
Suppose that
$\Omega$ is continuous at $\bar{\theta} \in \Theta$,
$\varphi$ is continuous on $\Omega(\bar{\theta}) \times \{\bar{\theta}\}$,
$\Xi^{\ast}$ is nonempty and uniformly compact near $\bar{\theta}$, and
$\Xi^{\ast}(\bar{\theta})$ is a singleton.
Then, $\Xi^{\ast}$ is continuous at $\bar{\theta}$.
\end{lemma}

\begin{proof}[Proof of Theorem \ref{thm: existence}]
Since $S_{1}$ is compact from Assumption \ref{asp: existence}-\ref{enumi: compact},
it follows from the Weierstrass theorem that
a sufficient condition for the existence of minimizers is that $F_{1}$ is continuous.
Since $F_{1}$ is a composition of function $f_{1}$,
which is continuous from Assumption \ref{asp: existence}-\ref{enumi: continuous},
and functions $x_{2}^{\ast}, \dots, x_{n}^{\ast}$,
we inductively show the continuity of $x_{2}^{\ast}, \dots, x_{n}^{\ast}$ in what follows.


For any $i = 2, \dots, n$,
we derive the continuity of $x_{i}^{\ast}$ assuming the continuity of $x_{i + 1}^{\ast}, \dots, x_{n}^{\ast}$.
Since $f_{i}$ is continuous from Assumption \ref{asp: existence}-\ref{enumi: continuous} and
so are $x_{i + 1}^{\ast}, \dots, x_{n}^{\ast}$ from the induction hypothesis,
$F_{i}$ is continuous.
In addition, we define $X_{i}^{\ast}\colon S_{1} \times \bbR^{d_{2}} \times \dots \times \bbR^{d_{i - 1}} \to 2^{\bbR^{d_{i}}}$ as
\begin{equation}
X_{i}^{\ast}(x_{1}, \dots, x_{i - 1})
\coloneqq \left\{x_{i} \in \bbR^{d_{i}}\colon
F_{i}(x_{1}, \dots, x_{i - 1}, x_{i}) = \inf_{x_{i}' \in \bbR^{d_{n}}} F_{i}(x_{1}, \dots, x_{i - 1}, x_{i}')\right\}.
\end{equation}
For any $(\bar{x}_{1}, \dots, \bar{x}_{i - 1}) \in S_{1} \times \bbR^{d_{2}} \times \dots \bbR^{d_{i - 1}}$,
$X_{i}^{\ast}$ is nonempty near $(\bar{x}_{1}, \dots, \bar{x}_{i - 1})$
and $X_{i}^{\ast}(\bar{x}_{1}, \dots, \bar{x}_{i - 1})$ is a singleton from Assumption \ref{asp: existence}-\ref{enumi: singleton}; and
$X_{i}^{\ast}$ is uniformly bounded near $(\bar{x}_{1}, \dots, \bar{x}_{i - 1})$ from Assumption \ref{asp: existence}-\ref{enumi: bounded}.
Therefore, from Lemma \ref{lem: Hog73}, $X_{i}^{\ast}$ is continuous at $(\bar{x}_{1}, \dots, \bar{x}_{i - 1})$.
From the arbitrariness of $(\bar{x}_{1}, \dots, \bar{x}_{i - 1})$,
$X_{i}^{\ast}$ is continuous on $S_{1} \times \bbR^{d_{2}} \times \dots \bbR^{d_{i - 1}}$.
Recalling Assumption \ref{asp: existence}-\ref{enumi: singleton} again,
we have $X_{i}^{\ast}(x_{1}, \dots, x_{i - 1}) = \{x_{i}^{\ast}(x_{1}, \dots, x_{i - 1})\}$.
The continuity (in the sense of point-to-set mappings) and single-valuedness of $X_{i}^{\ast}$
means the continuity (in the sense of functions) of $x_{i}^{\ast}$.
\end{proof}

\subsection{Proof of Theorem \ref{thm: convergence}}
\label{sect: convergence proof}
To prove Theorem \ref{thm: convergence}, we use the following lemma.

\begin{lemma}[{\cite[Theorem A.1]{FFSGP18}}]
\label{lem: FFSGP18}
Let $\varphi_{T}$ and $\varphi$ be continuous%
\footnote{In the original statement of \cite[Theorem A.1]{FFSGP18},
the lower-semicontinuity is assumed, but in the proof, the continuity is used.
\cite[Theorem A.1]{FFSGP18} would be rewriting of \cite[Theorem 14]{DZ93}.
In the statement of \cite[Theorem 14]{DZ93}, the continuity is assumed.}
functions defined on a compact set $\Omega$.
Suppose that $\varphi_{T}$ converges uniformly to $\varphi$ on $\Omega$ as $T \to \infty$.
Then,
\begin{enumerate}[label = (\alph*)]
	\item $\inf_{\xi \in \Omega} \varphi_{T}(\xi) \to \inf_{\xi \in \Omega} \varphi(\xi)$ as $T \to \infty$;
	\item $\argmin_{\xi \in \Omega} \varphi_{T}(\xi) \to \argmin_{\xi \in \Omega} \varphi(\xi)$ as $T \to \infty$,
	meaning that,
	for every $\{\xi_{T}\}_{T \in \bbN}$ such that $\xi_{T} \in \argmin_{\xi \in \Omega} \varphi_{T}(\xi)$,
	we have that:
	\begin{itemize}
		\item $\{\xi_{T}\}$ admits a convergent subsequence;
		\item for every subsequence $\{\xi_{T_{k}}\}$ such that
		$\xi_{T_{k}} \to \bar{\xi}$ as $k \to \infty$,
		we have $\bar{\xi} \in \argmin_{\xi \in \Omega} \varphi(\xi)$.
	\end{itemize}
\end{enumerate}
\end{lemma}

\begin{proof}[Proof of Theorem \ref{thm: convergence}]
From Assumption \ref{asp: convergence}-\ref{enumi: upper unif conv},
there exists $\nu > 0$ such that for every $T \in \bbN$ and every $x_{1} \in S_{1}$,
\begin{align}
|\tilde{F}_{1}(x_{1}) - F_{1}(x_{1})|
&= |f_{1}(x_{1}, x_{2}^{(T)}(x_{1}), \dots, x_{n}^{(T)}(x_{1})) - f_{1}(x_{1}, x_{2}^{\ast}(x_{1}), \dots, x_{n}^{\ast}(x_{1}))|\\
&\le \nu \sum_{i = 2}^{n} \|x_{i}^{(T)}(x_{1}) - x_{i}^{\ast}(x_{1})\|.
\end{align}
It follows from Assumption \ref{asp: convergence}-\ref{enumi: lower unif conv} that
$\tilde{F}_{1}$, which depends on $T$, converges to $F_{1}$ uniformly on $S_{1}$ as $T \to \infty$.
Then, the statement follows from Lemma \ref{lem: FFSGP18}.
Note that the continuity of $F_{1}$ has already shown in Theorem \ref{thm: existence}.
\end{proof}

\subsection{Proof of Theorem \ref{thm: grad formula}}
\label{sect: grad formula proof}

\begin{proof}
By using the chain rule, we obtain
\begin{equation}
\nabla_{x_{1}} \tilde{F}_{1}(x_{1}) = \nabla_{x_{1}} f_{1} + \sum_{i = 2}^{n}
\nabla_{x_{1}} x_{i}^{(T_{i})}(x_{1}, x_{2}^{(T_{2})}(x_{1}), \dots,
x_{i - 1}^{(T_{i - 1})}(x_{k - 2}^{(T_{k - 2})}(\cdots(x_{2}^{(T_{2})}(x_{1})))))
\nabla_{x_{i}} f_{1},
\end{equation}
where we denoted
$\nabla_{x_{i}} f_{1} = \nabla_{x_{i}} f_{1}(x_{1}, x_{2}^{(T_{2})}, \dots, x_{n}^{(T_{n})})$ for $i = 2, \dots, n$ and
$\nabla_{x_{1}} f_{1} = \nabla_{x_{1}} f_{1}(x_{1}, x_{2}^{(T_{2})}, \dots, x_{n}^{(T_{n})})$.
Hereinafter, for $i = 2, \dots, n$ and $t = 1, \dots, T_{i}$,
we denote $Y_{i}^{(t)} = \nabla_{x_{1}} x_{i}^{(t)}(x_{1}, x_{2}^{(T_{2})}(x_{1}), \dots,
x_{i - 1}^{(T_{i - 1})}(x_{k - 2}^{(T_{k - 2})}(\cdots(x_{2}^{(T_{2})}(x_{1})))))$,
and show $Y_{i}^{(T_{i})} = Z_{i}$ for all $i = 2, \dots, n$.
In what follows, we denote $\bar{B}_{i}^{(t)} = \sum_{j = 2}^{i - 1} Z_{j} C_{i j}^{(t)} + B_{i}^{(t)}$.

For $i = 3, \dots, n $ and $t = 1, \dots, T_{i}$,
by applying the chain rule to the update formula of $x_{i}$, i.e.,
$x_{i}^{(t)} = \Phi_{i}^{(t)}(x_{1}, x_{2}^{(T_{2})}, \dots, x_{i - 1}^{(T_{i - 1})}, x_{i}^{(t - 1)})$,
we obtain
\begin{align}
Y_{i}^{(t)}
&= Y_{i}^{(t - 1)} \nabla_{x_{i}} \Phi_{i}^{(t)}(x_{1}, x_{2}^{(T_{2})}, \dots, x_{i - 1}^{(T_{i - 1})}, x_{i}^{(t - 1)})\\
&\qquad + \sum_{j = 2}^{i - 1} Y_{j}^{(T_{j})}
\nabla_{x_{j}} \Phi_{i}^{(t)}(x_{1}, x_{2}^{(T_{2})}, \dots, x_{i - 1}^{(T_{i - 1})}, x_{i}^{(t - 1)})\\
&\qquad + \nabla_{x_{1}} \Phi_{i}^{(t)}(x_{1}, x_{2}^{(T_{2})}, \dots, x_{i - 1}^{(T_{i - 1})}, x_{i}^{(t - 1)})\\
&= Y_{i}^{(t - 1)} A_{i}^{(t)} + \sum_{j = 2}^{i - 1} Y_{j}^{(T_{j})} C_{i j}^{(t)} + B_{i}^{(t)}.
\label{eq: Y_{i}^{(t)}}
\end{align}
For $i = 2$ and $t = 1, \dots, T_{2}$,
note that the arguments of $\Phi_{2}^{(t)}$ are $x_{1}$ and $x_{2}^{(t - 1)}$.
For $t = 1, \dots, T_{2}$, by applying the chain rule, we obtain the following:
\begin{align}
Y_{2}^{(t)}
&= Y_{2}^{(t - 1)} \nabla_{x_{2}} \Phi_{2}^{(t)}(x_{1}, x_{2}^{(t - 1)}) + \nabla_{x_{1}} \Phi_{2}^{(t)}(x_{1}, x_{2}^{(t - 1)})\\
&= Y_{2}^{(t - 1)} A_{2}^{(t)} + B_{2}^{(t)}\\
&= Y_{2}^{(t - 1)} A_{2}^{(t)} + \bar{B}_{2}^{(t)}.
\label{eq: Y_{2}^{(t)}}
\end{align}
Note that, $Y_{i}^{(0)} = \nabla_{x_{1}} x_{i}^{(0)} = O$ for $i = 2, \dots, n$
since the initial point $x_{i}^{(0)}$ is independent to the other variables.
Using \eqref{eq: Y_{2}^{(t)}} repeatedly, we obtain
\begin{align}
Y_{2}^{(T_{2})}
&= Y_{2}^{(T_{2} - 1)} A_{2}^{(T_{2})} + \bar{B}_{2}^{(T_{2})}\\
&= (Y_{2}^{(T_{2} - 2)} A_{2}^{(T_{2} - 1)} + \bar{B}_{2}^{(T_{2} - 1)}) A_{2}^{(T_{2})} + \bar{B}_{2}^{(T_{2})}\\
&= ((Y_{2}^{(T_{2} - 3)} A_{2}^{(T_{2} - 2)} + \bar{B}_{2}^{(T_{2} - 2)}) A_{2}^{(T_{2} - 1)}
+ \bar{B}_{2}^{(T_{2} - 1)}) A_{2}^{(T_{2})} + \bar{B}_{2}^{(T_{2})}\\
&= \dots\\
&= (((\cdots(\bar{B}_{2}^{(1)} A_{2}^{(2)} + \bar{B}_{2}^{(2)})\cdots) A_{2}^{(T_{2} - 2)}
+ \bar{B}_{2}^{(T_{2} - 2)}) A_{2}^{(T_{2} - 1)} + \bar{B}_{2}^{(T_{2} - 1)}) A_{2}^{(T_{2})} + \bar{B}_{2}^{(T_{2})}\\
&= \sum_{t = 1}^{T_{2}} \bar{B}_{2}^{(t)} \prod_{s = t + 1}^{T_{2}} A_{2}^{(s)}\\
& = Z_{2}.
\label{eq: Z_{2}}
\end{align}
Next, for $t = 1, \dots, T_{3}$,
plugging Equation \eqref{eq: Y_{i}^{(t)}} with $i = 3$ into Equation \eqref{eq: Z_{2}},
we obtain
\begin{equation}
Y_{3}^{(t)}
= Y_{3}^{(t - 1)} A_{3}^{(t)} + Z_{2} C_{2 2}^{(t)} + B_{3}^{(t)}
= Y_{3}^{(t - 1)} A_{3}^{(t)} + \bar{B}_{3}^{(t)}.
\label{eq: Y_{3}^{(t)}}
\end{equation}
Repeatedly using Equation \eqref{eq: Y_{3}^{(t)}} similarly to the derivation of Equation \eqref{eq: Z_{2}},
we obtain
\begin{equation}
Y_{3}^{(T_{3})}
= \sum_{t = 1}^{T_{3}} \bar{B}_{3}^{(t)} \prod_{s = t + 1}^{T_{3}} A_{3}^{(s)}
= Z_{3}.
\end{equation}
For $i \ge 4$, we can use a similar argument incrementing $i$:
When we evaluate $Y_{i}^{({T_{i}})}$,
since we already have $Y_{j}^{(T_{j})} = Z_{j}$ for $l = 2, \dots, k - 1$,
we can derive the followings:
\begin{align}
Y_{i}^{(t)}
&= Y_{i}^{(t - 1)} A_{i}^{(t)} + \sum_{j = 2}^{i - 1} Y_{j}^{(T_{j})} C_{i j}^{(t)} + B_{i}^{(t)}\\
&= Y_{i}^{(t - 1)} A_{i}^{(t)} + \sum_{j = 2}^{i - 1} Z_{j} C_{i j}^{(t)} + B_{i}^{(t)}\\
&= Y_{i}^{(t - 1)} A_{i}^{(t)} + \bar{B}_{i}^{(t)}.
\label{eq: Y_{i}^{(t)}-bar{B}}
\end{align}
Repeatedly using Equation \eqref{eq: Y_{i}^{(t)}-bar{B}} similarly to the derivation of Equation \eqref{eq: Z_{2}},
we obtain $Y_{i}^{(T_{i})} = Z_{i}$ for $i = 2, \dots, n$.
\end{proof}

\subsection{Proof of Theorem \ref{thm: complexity recursive}}
\label{sec:complexity}

To prove Theorem \ref{thm: complexity recursive}, we introduce the following lemma.

\begin{lemma}[Complexity of Jacobian-vector products \cite{BPRS17, FDFP17}]
\label{thm: Jacobi-vec prod comp}
Let $\varphi\colon \bbR^{d} \to \bbR^{m}$ be a differentiable function
and suppose it can be evaluated in time $c(d, m)$ and requires space $s(d, m)$.
Denote let $J_{\varphi} \in \bbR^{m \times d}$ be the Jacobian matrix of $\varphi$.
Then, for any vector $\xi \in \bbR^{d}$,
the product $J_{\varphi} \xi$ can be computed within $\bigO(c(d, m))$ time complexity and $\bigO(s(d, m))$ space complexity.
\end{lemma}

On the complexity of Algorithm \ref{algo: grad comp}, the following lemma holds.
Hereinafter, $n$ times nested $\bigO(\cdot)$ is denoted by $\bigO^{n}(\cdot)$,
that is, $\bigO^{n}(\cdot) = \bigO(\bigO(\cdots\bigO(\cdot)\cdots))$.
Note that, for a parameter $n$, $\bigO^{n}$ is not equivalent to $\bigO$
while $\bigO^{c}$ is equivalent to $\bigO$ for a constant $c$.

\begin{lemma}[The complexity of Algorithm \ref{algo: grad comp}]
\label{thm: complexity}
Let the time and space complexity for computing $\nabla_{x_{i}} \tilde{F}_{i}(x_{i})$ be $c_{i}$ and $s_{i}$, respectively.
We assume that those for evaluating $\Phi_{i}^{(t)}$ are $\bigO(c_{i})$ and $\bigO(s_{i})$.
In addition, we assume that the time and space complexity of $\nabla_{x_{1}} f_{1}$ and $\nabla_{x_{i}} f_{1}$ for $i = 1, \dots, n$
are smaller in the sense of the order than those of \texttt{for} loops in lines 2--9 in Algorithm \ref{algo: grad comp}.
Then, the time and space complexity of Algorithm \ref{algo: grad comp} can be written as
\begin{alignat}{2}
c_{1}
&= \sum_{i = 2}^{n} \bigO(T_{i} i d_{1} \bigO(c_{i})),&\
s_{1}
&= \max_{i = 2}^{n} \bigO^{2}(s_{i}).
\label{eq: complexity}
\end{alignat}
\end{lemma}

\begin{proof}
For each $i = 2, \dots, n$ and $t = 1, \dots, T_{i}$,
the complexity of lines 5--7 of Algorithm \ref{algo: grad comp} can be evaluated as follows:
\begin{itemize}
	\item From the assumption, the time and space complexity of the 5th line are $\bigO(c_{i})$ and $\bigO(s_{i})$, respectively.
	\item The 6th line computes $i - 2$ products of Jacobian matrices and a $d_{i} \times d_{1}$ matrices.
	One product is obtained by computing the product of a Jacobian matrix and a vector $d_{1}$ times
	and hence its time and space complexity are $\bigO(d_{1} \bigO(c_{i}))$ and $\bigO(\bigO(s_{i}))$, respectively.
	Since the 6th line compute $i - 2$ products, the overall time and space complexity of the 6th line are
	$\bigO(i d_{1} \bigO(c_{i}))$ and $\bigO(\bigO(s_{i}))$, respectively.
	\item The 7th line computes the product of a Jacobian matrix and a $d_{i} \times d_{1}$ matrix.
	Its time and space complexity are $\bigO(d_{1} \bigO(c_{i}))$ and $\bigO(\bigO(s_{i}))$, respectively.
\end{itemize}
Hence, for each $i = 2, \dots, n$, the time and space complexity of the single run of the \texttt{for} loop in lines 2--7 are
$\bigO(T_{i} i d_{1} \bigO(c_{i}))$, $\bigO(\bigO(s_{i}))$ respectively.
Therefore, we obtain Equation \eqref{eq: complexity} as the overall time and space complexity of Algorithm \ref{algo: grad comp}.
\end{proof}

\begin{proof}[Proof of Theorem \ref{thm: complexity recursive}]
From the assumption, $c_{i}$, $s_{i}$ are the time and space complexity for computing $\nabla_{x_{i}} \tilde{F}_{i}(x_{i})$.
Hence, from Lemma \ref{thm: complexity}, $c_{i}$ and $s_{i}$ are written as
\begin{alignat}{2}
c_{i} &= \sum_{j = i + 1}^{n} \bigO(T_{j} {j} d_{j} \bigO(c_{j})),&\quad
s_{i} &= \max_{j = i + 1}^{n} \bigO^{2}(s_{j}).
\end{alignat}
From the equation above,
$c_{i}$ ($i = 1, \dots, n - 1$) can be written by using $c_{n}$ as follows:
\begin{align}
c_{n - 1}
&= \bigO(T_{n} {n} d_{n - 1} \bigO(c_{n})),\\
c_{n - 2}
&= \bigO(T_{n - 1} (n - 1) d_{n - 2} \bigO(c_{n - 1})) + \bigO(T_{n} n d_{n - 2} \bigO(c_{n}))\\
&= \bigO(T_{n - 1} (n - 1) d_{n - 2} \bigO(\bigO(T_{n} n d_{n - 1} \bigO(c_{n}))))\\
&= \bigO(T_{n} T_{n - 1} n (n - 1) d_{n - 1} d_{n - 2} \bigO(\bigO(\bigO(c_{n})))),\\
&\vdots\\
c_{i}
&= \bigO(T_{n} T_{n - 1} \cdots T_{i + 1} n (n - 1) \cdots (i + 1) d_{n - 1} d_{n - 2} \cdots d_{i} \bigO^{2 (n - i) - 1}(c_{n})),\\
&\vdots\\
c_{1}
&= \bigO(T_{n} T_{n - 1} \cdots T_{2} n (n - 1) \cdots 2 d_{n - 1} d_{n - 2} \cdots d_{1} \bigO^{2 n - 3}(c_{n}))\\
&= \bigO\left(\bigO^{2 n - 3}(c_{n}) n! \prod_{i = 1}^{n - 1} (T_{i + 1} d_{i})\right).
\end{align}
Similarly, $s_{i}$ ($i = 1, \dots, n - 1$) can be written by using $s_{n}$ as follows:
\begin{align}
s_{n - 1}
&= \bigO^{2}(s_{n}),\\
s_{n - 2}
&= \max_{l = n - 1}^{n} \bigO^{2}(s_{l})
= \max\{\bigO^{2}(s_{l}), \bigO^{4}(s_{l})\}
= \bigO^{4}(s_{n}), \\
&\vdots\\
s_{i} &=
\bigO^{2 (n - i)}(s_{n}),\\
&\vdots\\
s_{1}
&= \bigO^{2 n - 2}(s_{n}).
\end{align}
Note that, since the gradient of the objective function in the lower level is used to
calculate the gradient of the objective function in the upper level, the amount of
calculation of the gradient at each level increases from the bottom to the top.
Therefore, for some $p, q > 1$,
it holds that $\bigO^{2 n - 3}(c_{n}) = \bigO(p^{n} c_{n})$ and $\bigO^{2 n - 2}(s_{n}) = \bigO(q^{n} s_{n})$.
Hence, we obtain Equation \eqref{eq: complexity recursive}
as the overall time and space complexity for computing $\nabla_{x_{1}} \tilde{F}_{1}(x_{1})$
recursively calling Algorithm \ref{algo: grad comp}.
\end{proof}

\subsection{Proof of Theorem \ref{thm: Lip cont grad}}
\label{sect: Lip cont grad proof}

\begin{proof}
In this proof,
we inductively show the Lipschitz continuity of $\nabla_{x_{i}} \tilde{F}_{i}$ for all $i = 1, \dots, n$
by using the following notations:
When a function $\varphi$ is Lipschitz continuous on $\Omega$,
its Lipschitz constant is denoted by $L_{\varphi}$, i.e.,
$\|\varphi(\hat{\xi}) - \varphi(\check{\xi})\| \le L_{\varphi} \|\hat{\xi} - \check{\xi}\|$ for any $\hat{\xi}, \check{\xi} \in \Omega$.
When $\varphi$ is bounded on $\Omega$,
its bound is denoted by $M_{\varphi}$, i.e.,
$\|\varphi(\xi)\| \le M_{\varphi}$ for any $\xi \in \Omega$.

As the base case, we consider the case of $i = n$.
Since $\tilde{F}_{n} = f_{n}$,
the Lipschitz continuity of $\nabla_{x_{n}} \tilde{F}_{n}$ follows
from the assumption of the Lipschitz continuity of $\nabla_{x_{i}} f_{n}$.

In what follows, as the induction step,
we derive the Lipschitz continuity of $\nabla_{x_{i}} \tilde{F}_{i}$
by assuming the Lipschitz continuity of $\nabla_{x_{j}} \tilde{F}_{j}$ for $j = i + 1, \dots, n$.
To do so, for any $(\hat{x}_{1}, \dots, \hat{x}_{i})$ and $(\check{x}_{1}, \dots, \check{x}_{i})$
in $S_{1} \times \bbR^{d_{2}} \times \dots \times\bbR^{d_{i}}$,
we evaluate $\|\nabla_{x_{i}} \tilde{F}_{i}(\hat{x}_{1}, \dots, \hat{x}_{i})
- \nabla_{x_{i}} \tilde{F}_{i}(\check{x}_{1}, \dots, \check{x}_{i})\|$.
In what follows, constants related to $(\hat{x}_{1}, \dots, \hat{x}_{i})$ and $(\check{x}_{1}, \dots, \check{x}_{i})$
are denoted by ones with hat and check sign, e.g., $\hat{Z}_{j}$ and $\check{Z}_{j}$.
By applying Theorem \ref{thm: grad formula} to $\tilde{F}_{i}$, we have
\begin{align}
&\|\nabla_{x_{i}} \tilde{F}_{i}(\hat{x}_{1}, \dots, \hat{x}_{i}) - \nabla_{x_{i}} \tilde{F}_{i}(\check{x}_{1}, \dots, \check{x}_{i})\|\\
{} \le {}
&\|\nabla_{x_{i}} f_{i}(\hat{x}_{1}, \dots, \hat{x}_{i}, \hat{x}_{i + 1}^{(T_{i + 1})}, \dots, \hat{x}_{n}^{(T_{n})})
- \nabla_{x_{i}} f_{i}(\check{x}_{1}, \dots, \check{x}_{i}, \check{x}_{i + 1}^{(T_{i + 1})}, \dots, \check{x}_{n}^{(T_{n})})\|\\
&{} + \sum_{j = 2}^{n}
\|\hat{Z}_{j} \nabla_{x_{j}} f_{i}(\hat{x}_{1}, \dots, \hat{x}_{i}, \hat{x}_{i + 1}^{(T_{i + 1})}, \dots, \hat{x}_{n}^{(T_{n})})
- \check{Z}_{j} \nabla_{x_{j}} f_{i}(\check{x}_{1}, \dots, \check{x}_{i}, \check{x}_{i + 1}^{(T_{i + 1})}, \dots, \check{x}_{n}^{(T_{n})})\|
\label{eq: Lipschitz RHS}
\end{align}

First, we evaluate the first term of the right-hand-side of Equation \eqref{eq: Lipschitz RHS}.
From the assumption of the Lipschitz continuity of $\nabla_{x_{i}} f_{i}$, we have
\begin{align}
&\|\nabla_{x_{i}} f_{i}(\hat{x}_{1}, \dots, \hat{x}_{i}, \hat{x}_{i + 1}^{(T_{i + 1})}, \dots, \hat{x}_{n}^{(T_{n})})
- \nabla_{x_{i}} f_{i}(\check{x}_{1}, \dots, \check{x}_{i}, \check{x}_{i + 1}^{(T_{i + 1})}, \dots, \check{x}_{n}^{(T_{n})})\|\\
{} \le {}
&L_{\nabla_{x_{i}} f_{i}} \left(\|(\hat{x}_{1}, \dots, \hat{x}_{i}) - (\check{x}_{1}, \dots, \check{x}_{i})\|
+ \sum_{j = i + 1}^{n} \|\hat{x}_{j}^{(T_{j})} - \check{x}_{j}^{(T_{j})}\|\right).
\end{align}
For any $j = i + 1, \dots, n$ and $t = 1, \dots, T_{j}$,
recalling the definition of $\hat{x}_{j}^{(t)}$ and $\check{x}_{j}^{(t)}$
and invoking the Lipschitz continuity of $\nabla_{x_{j}} \tilde{F}_{j}$,
we have
\begin{align}
&\|\hat{x}_{j}^{(t)} - \check{x}_{j}^{(t)}\|\\
{} = {}
&\|\Phi_{j}^{(t)}(\hat{x}_{1}, \dots, \hat{x}_{i},
\hat{x}_{i + 1}^{(T_{i + 1})}, \dots, \hat{x}_{j - 1}^{(T_{j - 1})}, \hat{x}_{j}^{(t - 1)})
- \Phi_{j}^{(t)}(\check{x}_{1}, \dots, \check{x}_{i},
\check{x}_{i + 1}^{(T_{i + 1})}, \dots, \check{x}_{j - 1}^{(T_{j - 1})}, \check{x}_{j}^{(t - 1)})\|\\
{} \le {}
&\|\hat{x}_{j}^{(t - 1)} - \check{x}_{j}^{(t - 1)}\|
+ \alpha_{j}^{(t - 1)} \left\|
\begin{multlined}
\nabla_{x_{j}} \tilde{F}_{j}(\hat{x}_{1}, \dots, \hat{x}_{i},
\hat{x}_{i + 1}^{(T_{i + 1})}, \dots, \hat{x}_{j - 1}^{(T_{j - 1})}, \hat{x}_{j}^{(t - 1)})\\
{} - \nabla_{x_{j}} \tilde{F}_{j}(\check{x}_{1}, \dots, \check{x}_{i},
\check{x}_{i + 1}^{(T_{i + 1})}, \dots, \check{x}_{j - 1}^{(T_{j - 1})}, \check{x}_{j}^{(t - 1)})
\end{multlined}
\right\|\\
{} \le {}
&\|\hat{x}_{j}^{(t - 1)} - \check{x}_{j}^{(t - 1)}\|
+ \alpha_{j}^{(t - 1)} L_{\nabla_{x_{j}} \tilde{F}_{j}}
\left(
\begin{multlined}
\|(\hat{x}_{1}, \dots, \hat{x}_{i}) - (\check{x}_{1}, \dots, \check{x}_{i})\|\\
{} + \sum_{k = i + 1}^{j - 1} \|\hat{x}_{k}^{(T_{k})} - \check{x}_{k}^{(T_{k})}\|
+ \|\hat{x}_{j}^{(t - 1)} - \check{x}_{j}^{(t - 1)}\|
\end{multlined}\right).
\end{align}
In addition, for any $j = i + 1, \dots, n$, we have
\begin{align}
&\|\hat{x}_{j}^{(1)} - \check{x}_{j}^{(1)}\|\\
{} = {}
&\|\Phi_{j}^{(1)}(\hat{x}_{1}, \dots, \hat{x}_{i},
\hat{x}_{i + 1}^{(T_{i + 1})}, \dots, \hat{x}_{j - 1}^{(T_{j - 1})}, x_{j}^{(0)})
- \Phi_{i}^{(1)}(\check{x}_{1}, \dots, \check{x}_{i},
\check{x}_{i + 1}^{(T_{i + 1})}, \dots, \check{x}_{j - 1}^{(T_{j - 1})}, \check{x}_{j}^{(0)})\|\\
{} \le {}
&\underbrace{\|x_{j}^{(0)} - \check{x}_{j}^{(0)}\|}_{{} = 0} + \alpha_{j}^{(0)} L_{\nabla_{x_{i}} \tilde{F}_{i}}
\left(
\begin{multlined}
\|(\hat{x}_{1}, \dots, \hat{x}_{i}) - (\check{x}_{1}, \dots, \check{x}_{i})\|\\
{} + \sum_{k = i + 1}^{j - 1} \|\hat{x}_{k}^{(T_{k})} - \check{x}_{k}^{(T_{k})}\|
+ \underbrace{\|\hat{x}_{j}^{(0)} - \check{x}_{j}^{(0)}\|}_{{} = 0}
\end{multlined}
\right).
\end{align}
Using these inequalities repeatedly, we can derive
\begin{align}
&\|\nabla_{x_{i}} f_{i}(\hat{x}_{1}, \dots, \hat{x}_{i}, \hat{x}_{i + 1}^{(T_{i + 1})}, \dots, \hat{x}_{n}^{(T_{n})})
- \nabla_{x_{i}} f_{i}(\check{x}_{1}, \dots, \check{x}_{i}, \check{x}_{i + 1}^{(T_{i + 1})}, \dots, \check{x}_{n}^{(T_{n})})\|\\
{} \le {}
&N_{1} \|(\hat{x}_{1}, \dots, \hat{x}_{i}) - (\check{x}_{1}, \dots, \check{x}_{i})\|
\end{align}
for some positive constant $N_{1}$ depending only on
$\{\alpha_{j}^{(t - 1)}\}$, $L_{\nabla_{x_{i}} f_{i}}$ and $L_{\nabla_{x_{i}} \tilde{F}_{j}}$.

Next, we evaluate each summand in the second term of the right-hand-side of Equation \eqref{eq: Lipschitz RHS}.
To do so, we notice that $\hat{Z}_{j}$ ($\check{Z}_{j}$, respectively) is the sum of products of
$\hat{A}_{j}^{(t)}$'s, $\hat{B}_{j}^{(t)}$'s, and $\hat{C}_{j k}^{(t)}$'s
($\check{A}_{j}^{(t)}$'s, $\check{B}_{j}^{(t)}$'s, and $\check{C}_{j k}^{(t)}$'s, respectively).
Formally, we can decompose $\hat{Z}_{j} = \sum_{p = 1}^{P} \prod_{q = 1}^{Q_{p}} \hat{X}_{j p q}$,
where $\hat{X}_{j p q} \in \bigcup_{t = 1}^{T_{j}}
(\{\hat{A}_{j}^{(t)}, \hat{B}_{j}^{(t)}\} \cup \bigcup_{k = 1}^{j - 1} \{\hat{C}_{j k}^{(t)}\})$.
Similarly, by appropriately defining $\check{X}_{j p q}$,
we can rewrite $\check{Z}_{j} = \sum_{p = 1}^{P} \prod_{q = 1}^{Q_{p}} \check{X}_{j p q}$,
where $\check{X}_{j p q} \in \bigcup_{t = 1}^{T_{j}}
(\{\check{A}_{j}^{(t)}, \check{B}_{j}^{(t)}\} \cup \bigcup_{k = 1}^{j - 1} \{\check{C}_{j k}^{(t)}\})$ and
$\check{X}_{j p q} = \check{A}_{j}^{(t)}, \check{B}_{j}^{(t)}, \check{C}_{j k}^{(t)}$
if $\hat{X}_{j p q} = \hat{A}_{j}^{(t)}, \hat{B}_{j}^{(t)}, \hat{C}_{j k}^{(t)}$, respectively.
Using this notation, we have
\begin{align}
&\|\hat{Z}_{j} \nabla_{x_{j}} f_{i}(\hat{x}_{1}, \dots, \hat{x}_{i}, \hat{x}_{i + 1}^{(T_{i + 1})}, \dots, \hat{x}_{n}^{(T_{n})})
- \check{Z}_{j} \nabla_{x_{j}} f_{i}(\check{x}_{1}, \dots, \check{x}_{i}, \check{x}_{i + 1}^{(T_{i + 1})}, \dots, \check{x}_{n}^{(T_{n})})\|
\\
{} \le {}
&\sum_{p = 1}^{P}
\left\|
\prod_{q = 1}^{Q_{p}} \hat{X}_{j p q}
\nabla_{x_{j}} f_{i}(\hat{x}_{1}, \dots, \hat{x}_{i}, \hat{x}_{i + 1}^{(T_{i + 1})}, \dots, \hat{x}_{n}^{(T_{n})})
- \prod_{q = 1}^{Q_{p}} \check{X}_{j p q}
\nabla_{x_{j}} f_{i}(\check{x}_{1}, \dots, \check{x}_{i}, \check{x}_{i + 1}^{(T_{i + 1})}, \dots, \check{x}_{n}^{(T_{n})})
\right\|\\
{} \le {}
&\sum_{p = 1}^{P} \left(
\begin{multlined}
\sum_{r = 1}^{Q_{p}} \left\|
\begin{multlined}
\prod_{q = 1}^{r - 1} \check{X}_{j p q} \prod_{q = r}^{Q_{p}} \hat{X}_{j p q}
\nabla_{x_{j}} f_{i}(\hat{x}_{1}, \dots, \hat{x}_{i}, \hat{x}_{i + 1}^{(T_{i + 1})}, \dots, \hat{x}_{n}^{(T_{n})})\\
{} - \prod_{q = 1}^{r} \check{X}_{j p q} \prod_{q = r + 1}^{Q_{p}} \hat{X}_{j p q}
\nabla_{x_{j}} f_{i}(\hat{x}_{1}, \dots, \hat{x}_{i}, \hat{x}_{i + 1}^{(T_{i + 1})}, \dots, \hat{x}_{n}^{(T_{n})})
\end{multlined}
\right\|\\
{} + \left\|
\begin{multlined}\prod_{q = 1}^{Q_{p}} \check{X}_{j p q}
\nabla_{x_{j}} f_{i}(\hat{x}_{1}, \dots, \hat{x}_{i}, \hat{x}_{i + 1}^{(T_{i + 1})}, \dots, \hat{x}_{n}^{(T_{n})})\\
{} - \prod_{q = 1}^{Q_{p}} \check{X}_{j p q}
\nabla_{x_{j}} f_{i}(\check{x}_{1}, \dots, \check{x}_{i}, \check{x}_{i + 1}^{(T_{i + 1})}, \dots, \check{x}_{n}^{(T_{n})})
\end{multlined}\right\|
\end{multlined}
\right)\\
{} \le {}
&\sum_{p = 1}^{P} \left(
\begin{multlined}
\sum_{r = 1}^{Q_{p}}
\|\hat{X}_{j p r} - \check{X}_{j p r}\|
\|\nabla_{x_{j}} f_{i}(\hat{x}_{1}, \dots, \hat{x}_{i}, \hat{x}_{i + 1}^{(T_{i + 1})}, \dots, \hat{x}_{n}^{(T_{n})})\|
\prod_{q = 1}^{r - 1} \|\check{X}_{j p q}\| \prod_{q = r + 1}^{Q_{p}} \|\hat{X}_{j p q}\|\\
{} + \prod_{q = 1}^{Q_{p}} \|\check{X}_{j p q}\|
\left\|\nabla_{x_{j}} f_{i}(\hat{x}_{1}, \dots, \hat{x}_{i}, \hat{x}_{i + 1}^{(T_{i + 1})}, \dots, \hat{x}_{n}^{(T_{n})})
- \nabla_{x_{j}} f_{i}(\check{x}_{1}, \dots, \check{x}_{i}, \check{x}_{i + 1}^{(T_{i + 1})}, \dots, \check{x}_{n}^{(T_{n})})\right\|
\end{multlined}
\right)\\
{} \le {}
&\underbrace{\sum_{p = 1}^{P} \left(\sum_{r = 1}^{Q_{p}} L_{X_{j p r}} M_{\nabla_{x_{j}} f_{i}} \prod_{q \ne r} M_{X_{j p q}}
+ L_{\nabla_{x_{j}} f_{i}} \prod_{q = 1}^{Q_{p}} M_{X_{j p q}}\right)}_{{} = N_{2}}
\|(\hat{x}_{1}, \dots, \hat{x}_{i}) - (\check{x}_{1}, \dots, \check{x}_{i})\|,
\end{align}
where, in the last inequality, we used the Lipschitz continuity and the boundedness of
$\nabla_{x_{i}} \Phi_{i}^{(t)}$, $\nabla_{x_{1}} \Phi_{i}^{(t)}$, and $\nabla_{x_{j}} \Phi_{i}^{(t)}$.
Therefore, we obtain
\begin{equation}
\|\nabla_{x_{i}} \tilde{F}_{i}(\hat{x}_{1}, \dots, \hat{x}_{i}) - \nabla_{x_{i}} \tilde{F}_{i}(\check{x}_{1}, \dots, \check{x}_{i})\|
\le \left(N_{1} + \sum_{j = 2}^{n} N_{2 j}\right) \|(\hat{x}_{1}, \dots, \hat{x}_{i}) - (\check{x}_{1}, \dots, \check{x}_{i})\|.
\end{equation}
This means the Lipschitz continuity of $\nabla_{x_{i}} \tilde{F}_{i}$,
and hence, that of $\nabla_{x_{1}} \tilde{F}_{1}$ is obtained by induction.
\end{proof}

\section{An example of applying Algorithm \ref{algo: grad comp}}
\label{sect: trilevel grad comp}
We explain an example of approximated problems to be solved by Algorithm \ref{algo: grad comp}
(i.e., problem \ref{eq: approx} with $n = 3$) by assuming a simple setting,
where we apply the steepest descent method for the lower-level problems
with the same iteration number $T$ and step size $\alpha$
for all levels as follows:
\begin{equation}
\begin{alignedat}{2}
\min_{x_{1} \in S_{1}, \{x_{2}^{(t)}\}, \{x_{3}^{(t)}\}}{}&
	f_{1}(x_{1}, x_{2}^{(T)}, x_{3}^{(T)})\\
\text{s.t. }&
	x_{2}^{(t)}
	= x_{2}^{(t - 1)} - \alpha \nabla_{x_{2}} \tilde{F}_{2}(x_{1}, x_{2}^{(t - 1)}),&
		\quad& (t = 1, \dots, T),\\
& 	x_{3}^{(t)}
	= x_{3}^{(t - 1)} - \alpha \nabla_{x_{3}} \tilde{F}_{3}(x_{1}, x_{2}^{(T)}, x_{3}^{(t - 1)})&
		\quad&	(t = 1, \dots, T).
\end{alignedat}
\label{eq: trilevel instance}
\end{equation}
Here, $\nabla_{x_{2}} \tilde{F}_{2}(x_{1}, x_{2})$
is the gradient of the objective function of
$\min_{x_{2} \in \mathbb{R}^{d_{2}}} \tilde{F}_{2}(x_{1}, x_{2})$,
which is equivalent to the bilevel optimization problem as follows:
\begin{equation}
\begin{alignedat}{2}
\min_{x_{2}^{(t)}, \{x_{3}^{(t')}\}}{}&
	f_{2}(x_{1}, x_{2}^{(t)}, x_{3}^{(T)})\\
\text{s.t. }&
	x_{3}^{(t')}
	= x_{3}^{(t' - 1)} - \alpha \nabla_{x_{3}} \tilde{F}_{3}(x_{1}, x_{2}^{(t)}, x_{3}^{(t' - 1)})&
		\quad&	(t' = 1, \dots, T).
\end{alignedat}
\end{equation}
In addition, since the third level problem is the lowest level problem,
$\tilde{F}_{3} = f_{3}$ holds.
We replace $x_{2}^{(T)}$ in the objective function of \eqref{eq: trilevel instance} by
$x_{2}^{(T - 1)} - \alpha \nabla_{x_{2}} \tilde{F}_{2}(x_{1}, x_{2}^{(T - 1)})$
and then replace recursively $x_{2}^{(t)}$  using $\{x_{2}^{(t)}\}_{t = 0}^{T - 1}$.
By applying the same procedure for $x_{3}^{(T)}$,
we can reformulate \eqref{eq: trilevel instance} to $\min_{x_{1} \in S_{1}} \tilde{F}_{1}(x_{1})$.
Theorem \ref{thm: convergence} proves that this problem
converges to the trilevel optimization problem as $T \rightarrow \infty$.

\begin{figure}
\centering
\includegraphics[width = 0.95\textwidth]{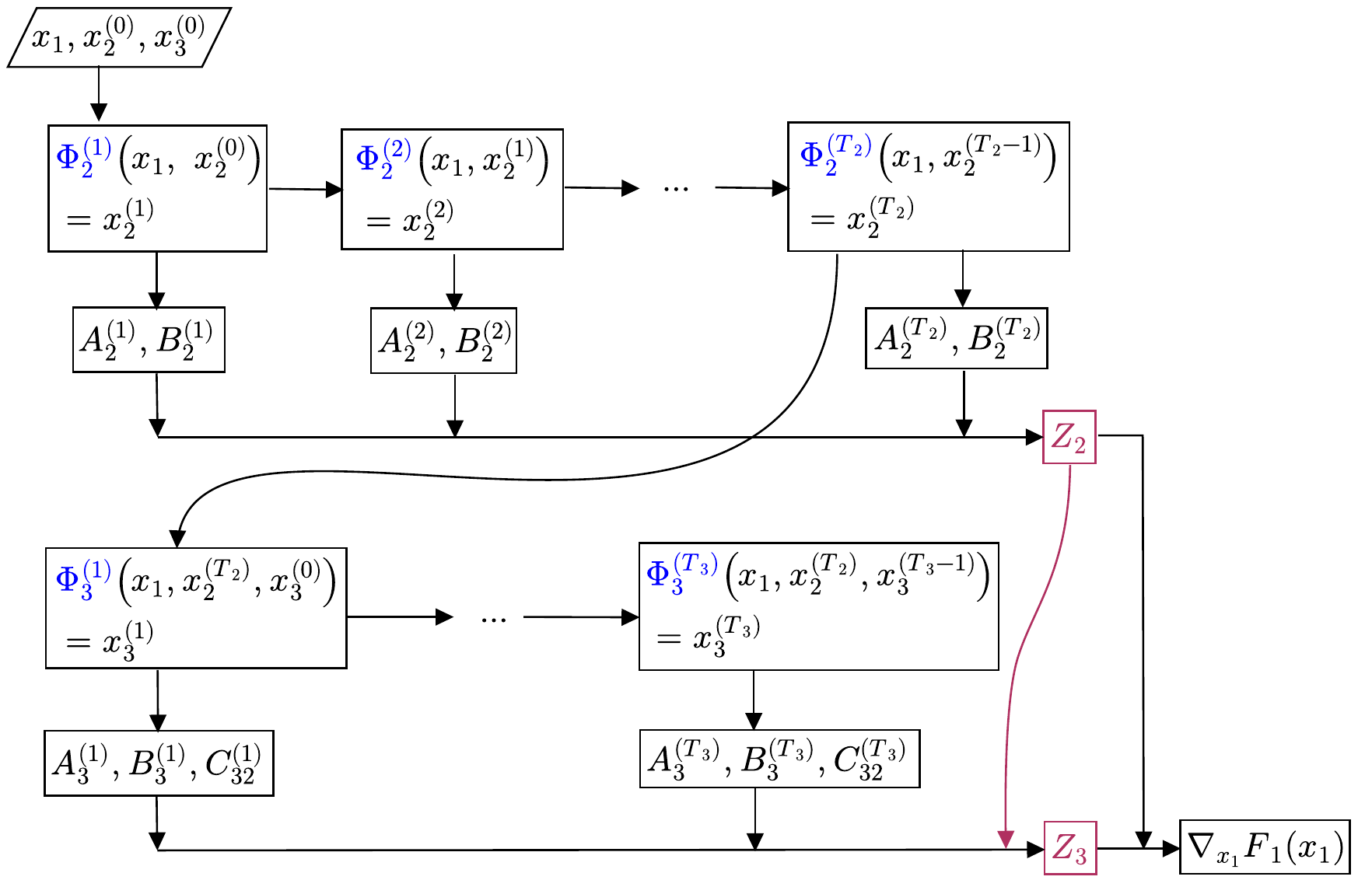}
\caption{Procedure when applying proposed algorithm to trilevel optimization problems.}
\label{fig: suppl flowchart}
\end{figure}

The gradient of $\tilde{F}_{1}(x_{1})$ can be calculated by applying
Algorithm \ref{algo: grad comp} to this problem.
The explicit formula of $\nabla_{x_{1}} \tilde{F}_{1}(x_{1})$ is given in Theorem \ref{thm: grad formula}
and hence we can compute it by using Algorithm \ref{algo: grad comp}.
To compute $\nabla_{x_{1}} \tilde{F}_{1}(x_{1})$ with Algorithm \ref{algo: grad comp},
the computation of $\nabla_{x_{2}} \tilde{F}_{2}(x_{1}, x_{2})$ is also required.
Its explicit formula is also given in Theorem \ref{thm: grad formula}
and hence we can compute it by using Algorithm \ref{algo: grad comp}.
To compute $\nabla_{x_{2}} \tilde{F}_{2}(x_{1}, x_{2})$ with Algorithm \ref{algo: grad comp},
the computation of $\nabla_{x_{3}} \tilde{F}_{3}(x_{1}, x_{2}, x_{3})$ is also required.
This computation is easy (recall $\tilde{F}_{3} = f_{3}$).
With these results, we can apply a gradient-based method
(e.g., the projected gradient method)
to minimize $\tilde{F}_{1}(x_{1})$ using $\nabla \tilde{F}_{1}$.
This procedure of gradient computation is depicted in Figure \ref{fig: suppl flowchart}.

\section{Complete numerical results}
\label{sect: experiments detail}

To validate the effectiveness of our proposed method,
we conducted numerical experiments on an artificial problem and a hyperparameter optimization problem arising from real data.
In our numerical experiments,
we implemented all codes with Python 3.9.2 and JAX 0.2.10 for automatic differentiation
and executed them on a computer with 12 cores of Intel Core i7-7800X CPU 3.50 GHz, 64 GB RAM, Ubuntu OS 20.04.2 LTS.

\subsection{Convergence to optimal solution}

We solved the following trilevel optimization problem by using Algorithm \ref{algo: grad comp} to evaluate the performance:
\begin{equation}
\begin{aligned}
\min_{x_{1} \in \bbR^{2}}{}&
	f_{1}(x_{1}, x_{2}^{\ast}, x_{3}^{\ast}) = \|x_{3}^{\ast} - x_{1}\|_{2}^{2} + \|x_{1}\|_{2}^{2} \st\\
&	\begin{aligned}
	x_{2}^{\ast} \in \argmin_{x_{2} \in \bbR^{2}}{}&
		f_{2}(x_{1}, x_{2}, x_{3}^{\ast}) = \|x_{2} - x_{1}\|_{2}^{2} \st\\
	&	x_{3}^{\ast} \in \argmin_{x_{3} \in \bbR^{2}} f_{3}(x_{1}, x_{2}, x_{3}) = \|x_{3} - x_{2}\|_{2}^{2}.
	\end{aligned}
\end{aligned}
\label{eq: suppl artificial instance}
\end{equation}
Clearly, the optimal solution for this problem is $x_{1} = x_{2} = x_{3} = (0, 0)^{\top}$.

\begin{figure}
\centering
\subfigure[Objective values in early iterations]{
\includegraphics[width = 0.49\textwidth]{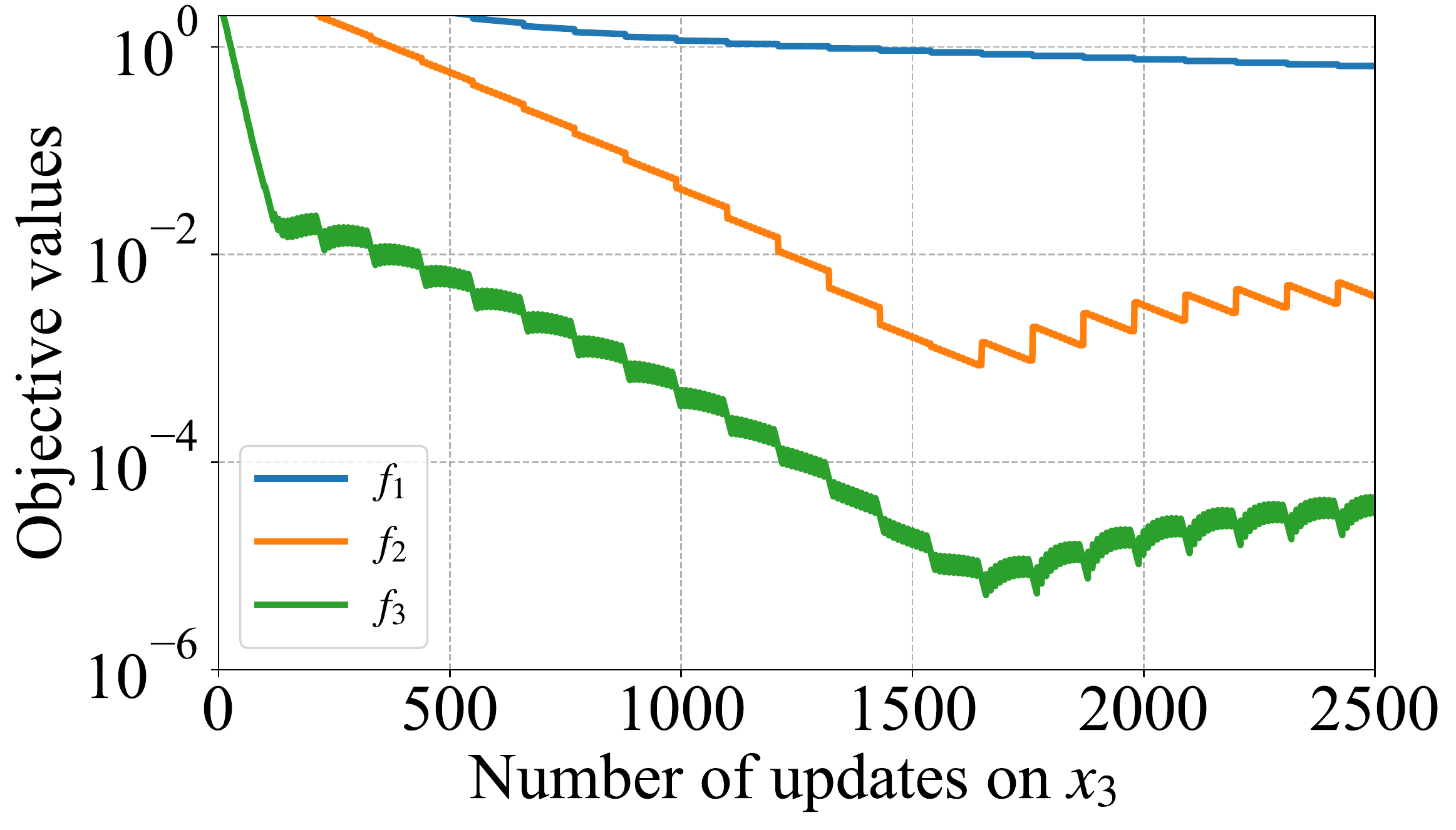}}
\subfigure[Objective values in later iterations]{
\includegraphics[width = 0.49\textwidth]{fig_3pts_obj/10_10_20000.pdf}%
\label{fig: suppl obj-10-10-20000}}
\caption{Performance of Algorithm \ref{algo: grad comp} for
	Problem \eqref{eq: suppl artificial instance} ($T_{2} = 10$, $T_{3} = 10$).}
\label{fig: suppl obj-10-10}
\end{figure}

\begin{figure}
\centering
\subfigure[Objective values in early iterations]{
\includegraphics[width = 0.49\textwidth]{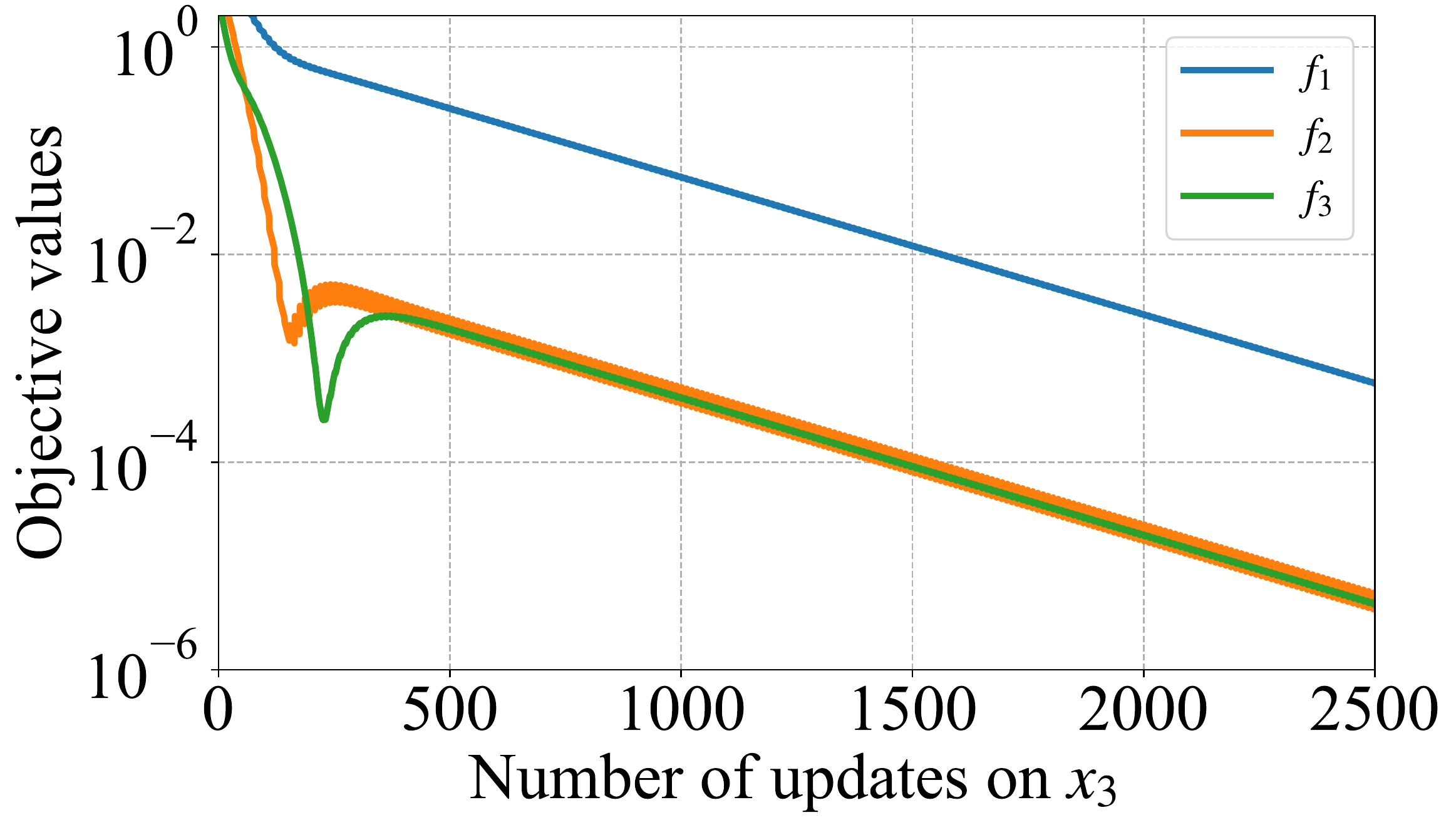}}
\subfigure[Objective values in later iterations]{
\includegraphics[width = 0.49\textwidth]{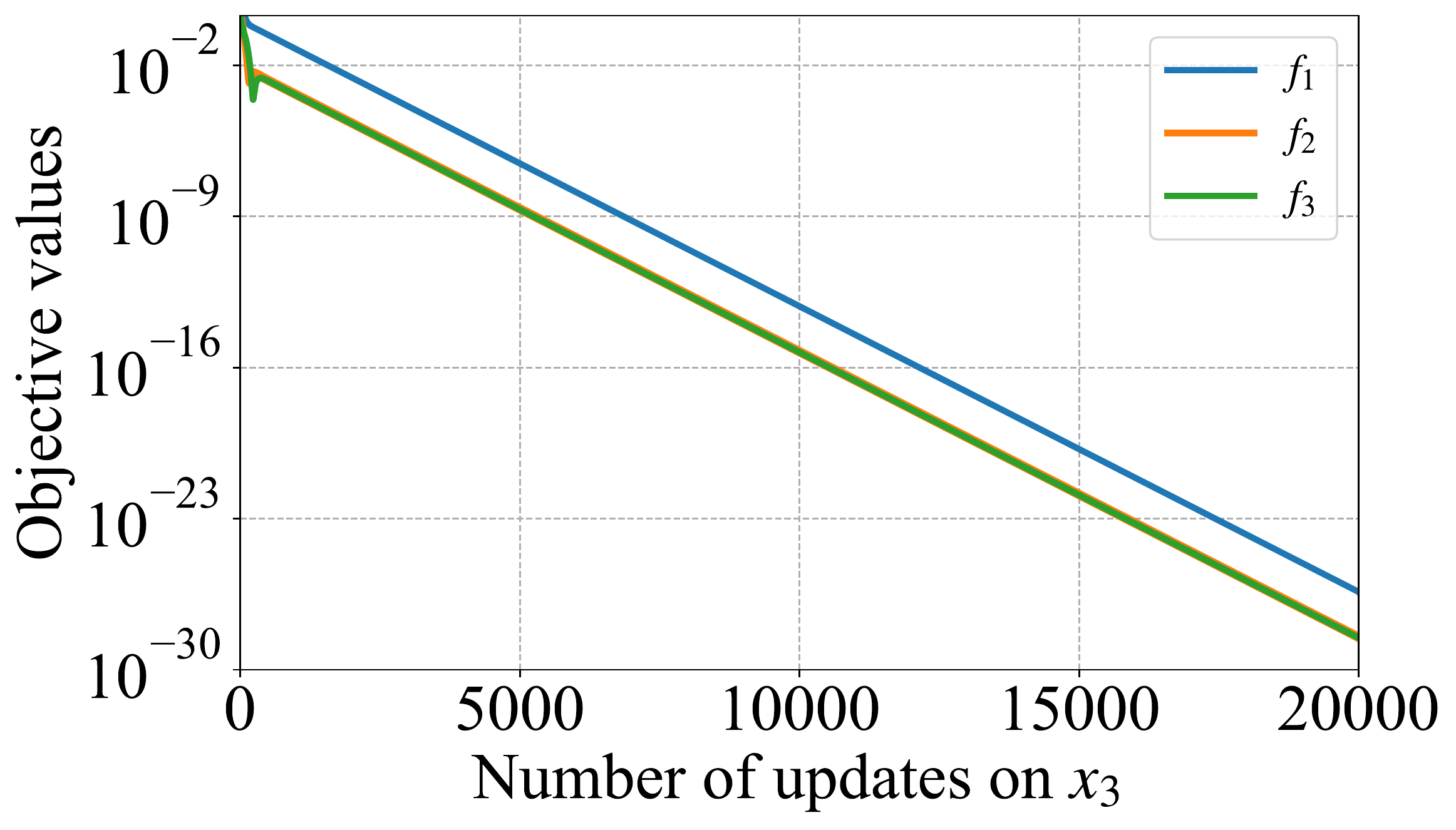}}
\caption{Performance of Algorithm \ref{algo: grad comp} for
	Problem \eqref{eq: suppl artificial instance} ($T_{2} = 10$, $T_{3} = 1$).}
\label{fig: suppl obj-10-1}
\end{figure}

\begin{figure}
\centering
\subfigure[Objective values in early iterations]{
\includegraphics[width = 0.49\textwidth]{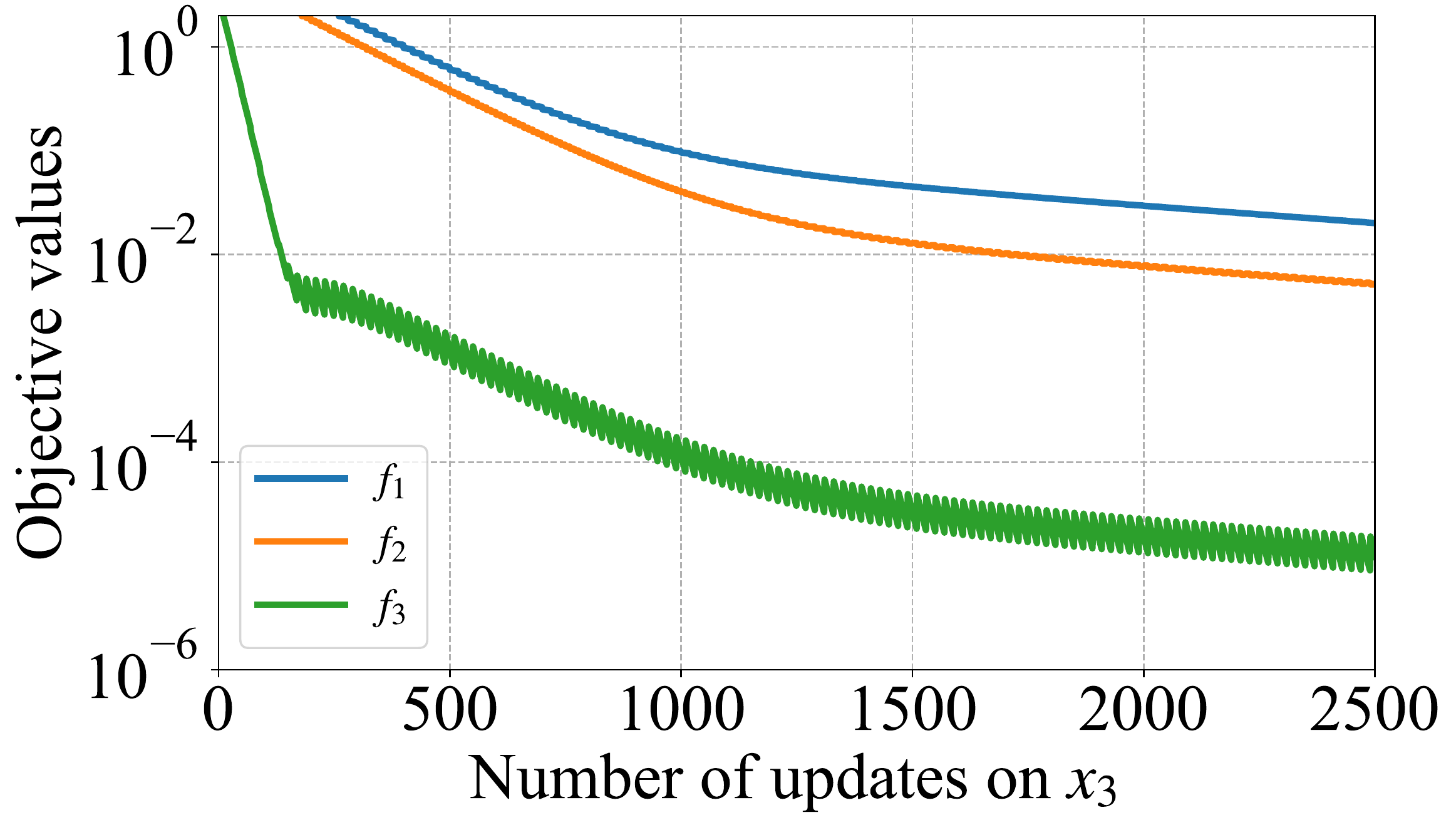}}
\subfigure[Objective values in later iterations]{
\includegraphics[width = 0.49\textwidth]{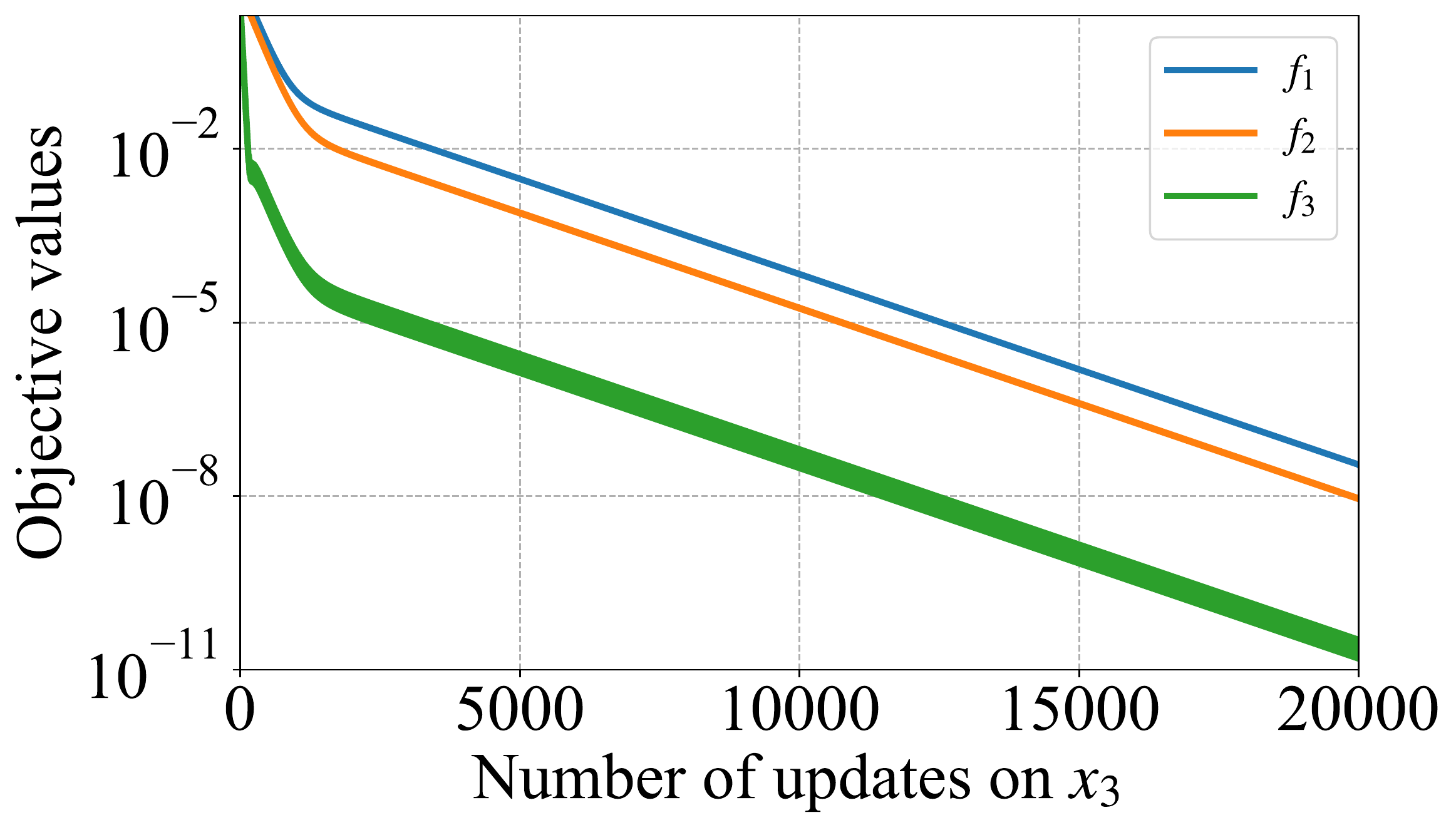}}
\caption{Performance of Algorithm \ref{algo: grad comp}
for Problem \eqref{eq: suppl artificial instance} ($T_{2} = 1$, $T_{3} = 10$).}
\label{fig: suppl obj-1-10}
\end{figure}

\begin{figure}
\centering
\subfigure[Objective values in early iterations]{
\includegraphics[width = 0.49\textwidth]{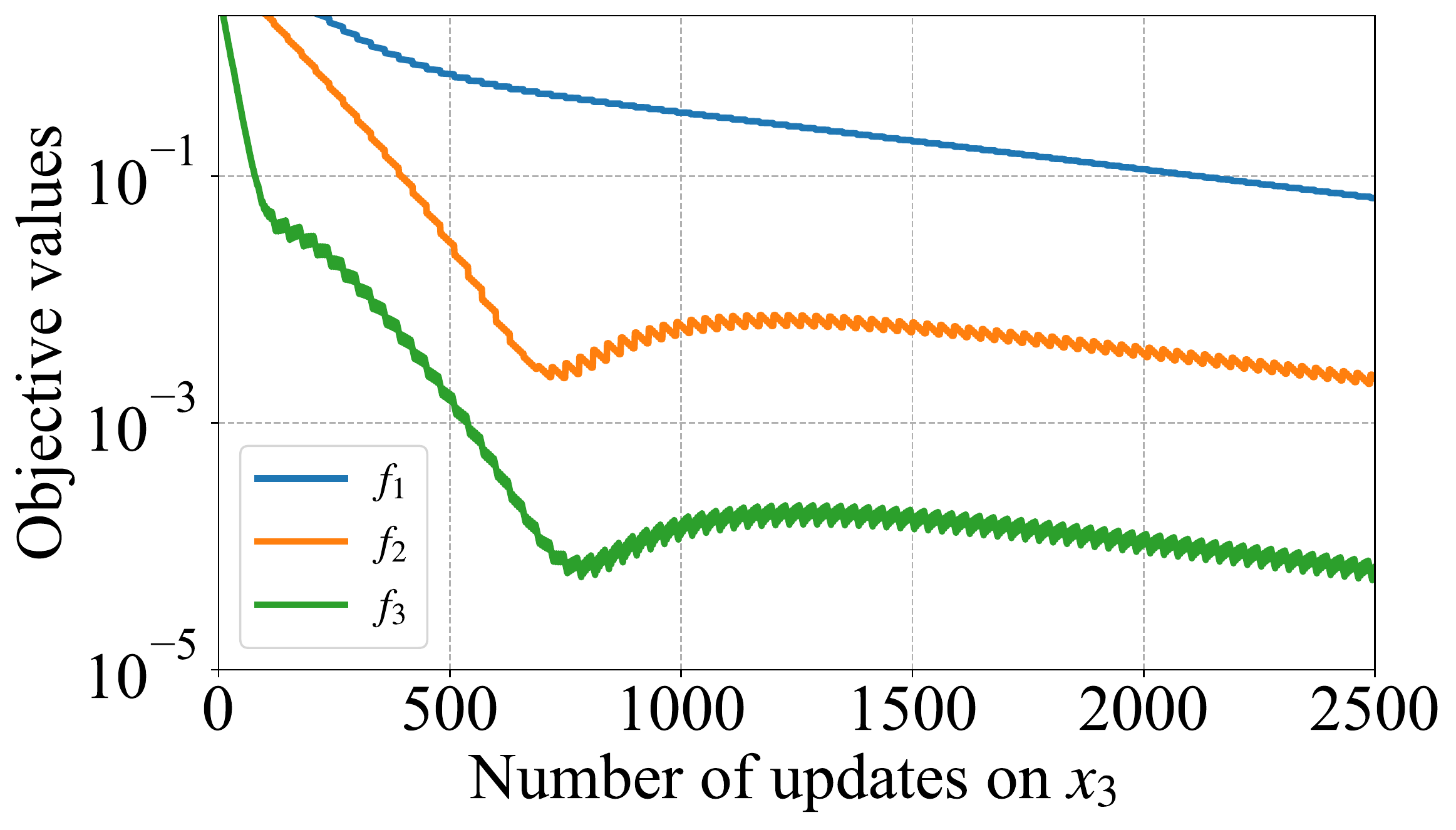}}
\subfigure[Objective values in later iterations]{
\includegraphics[width = 0.49\textwidth]{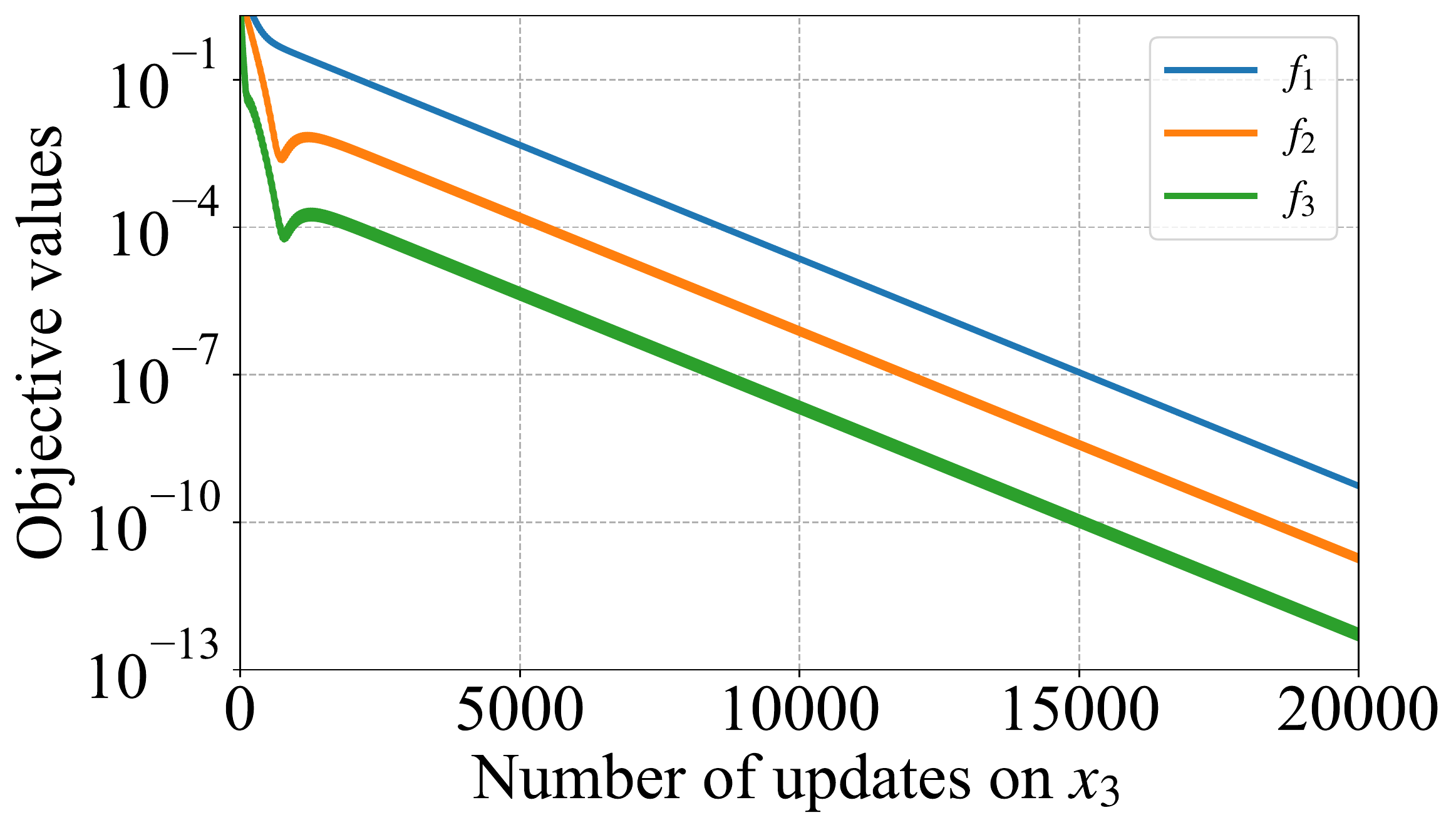}} 
\caption{Performance of Algorithm \ref{algo: grad comp}
	for Problem \eqref{eq: suppl artificial instance} ($T_{2} = 5$, $T_{3} = 5$).}
\label{fig: suppl obj-5-5}
\end{figure}

\begin{figure}
\centering
\subfigure[Objective values in early iterations]{
\includegraphics[width = 0.49\textwidth]{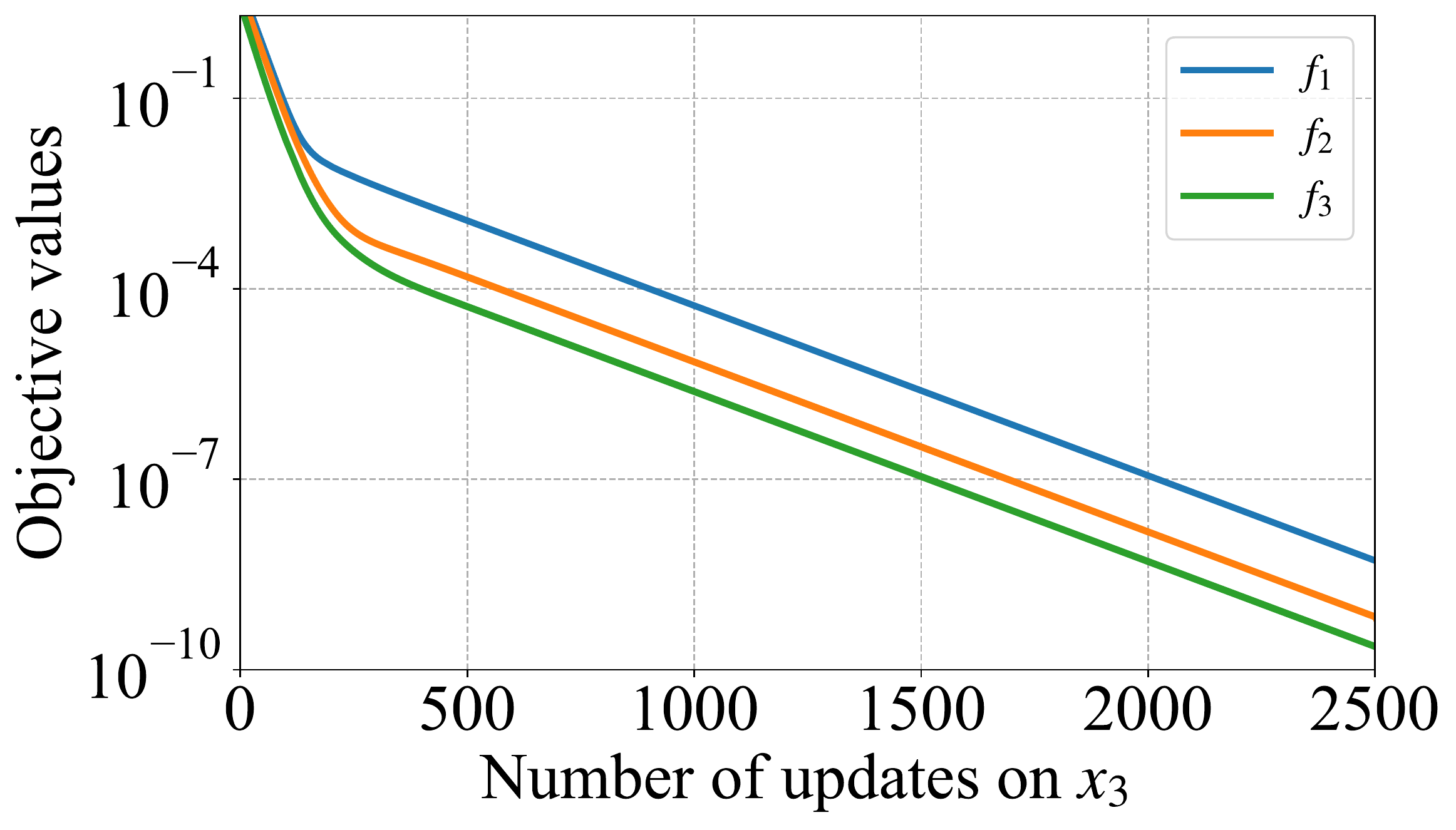}}
\subfigure[Objective values in later iterations]{
\includegraphics[width = 0.49\textwidth]{fig_3pts_obj/1_1_20000.pdf}%
\label{fig: suppl obj-1-1-20000}}
\caption{Performance of Algorithm \ref{algo: grad comp}
	for Problem \eqref{eq: suppl artificial instance} ($T_{2} = 1$, $T_{3} = 1$).}
\label{fig: suppl obj-1-1}
\end{figure}

\begin{figure}
\centering
\subfigure[$T_{2} = 10$, $T_{3} = 10$]{%
\includegraphics[width = 0.49\textwidth]{fig_3pts_traj/10_10_traj.pdf}}
\subfigure[$T_{2} = 10$, $T_{3} = 1$]{%
\includegraphics[width = 0.49\textwidth]{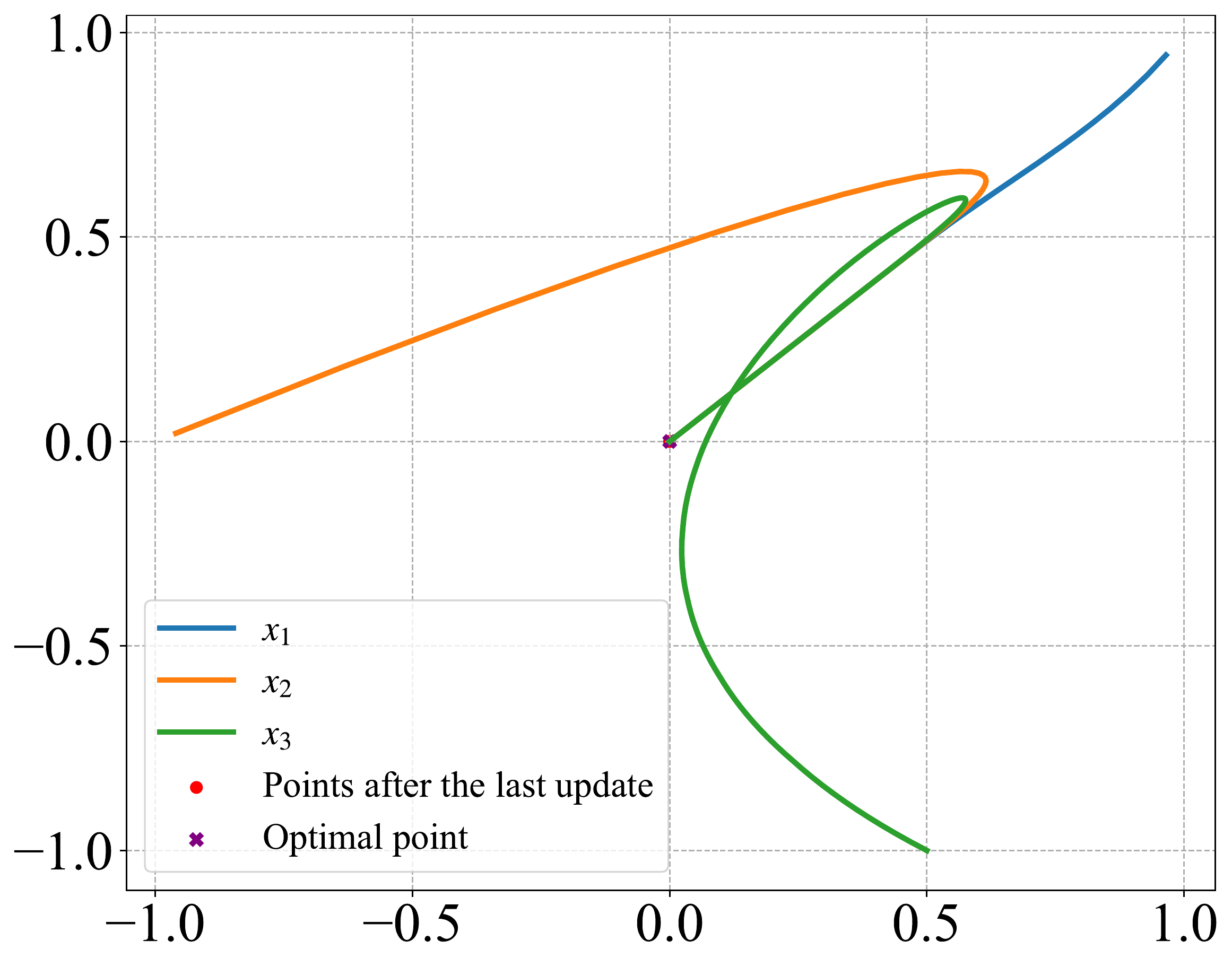}}
\subfigure[$T_{2} = 1$, $T_{3} = 10$]{%
\includegraphics[width = 0.49\textwidth]{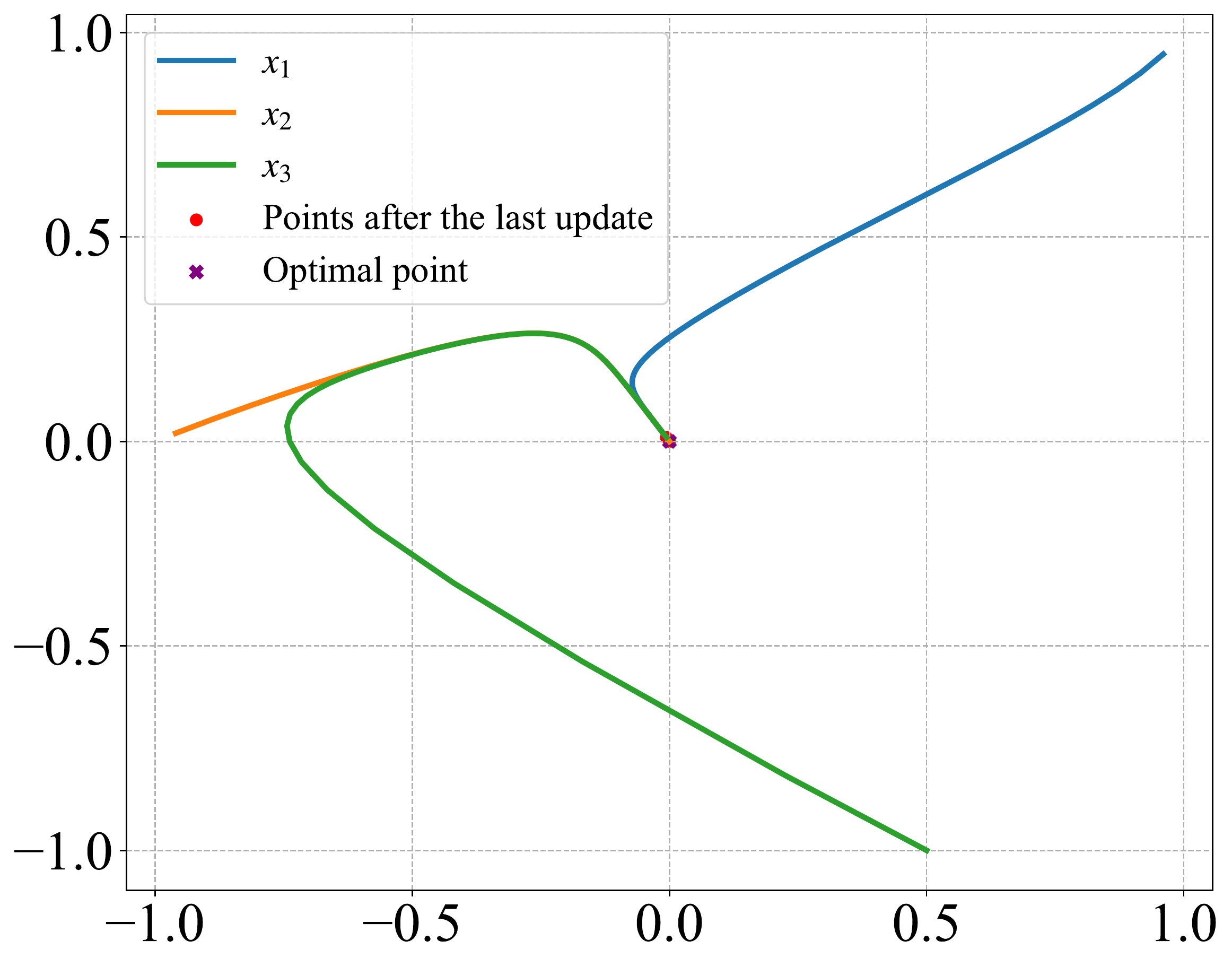}}
\subfigure[$T_{2} = 5$, $T_{3} = 5$]{%
\includegraphics[width = 0.49\textwidth]{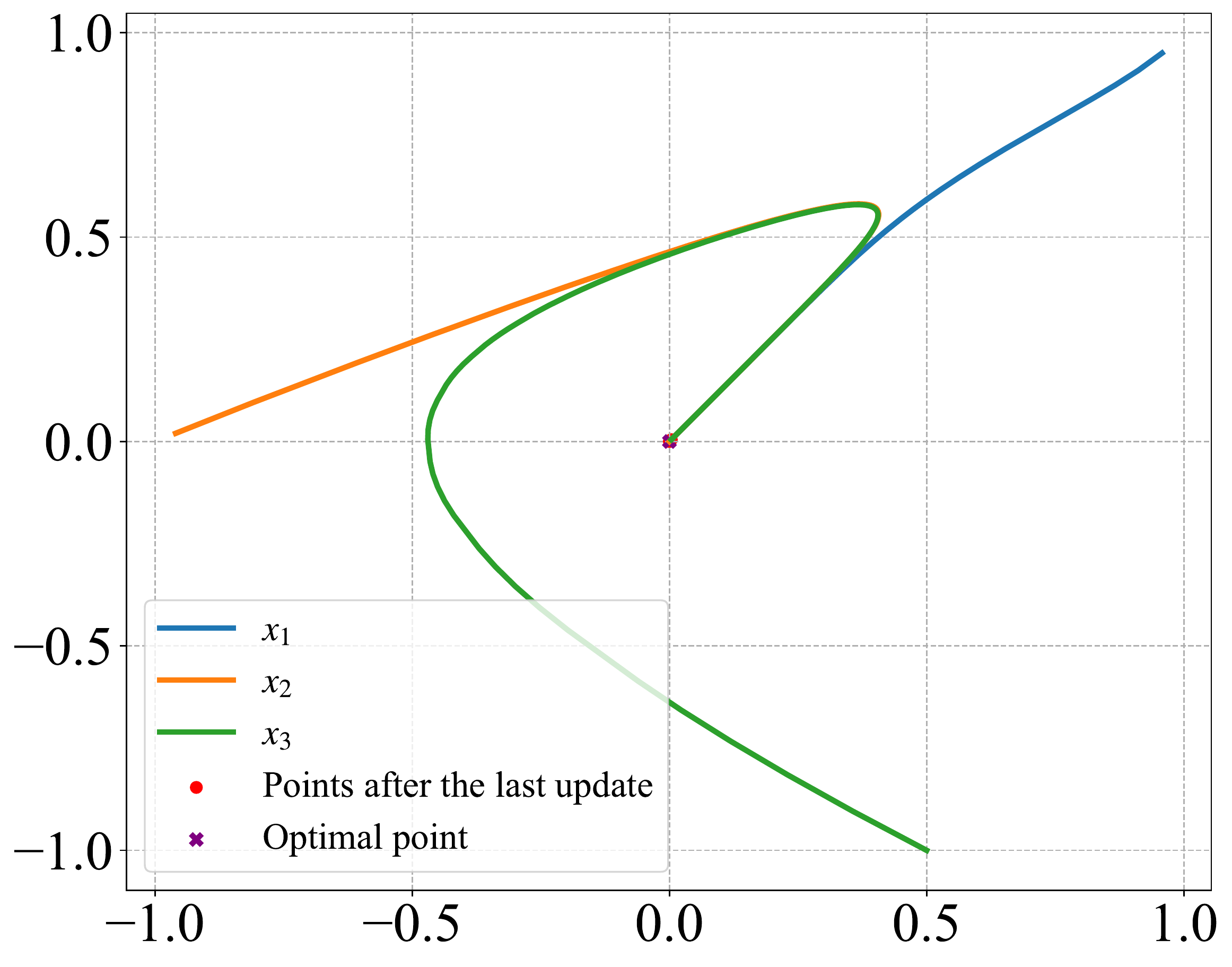}}
\subfigure[$T_{2} = 1$, $T_{3} = 1$]{%
\includegraphics[width = 0.49\textwidth]{fig_3pts_traj/1_1_traj.pdf}}
\caption{Trajectories of each variables.}
\label{fig: suppl traj}
\end{figure}

We solved \eqref{eq: artificial instance} with fixed constant step size and initialization but different $(T_{2}, T_{3})$.
For the iterative method in Algorithm \ref{algo: grad comp},
we employed the steepest descent method at all levels.
We show the transition of the value of the objective functions in
Figures \ref{fig: suppl obj-10-10}, \ref{fig: suppl obj-10-1}, \ref{fig: suppl obj-1-10}, \ref{fig: suppl obj-5-5},
and \ref{fig: suppl obj-1-1}
and the trajectories of each decision variables in Figure \ref{fig: suppl traj}.
We can confirm that the objective values $f_{1}$, $f_{2}$, and $f_{3}$ converge to the optimal value $0$ at all levels and
the gradient method outputs the optimal solution $(0, 0)^{\top}$
even if we use small iteration numbers $T_{2} = 1$ and $T_{3} = 1$ for the approximation problem \eqref{eq: approx}.
In Figures \ref{fig: suppl obj-10-10}, \ref{fig: suppl obj-10-1}, \ref{fig: suppl obj-1-10}, and \ref{fig: suppl obj-5-5},
there is oscillation of the objective values of lower levels,
as some decision variable moves away and closer to other decision variables.
In Figures \ref{fig: suppl obj-10-10}, \ref{fig: suppl obj-10-1}, and \ref{fig: suppl obj-5-5},
there is a temporary trade-off in the objective values between the levels.

\subsection{Comparison between Algorithm \ref{algo: grad comp} and an existing algorithm}
\label{sect: evolutionary}
We also conducted an experiments to compare the performance of Algorithm \ref{algo: grad comp}
with that of an existing method \cite{TKMO12}, which is based on evolutional strategy.
In this experiment, we used the same settings in Section \ref{exp: convergence}.

\begin{figure}
	\centering
	\subfigure[Transition of objective values (average of 100 times)]{
	\includegraphics[width = 0.49\textwidth]{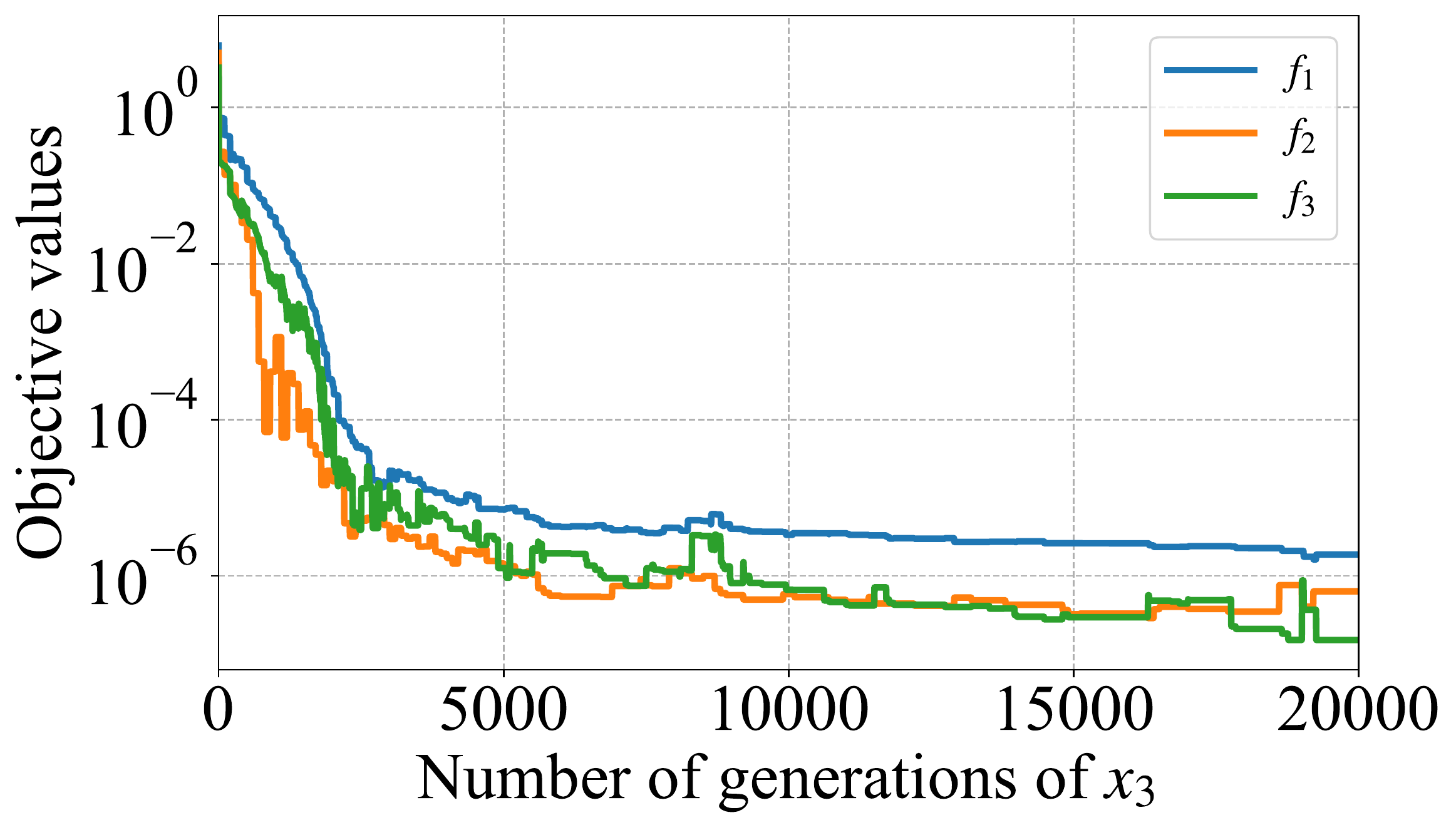}%
	\label{fig: suppl obj-evolutionary-001-obj}}
	\subfigure[Trajectory of each variables (one of 100 times)]{
	\includegraphics[width = 0.49\textwidth]{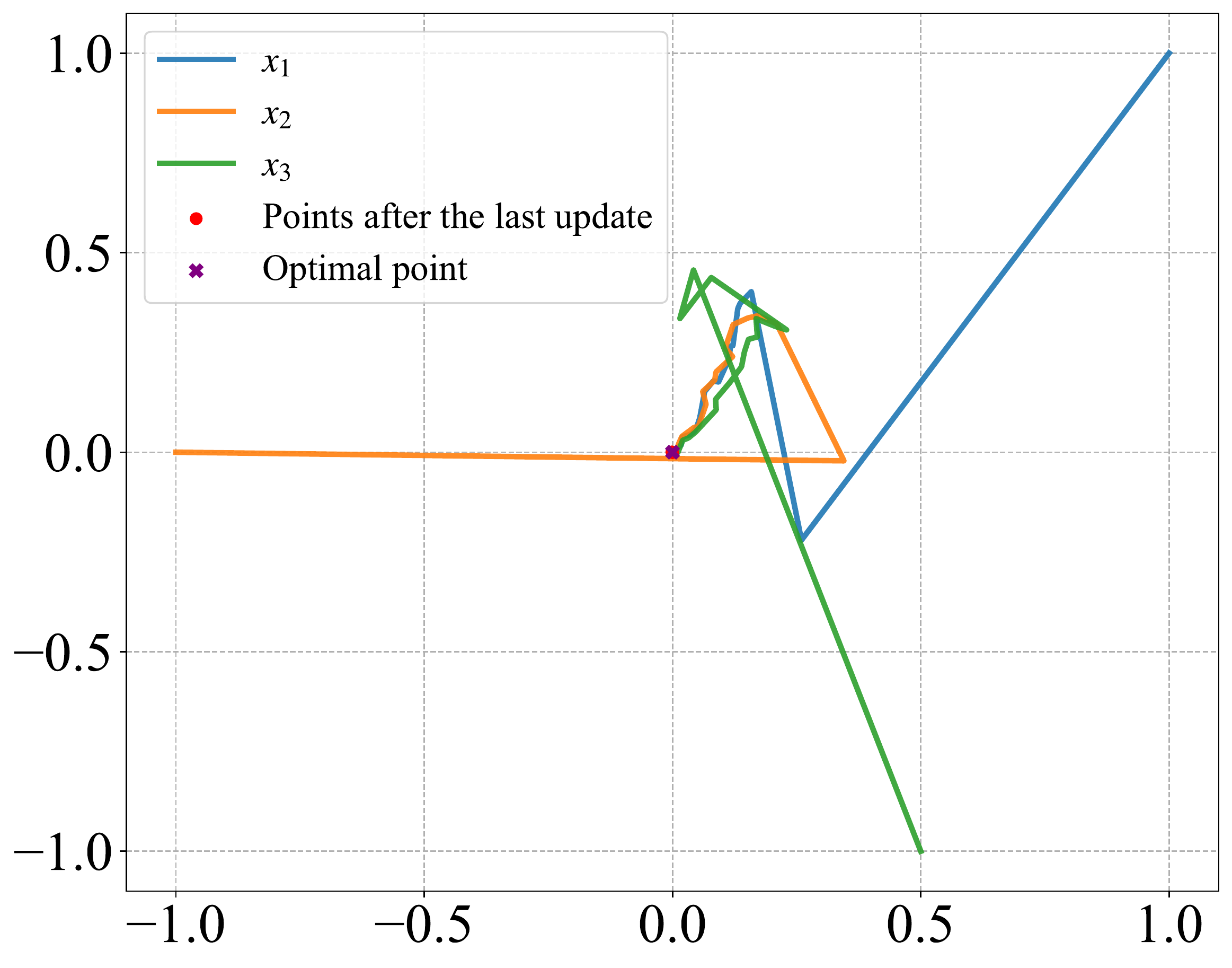}}
	\caption{Performance of the evolutional strategic algorithm 
		\cite{TKMO12} for Problem \eqref{eq: suppl artificial instance} $(\delta = 10^{-2})$.}
	\label{fig: suppl obj-evolutionary-001}
\end{figure}

\begin{figure}
	\centering
	\subfigure[Transition of objective values (average of 100 times)]{
	\includegraphics[width = 0.49\textwidth]{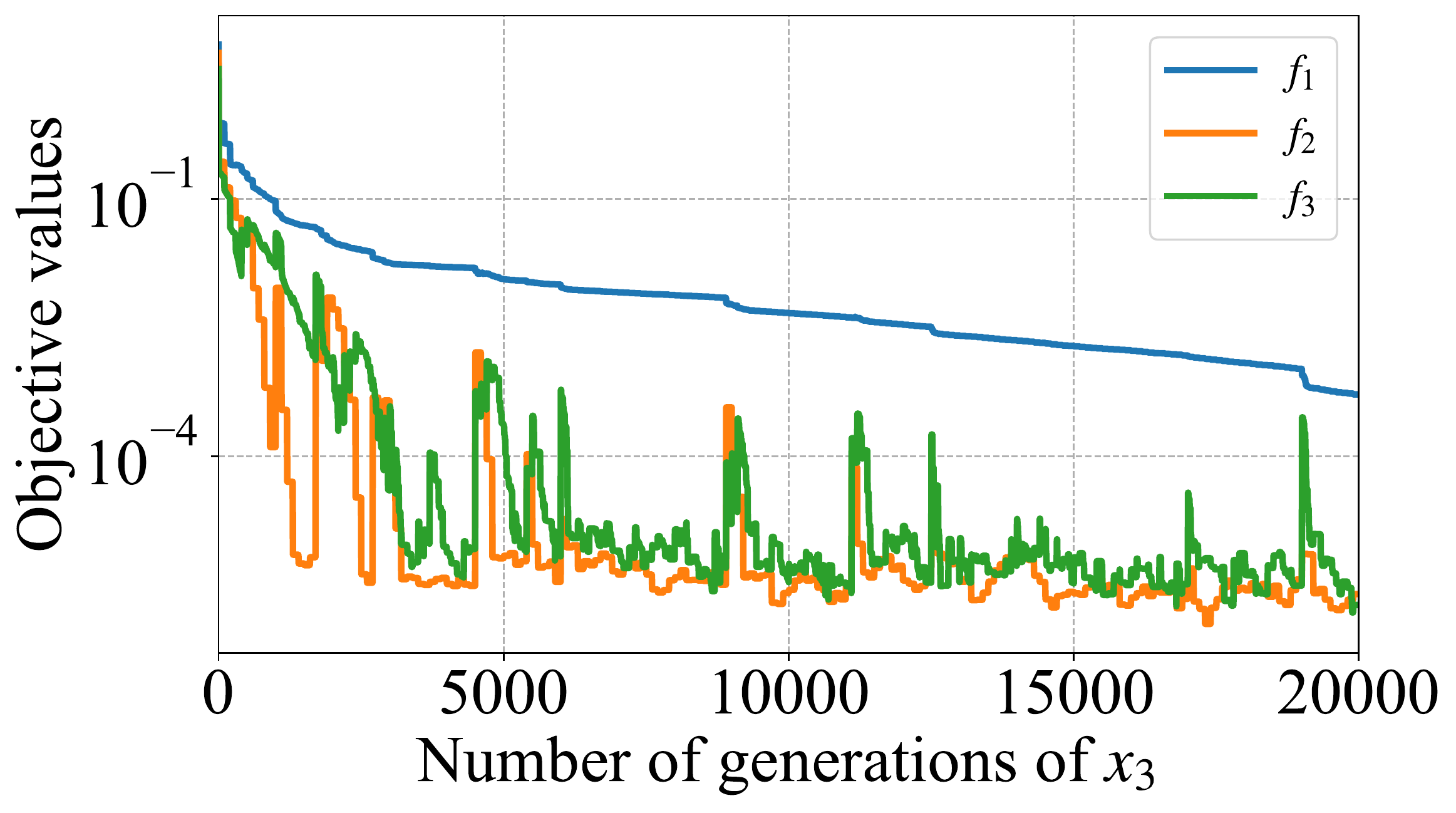}%
	\label{fig: suppl obj-evolutionary-00001-obj}}
	\subfigure[Trajectory of each variables (one of 100 times)]{
	\includegraphics[width = 0.49\textwidth]{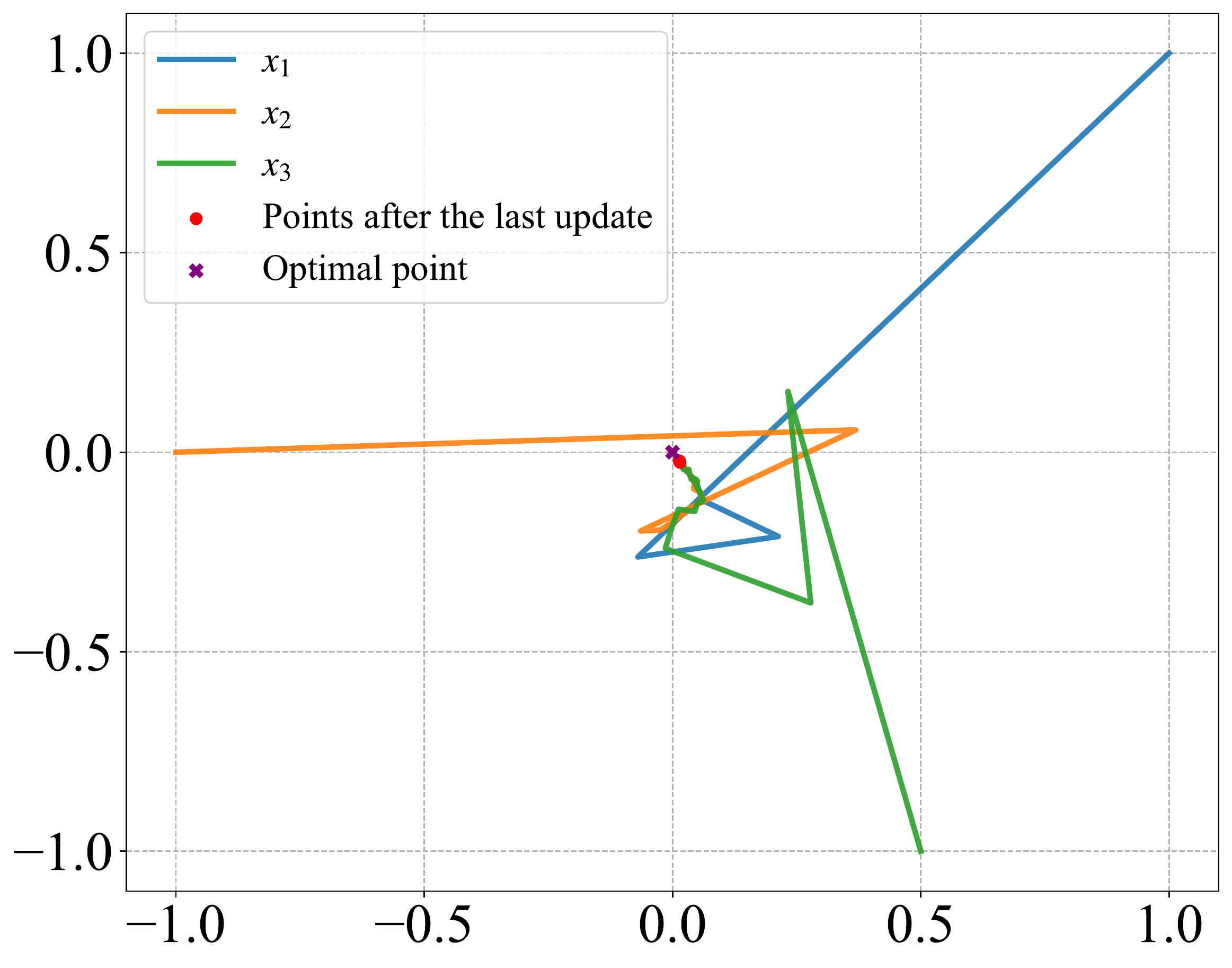}}
	\caption{Performance of the evolutional strategic algorithm 
		\cite{TKMO12} for Problem \eqref{eq: suppl artificial instance} $(\delta = 10^{-4})$.}
	\label{fig: suppl obj-evolutionary-00001}
\end{figure}

First, we applied the evolutional strategic algorithm \cite[Table 1 and Table 2]{TKMO12} to the problem \ref{eq: artificial instance}.
To clarify its difference between Algorithm \ref{algo: grad comp},
we picked the same initial solution as that of Section \ref{exp: convergence}.
In \cite[\texttt{EvolutionaryStrategy} in Table 2]{TKMO12}, the number of the generations of solution candidates
and the number of iterations of perturbations
were both set to be 10, the algorithm parameter $\delta$ was set to be $10^{-2}$ and $10^{-4}$,
and the solution candidates were randomly chosen from the standard normal distribution.
We conducted the experiment 100 times for each $\delta$, and
we show the transition of the average of the objective functions and
the trajectories of each decision variables in
Figures \ref{fig: suppl obj-evolutionary-001} and \ref{fig: suppl obj-evolutionary-00001}.

\begin{figure}
	\centering
	\includegraphics[width = 0.99\textwidth]{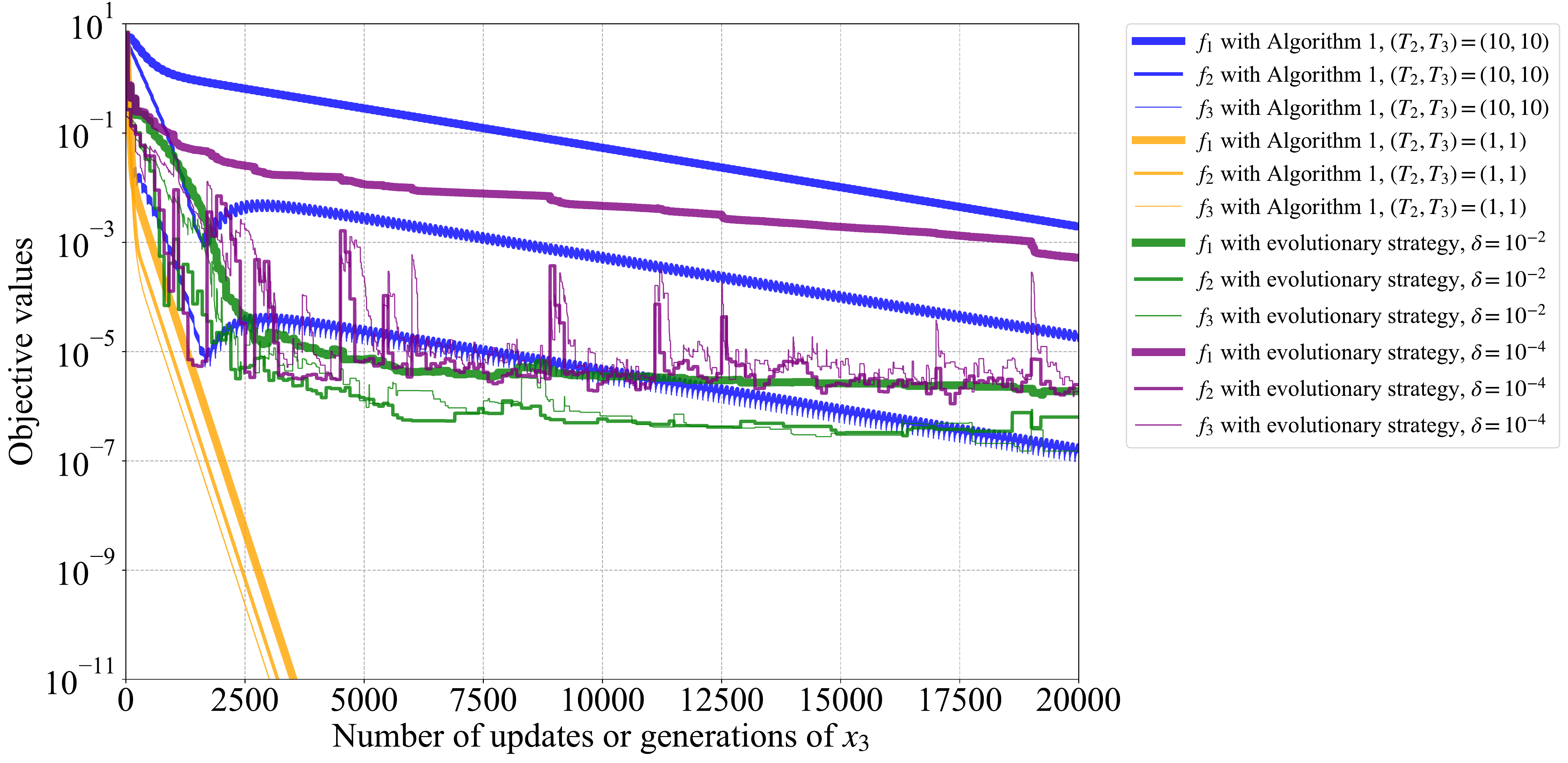}
	\caption{Comparison between the evolutional strategy \cite{TKMO12} and Algorithm \ref{algo: grad comp}
		for Problem \eqref{eq: suppl artificial instance}.
		Transition of the objective values is plotted.
		Note that the objective values with the evolutionary strategy are the average of 100 times.}
	\label{fig: suppl obj-compare}
\end{figure}

Next, we compared the performance of the existing method with that of Algorithm \ref{algo: grad comp}.
We show the transition of the objective values when using each algorithm in Figure \ref{fig: suppl obj-compare}
by overlaying Figures \ref{fig: suppl obj-10-10-20000}, \ref{fig: suppl obj-1-1-20000}, \ref{fig: suppl obj-evolutionary-001-obj}, and \ref{fig: suppl obj-evolutionary-00001-obj}.
In later iterations, we can confirm that the convergence of Algorithm \ref{algo: grad comp} with $(T_{2}, T_{3}) = (1, 1)$
is faster than the evolutional strategic algorithm.

\subsection{Application to hyperparameter optimization}
Next, we conducted experiments on real data to usefulness of our method for a machine learning problem.
We consider the problem of how to achieve a model robust to noise in input data and
formulate this problem by
a trilevel model by a model learner and an attacker.
The model learner decides the hyperparameter $\lambda$ to minimize the validation error,
while the attacker tries to poison training data so as to make the model less accurate.
This model is formulated as follows:
\begin{equation}
\begin{aligned}
\min_{\lambda}{}&
	\frac{1}{m} \|y_{\text{valid}} - X_{\text{valid}} \theta\|_{2}^{2} \st\\
&	\begin{aligned}
	P \in \argmax_{P'}{}&
		\frac{1}{n} \|y_{\text{train}} - (X_{\text{train}} + P') \theta\|_{2}^{2} - \frac{c}{n d}\|P'\|_{2}^{2} \st\\
	&	\theta \in \argmin_{\theta'} \frac{1}{n} \|y_{\text{train}} - (X_{\text{train}} + P') \theta'\|_{2}^{2}
		+ \exp(\lambda) \frac{\|\theta'\|_{1^{\ast}}}{d},
	\end{aligned}
\end{aligned}
\end{equation}
where $\theta$ denotes the parameter of the model,
$d$ denotes the dimension of $\theta$,
$n$ denotes the number of the training data $X_{\text{train}}$,
$m$ denotes the number of the validation data $X_{\text{val}}$,
$c$ denotes the penalty for the noise $P$,
and $\|\cdot\|_{1^{\ast}}$ is a smoothed $\ell_{1}$-norm,
which is a differentiable approximation of the $\ell_{1}$-norm.
Here, we use $\exp(\lambda)$ to express nonnegative penalty parameter
instead of the constraint $\lambda \ge 0$.

To validate the effectiveness of our proposed method,
we compared the results by the trilevel model with those of the following bilevel model:
\begin{equation}
\min_{\lambda} \frac{1}{m} \|y_{\text{valid}} - X_{\text{valid}} \theta\|_{2}^{2}
\st \theta \in \argmin_{\theta'}
\frac{1}{n} \|y_{\text{train}} - X_{\text{train}} \theta'\|_{2}^{2} + \exp(\lambda) \frac{\|\theta'\|_{1^{\ast}}}{d}.
\label{ho-2level}
\end{equation}
This model is equivalent to the trilevel model without the attacker's level.

We used Algorithm \ref{algo: grad comp} to compute the gradient of the objective function in
the trilevel and bilevel models with real datasets.
For the iterative method in Algorithm \ref{algo: grad comp},
we employed the steepest descent method at all levels.
We set $T_{2} = 30$ and $T_{3} = 3$ for the trilevel model and $T_{2} = 30$ for the bilevel model.
In each dataset, we used the same initialization and step sizes in the updates of $\lambda$ and $\theta$ in trilevel and bilevel models.
We compared these methods on the regression tasks with the following datasets:
the diabetes dataset \cite{data-diabetes},
the (red and white) wine quality datasets \cite{data-wine},
the Boston dataset \cite{data-Boston}.
For each dataset, we standardized each feature and the objective variable;
and randomly chose 40 rows as training data $(X_{\text{train}}, y_{\text{train}})$
other chose 100 rows as validation data $(X_{\text{valid}}, y_{\text{valid}})$,
and used the rest of the rows as test data.

\begin{figure}
	\centering
	\subfigure[$\sigma = 0.01$]{
	\includegraphics[width = 0.49\textwidth]{fig_diabetes/learning_curve_diabete_001.pdf}}
	\subfigure[$\sigma = 0.03$]{
	\includegraphics[width = 0.49\textwidth]{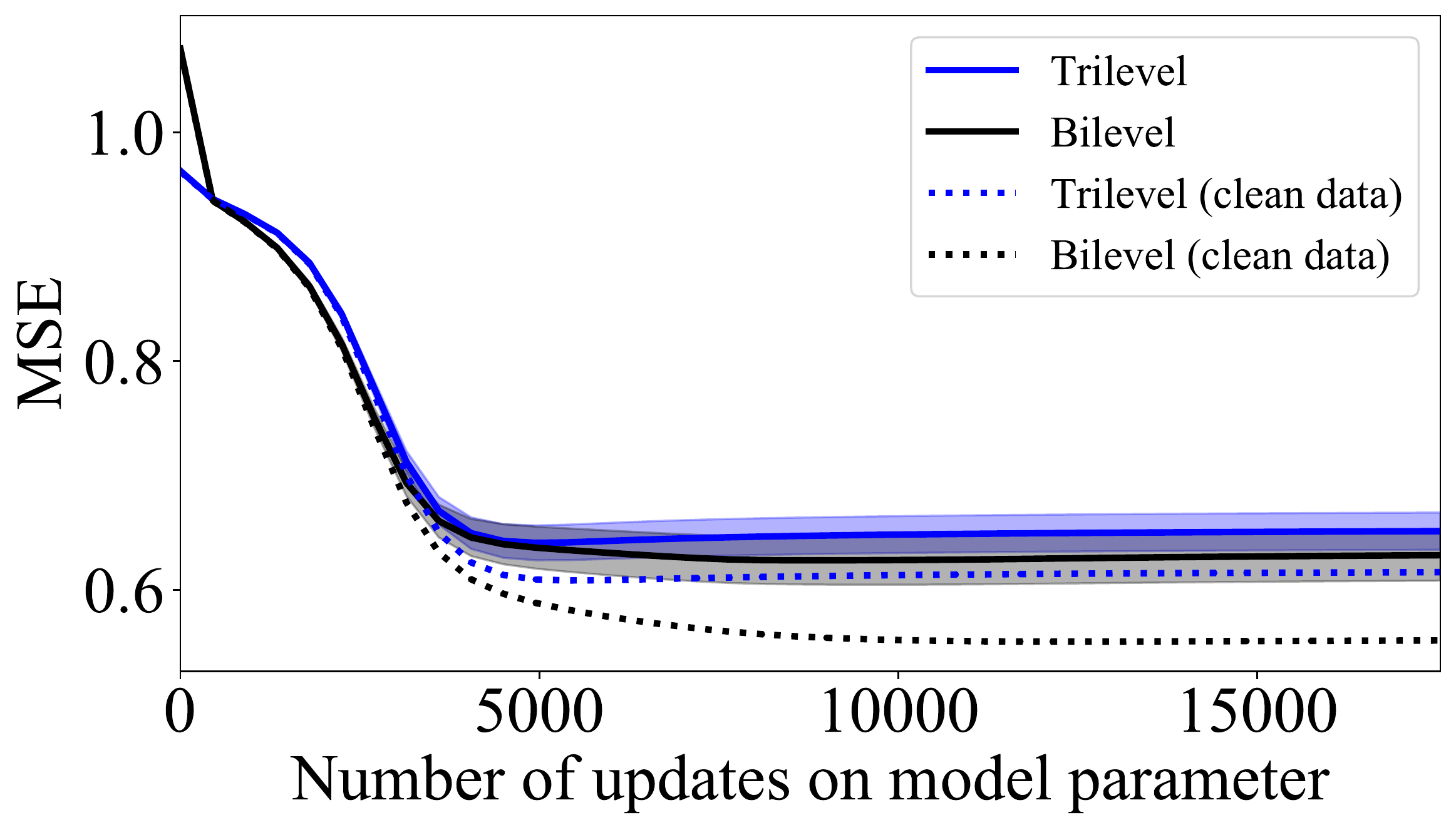}}
	\subfigure[$\sigma = 0.05$]{
	\includegraphics[width = 0.49\textwidth]{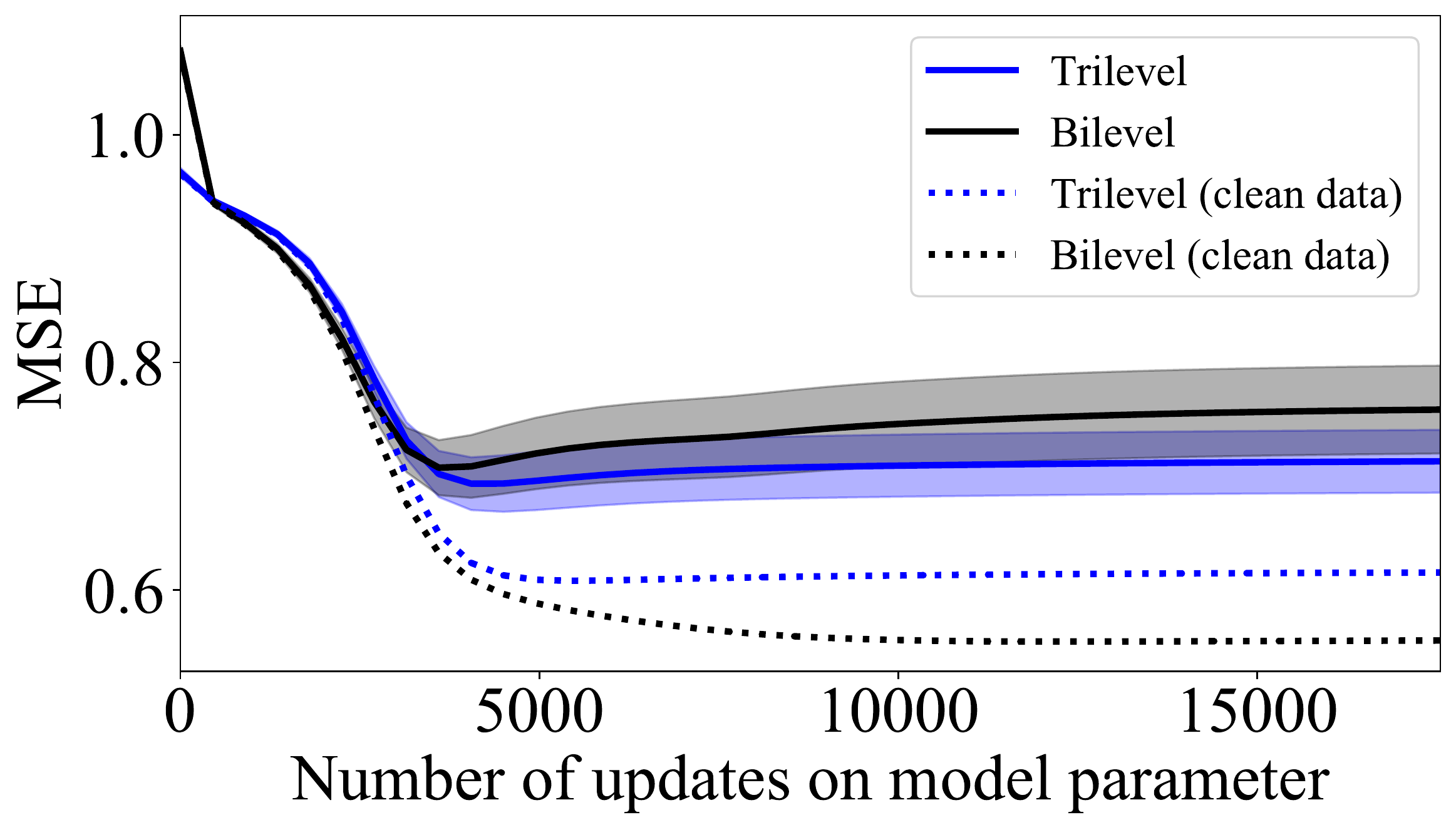}}
	\subfigure[$\sigma = 0.08$]{
	\includegraphics[width = 0.49\textwidth]{fig_diabetes/learning_curve_diabete_008.pdf}}
	\caption{MSE of test data with Gaussian noise with a standard deviation of $\sigma$ on diabetes dataset.}
	\label{fig: suppl iter-diabetes}
\end{figure}

\begin{figure}
	\centering
	\subfigure[$\sigma = 0.01$]{
	\includegraphics[width = 0.49\textwidth]{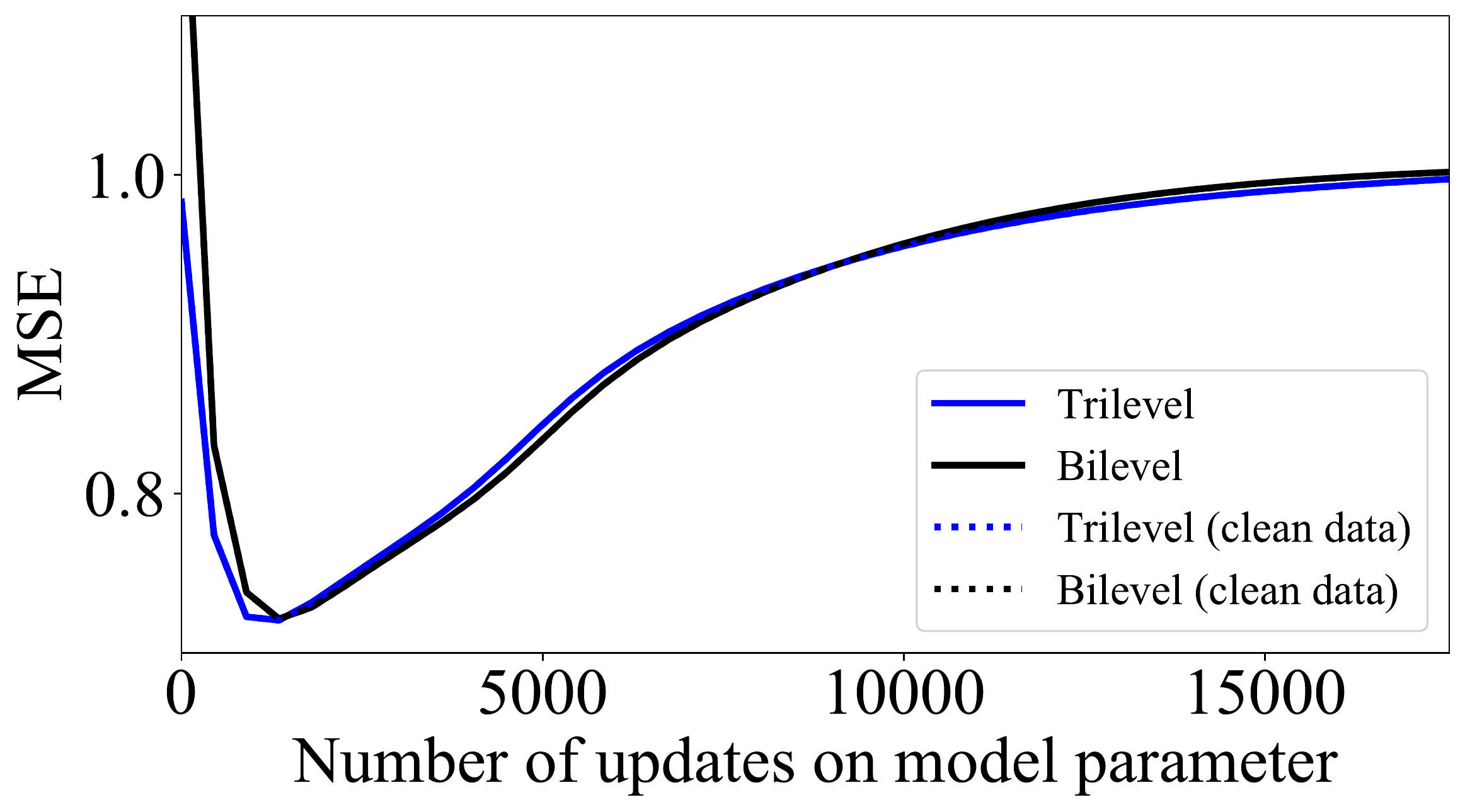}}
	\subfigure[$\sigma = 0.03$]{
	\includegraphics[width = 0.49\textwidth]{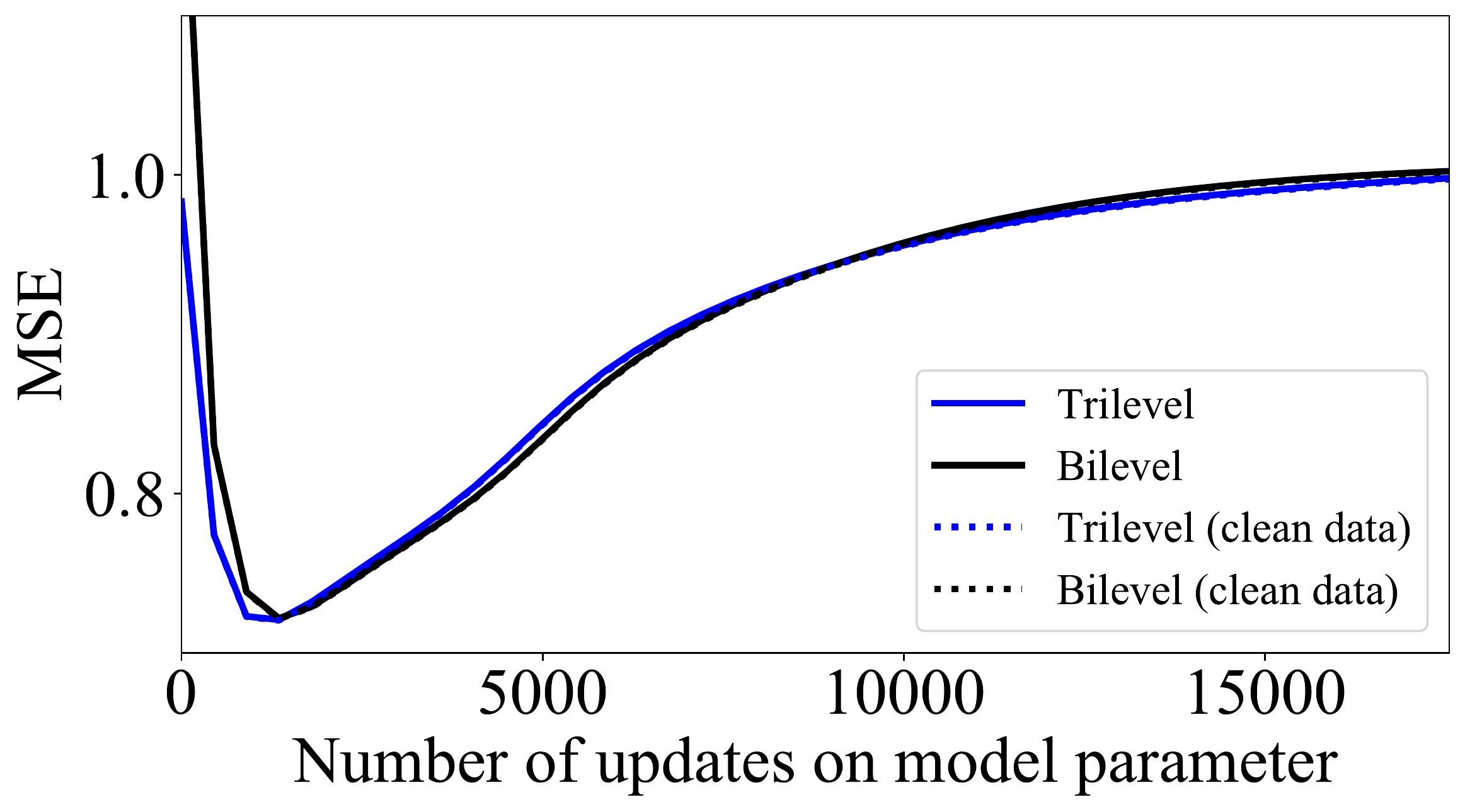}}
	\subfigure[$\sigma = 0.05$]{
	\includegraphics[width = 0.49\textwidth]{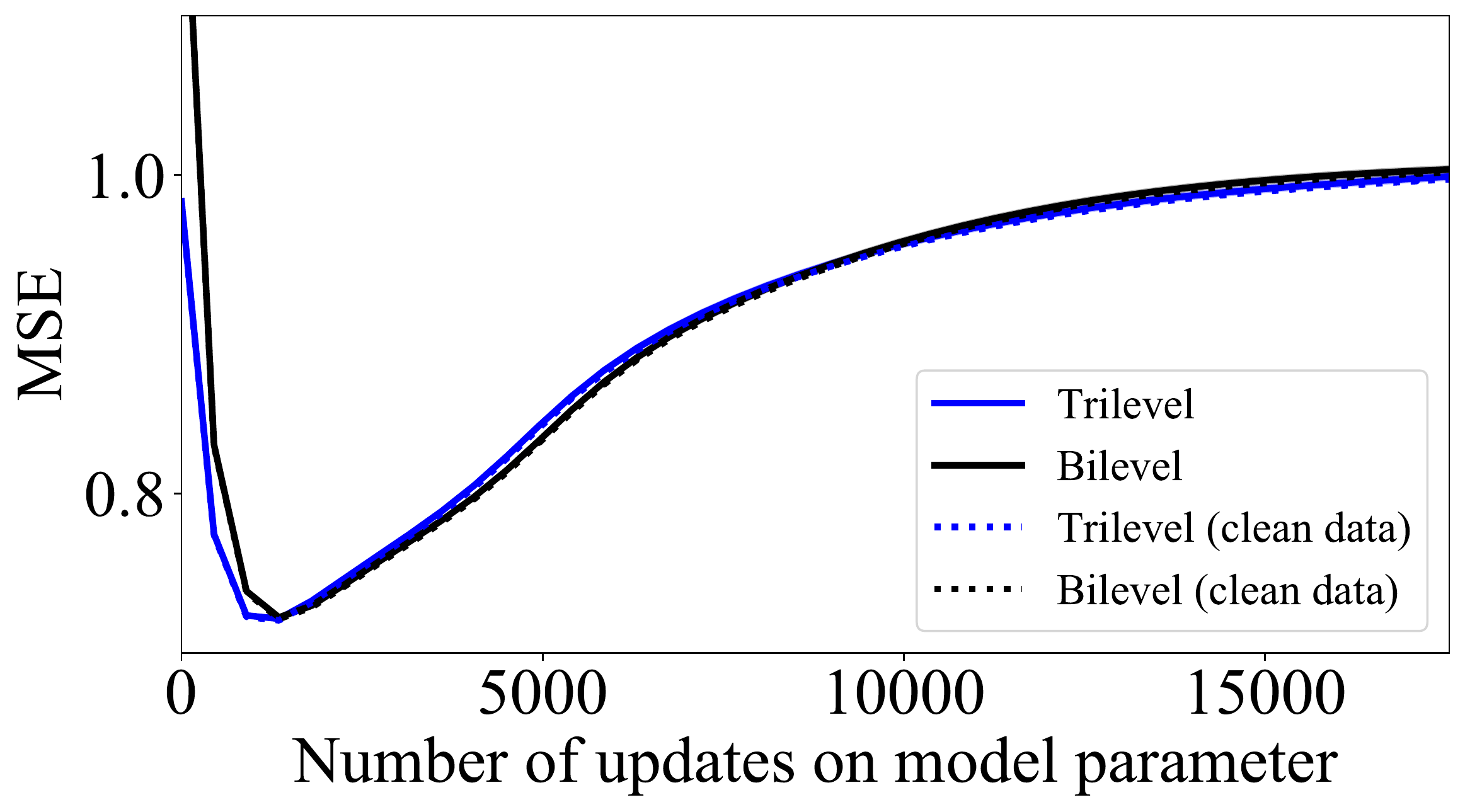}}
	\subfigure[$\sigma = 0.08$]{
	\includegraphics[width = 0.49\textwidth]{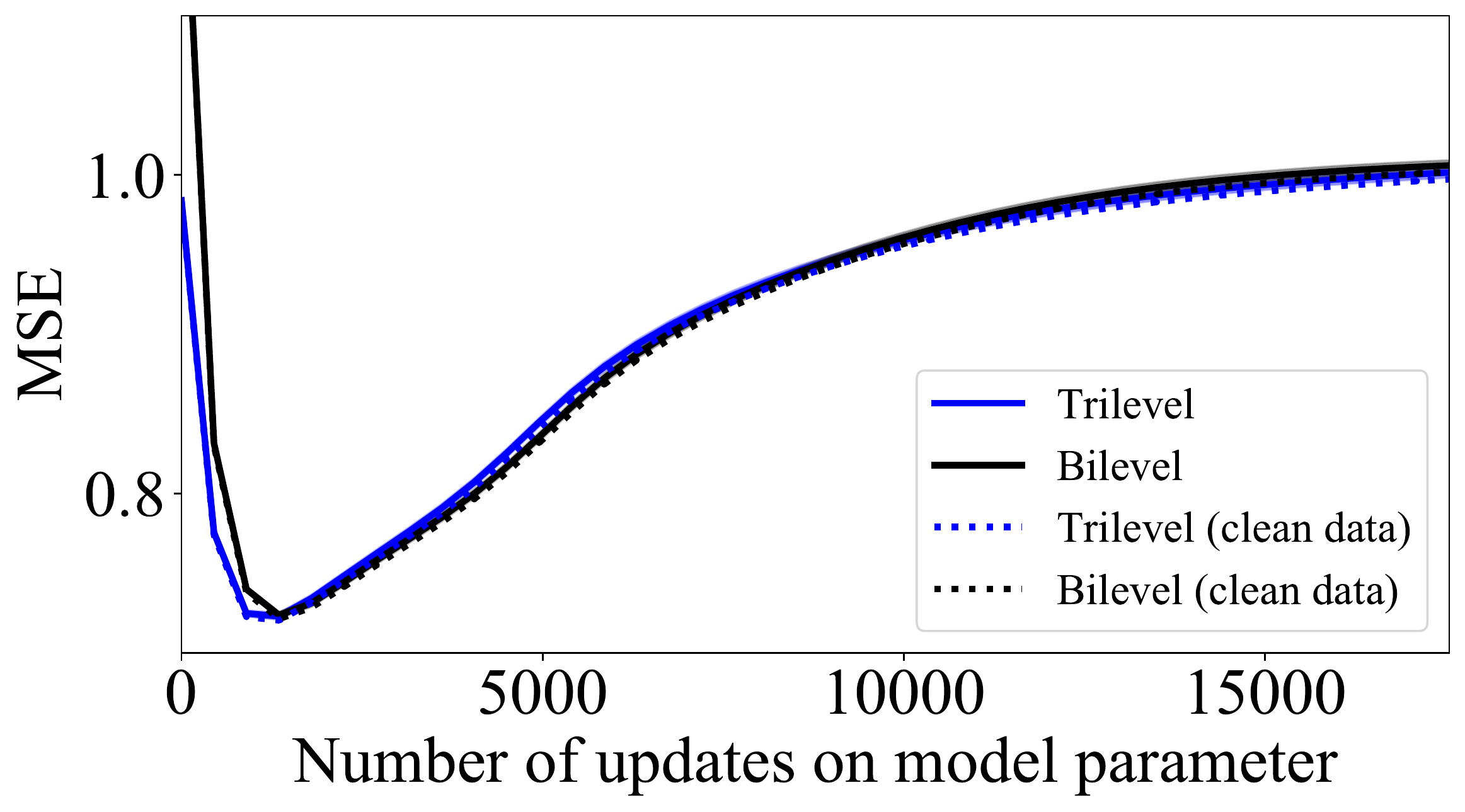}}
	\caption{MSE of test data with Gaussian noise with a standard deviation of $\sigma$ on red wine quality dataset.}
	\label{fig: suppl iter-wine_red}
\end{figure}

\begin{figure}
	\centering
	\subfigure[$\sigma = 0.01$]{
	\includegraphics[width = 0.49\textwidth]{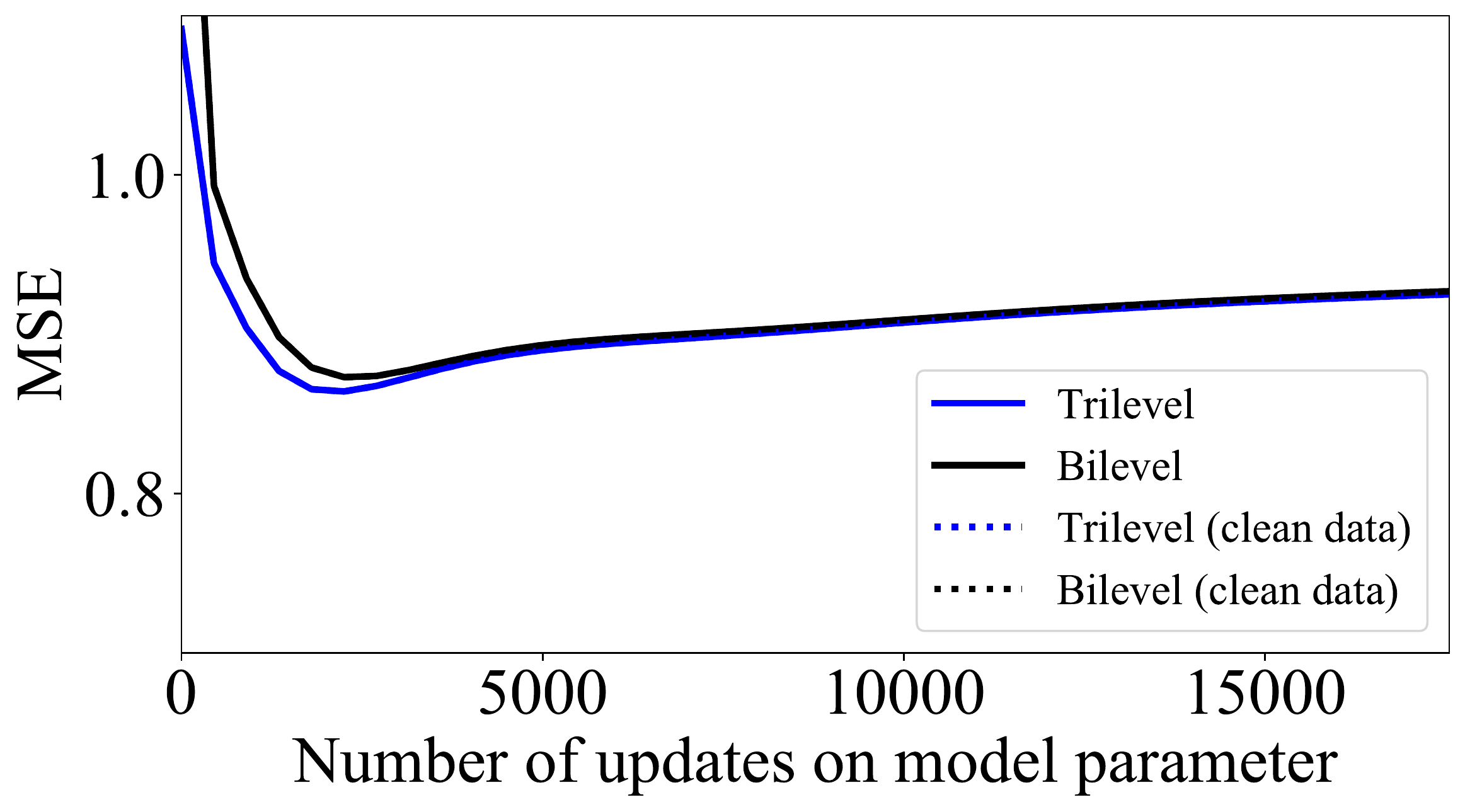}}
	\subfigure[$\sigma = 0.03$]{
	\includegraphics[width = 0.49\textwidth]{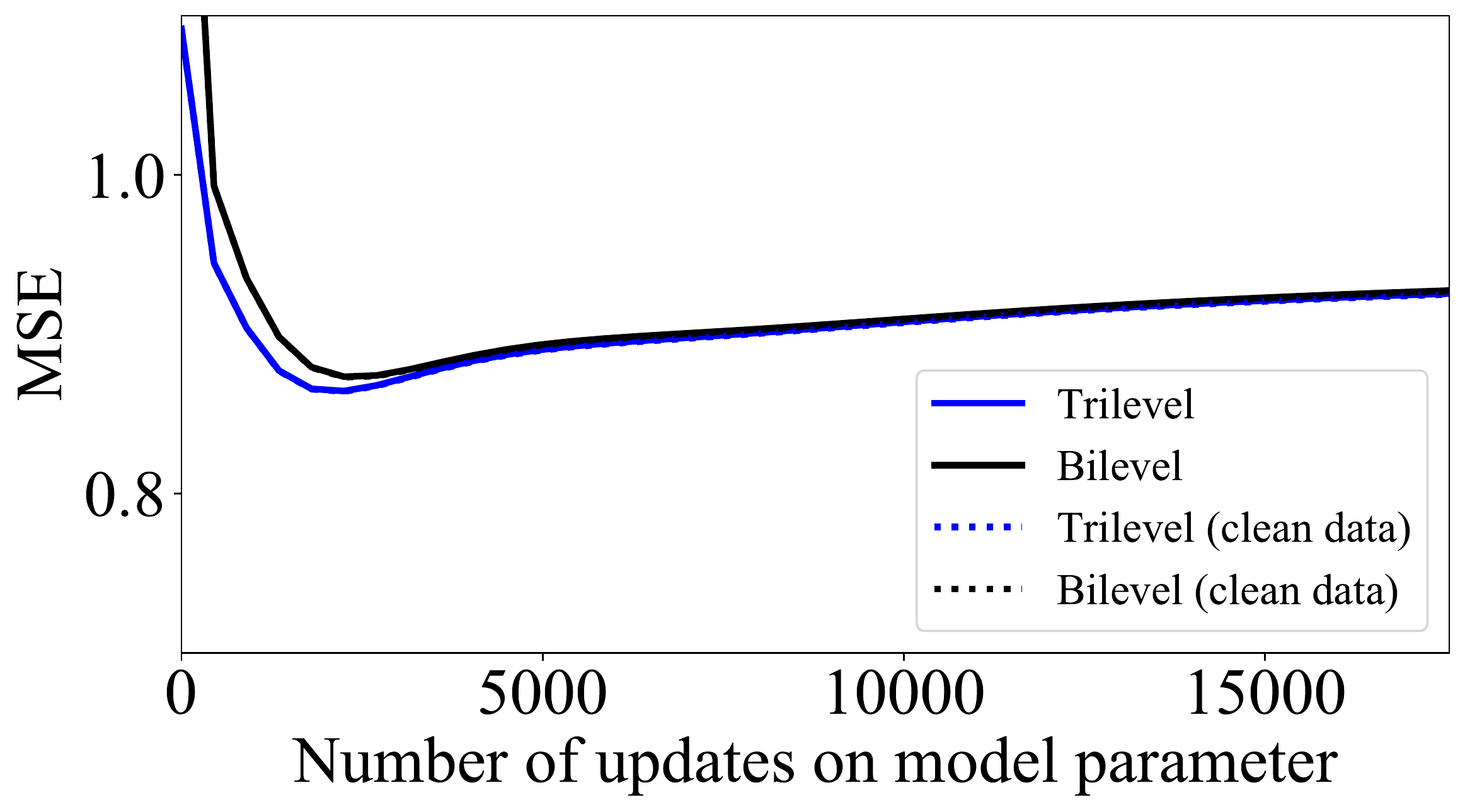}}
	\subfigure[$\sigma = 0.05$]{
	\includegraphics[width = 0.49\textwidth]{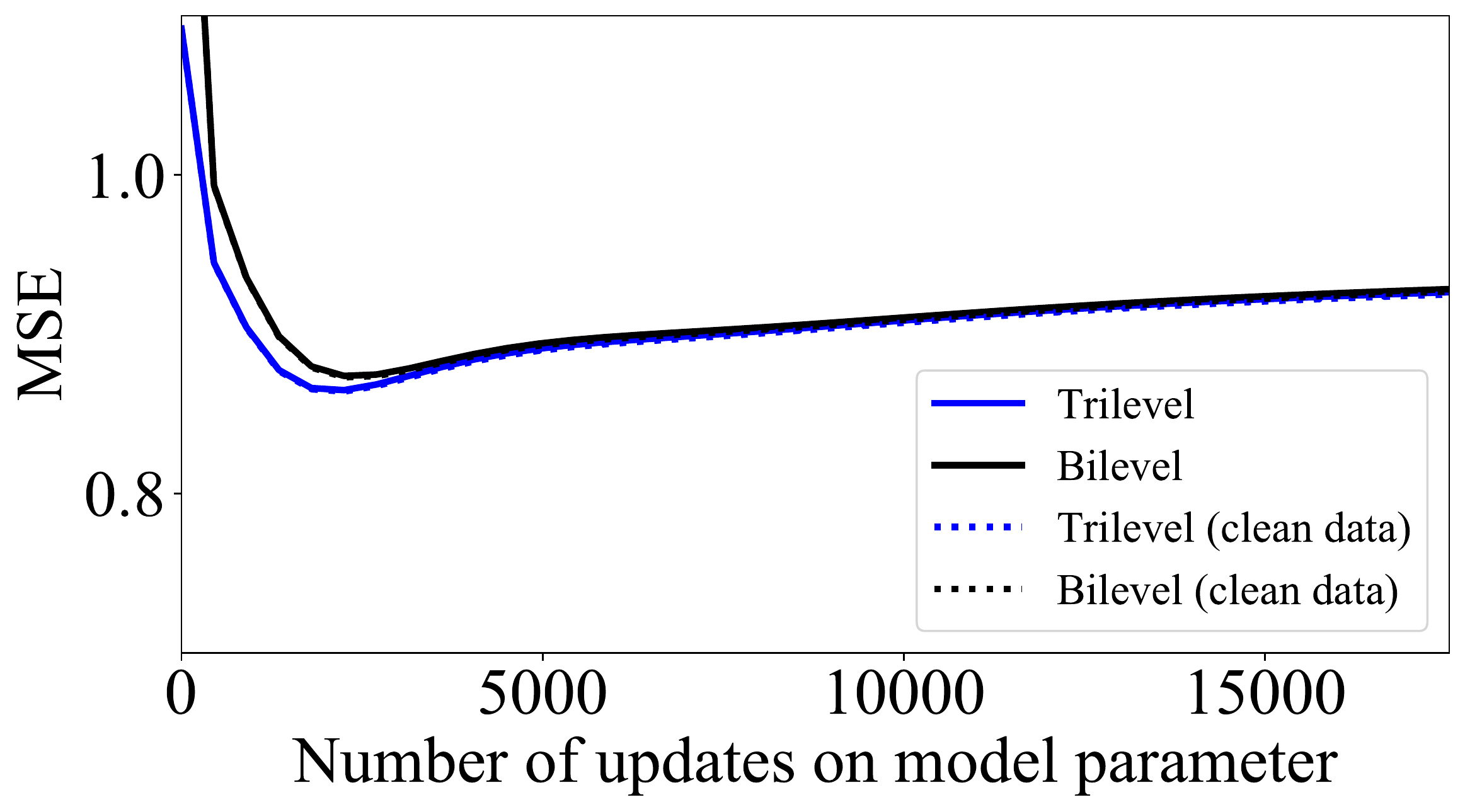}}
	\subfigure[$\sigma = 0.08$]{
	\includegraphics[width = 0.49\textwidth]{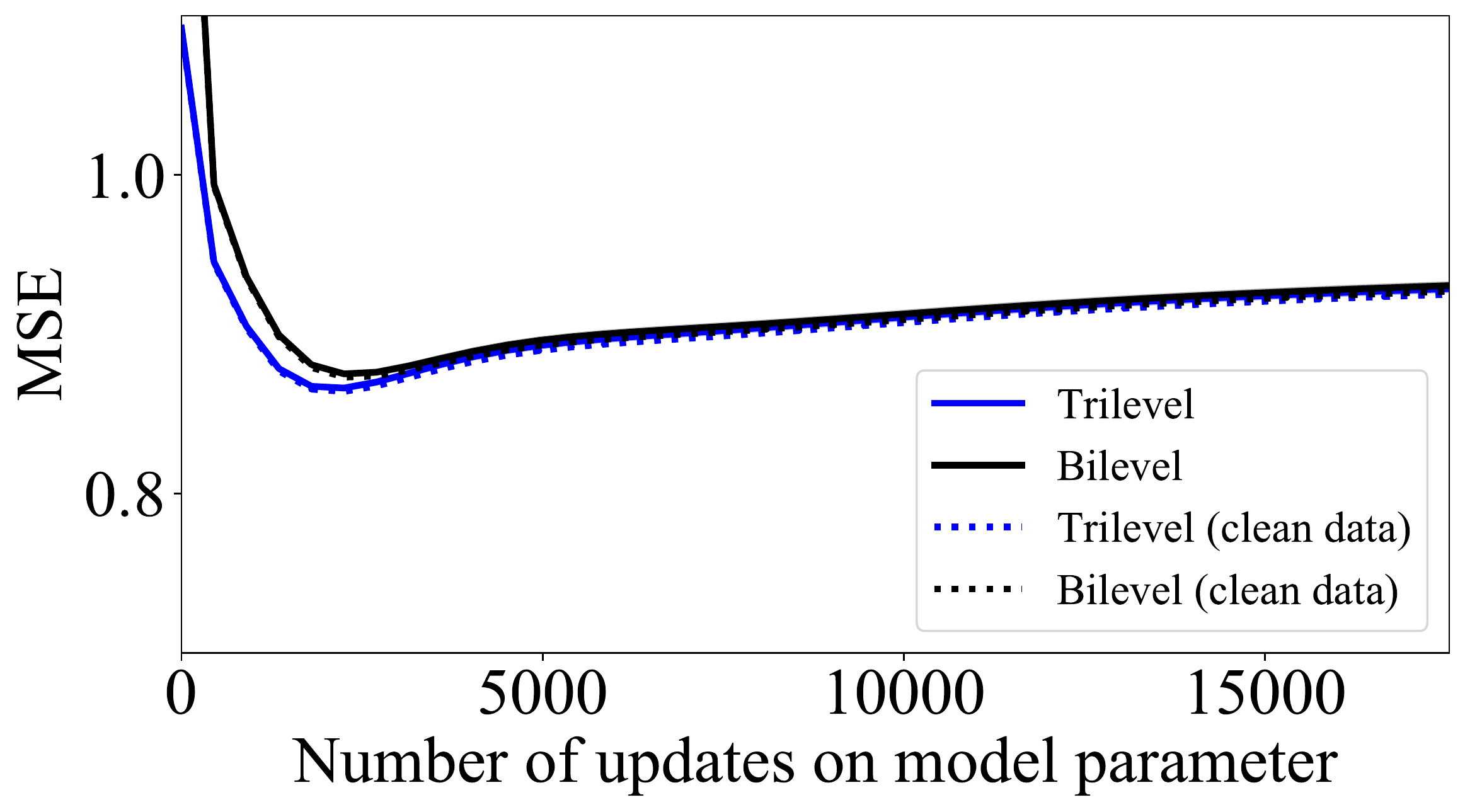}}
	\caption{MSE of test data with Gaussian noise with a standard deviation of $\sigma$ on white wine quality dataset.}
	\label{fig: suppl iter-wine_white}
\end{figure}

\begin{figure}
	\centering
	\subfigure[$\sigma = 0.01$]{
	\includegraphics[width = 0.49\textwidth]{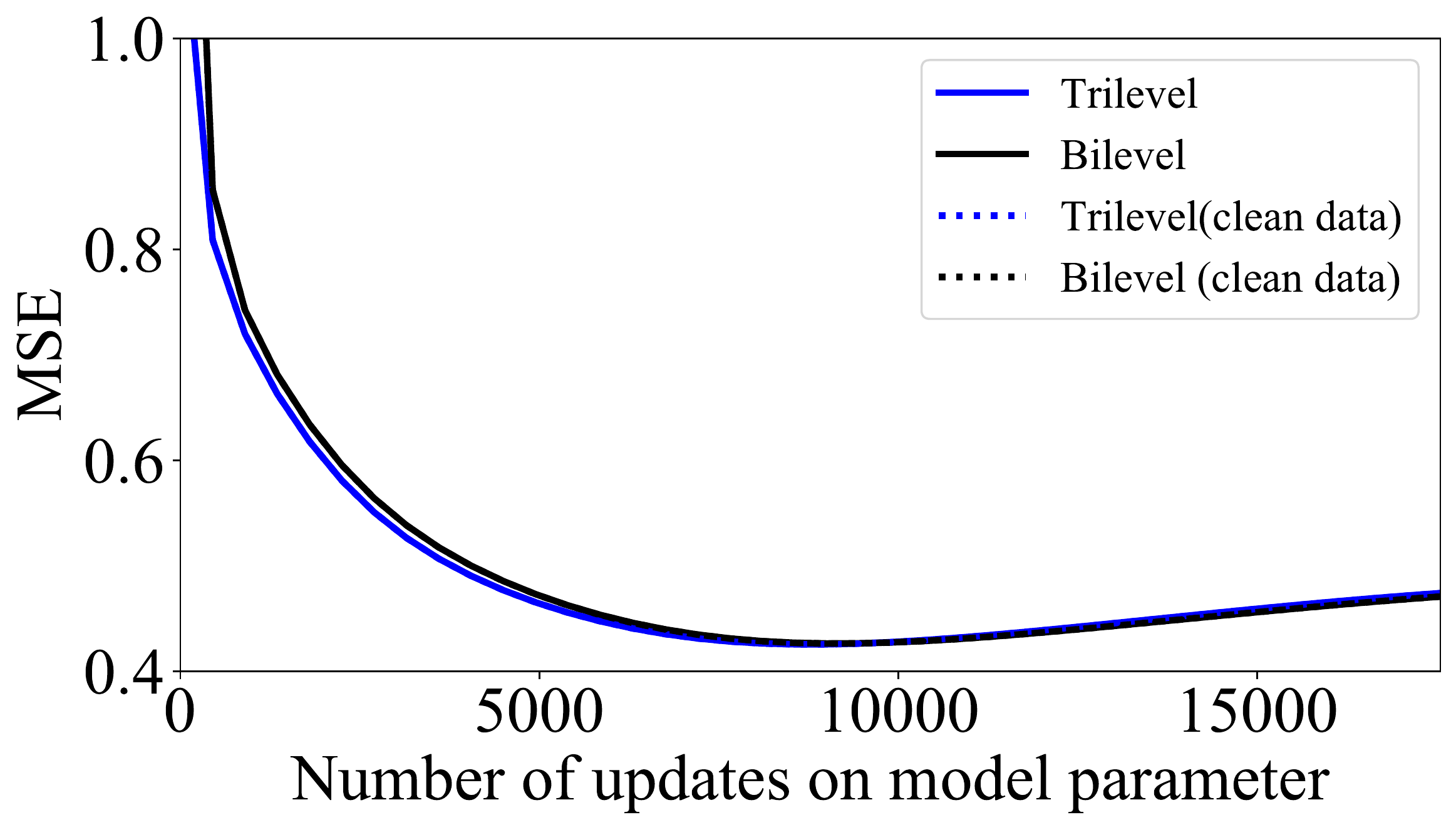}}
	\subfigure[$\sigma = 0.03$]{
	\includegraphics[width = 0.49\textwidth]{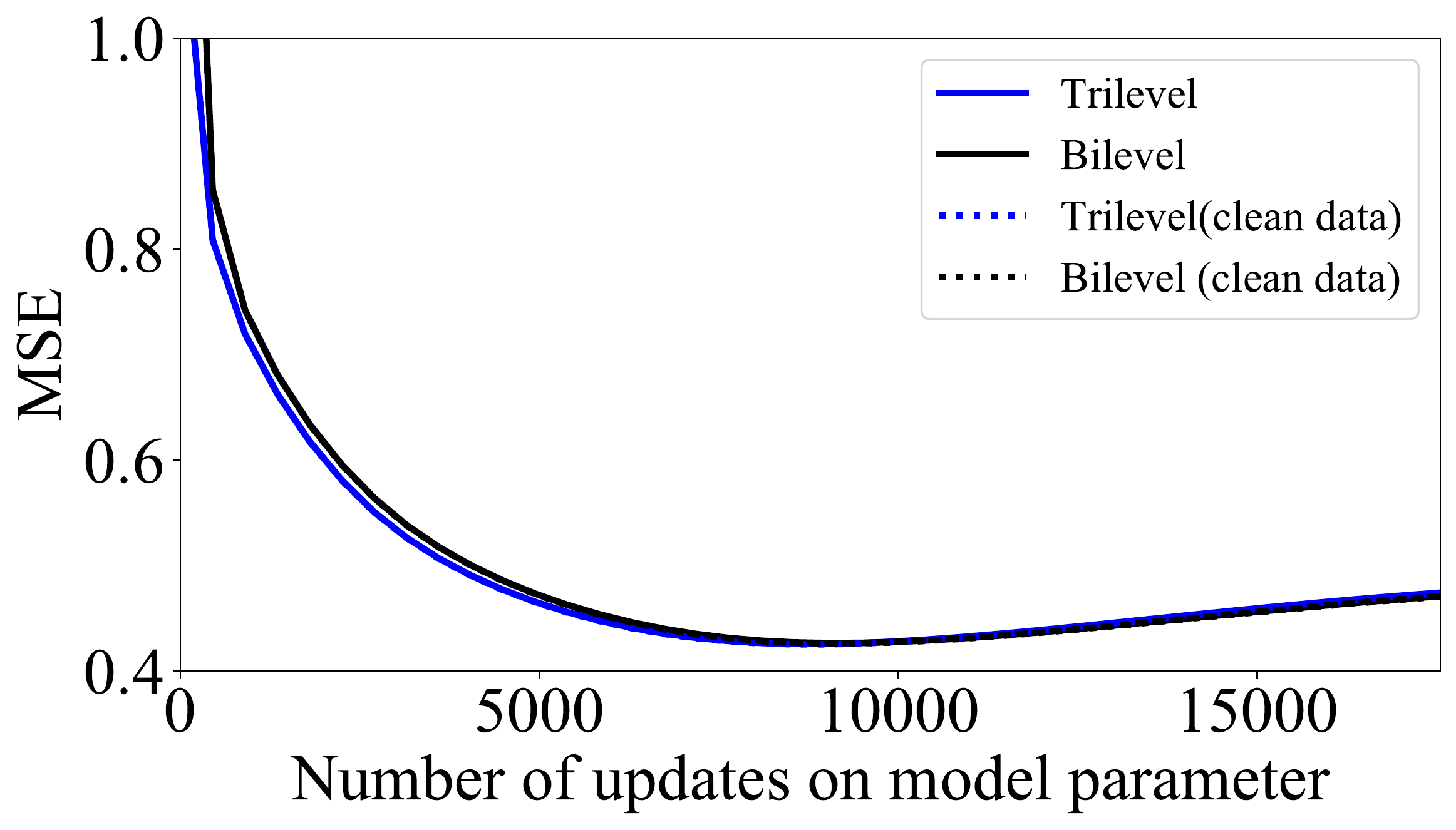}}
	\subfigure[$\sigma = 0.05$]{
	\includegraphics[width = 0.49\textwidth]{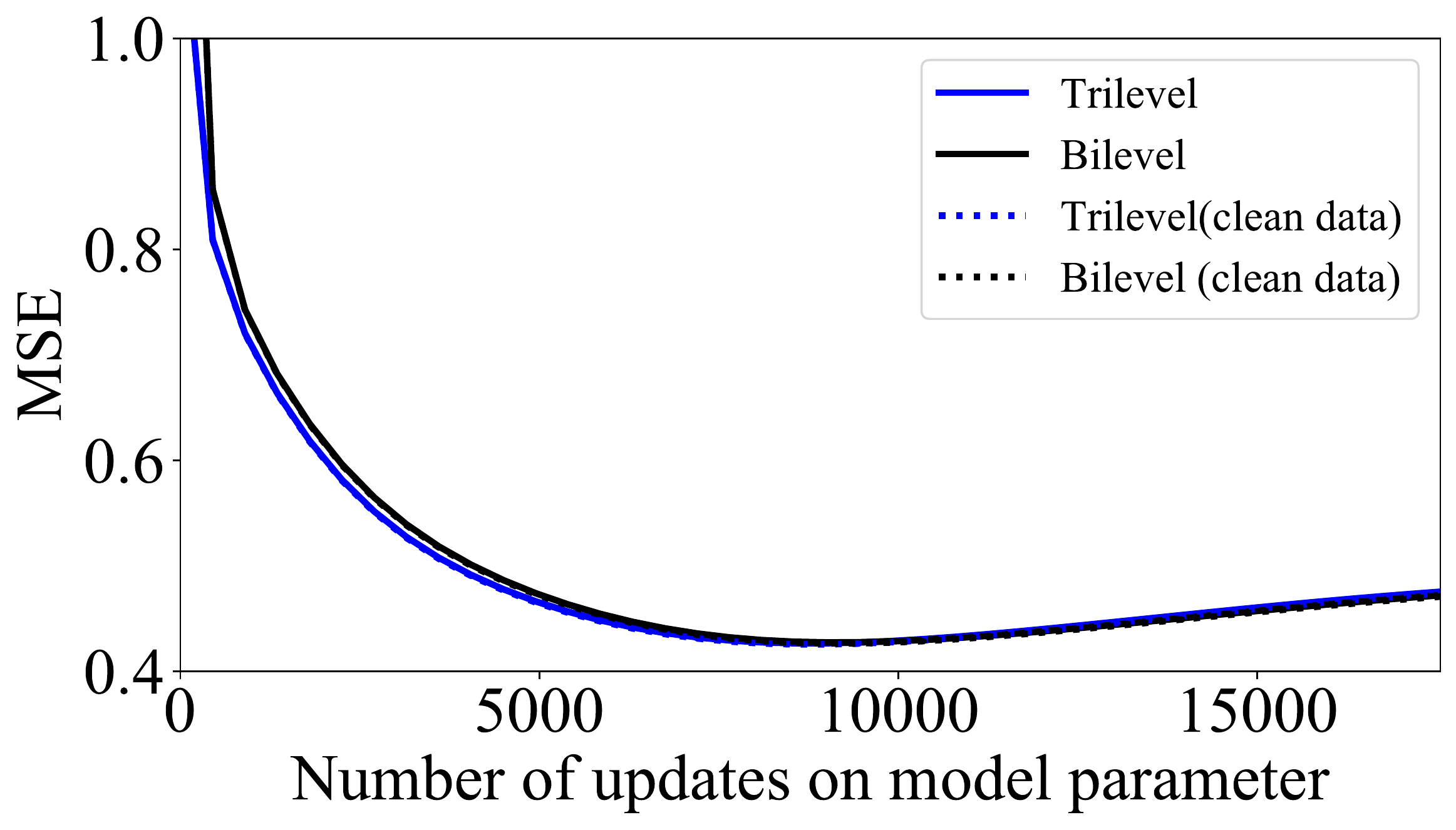}}
	\subfigure[$\sigma = 0.08$]{
	\includegraphics[width = 0.49\textwidth]{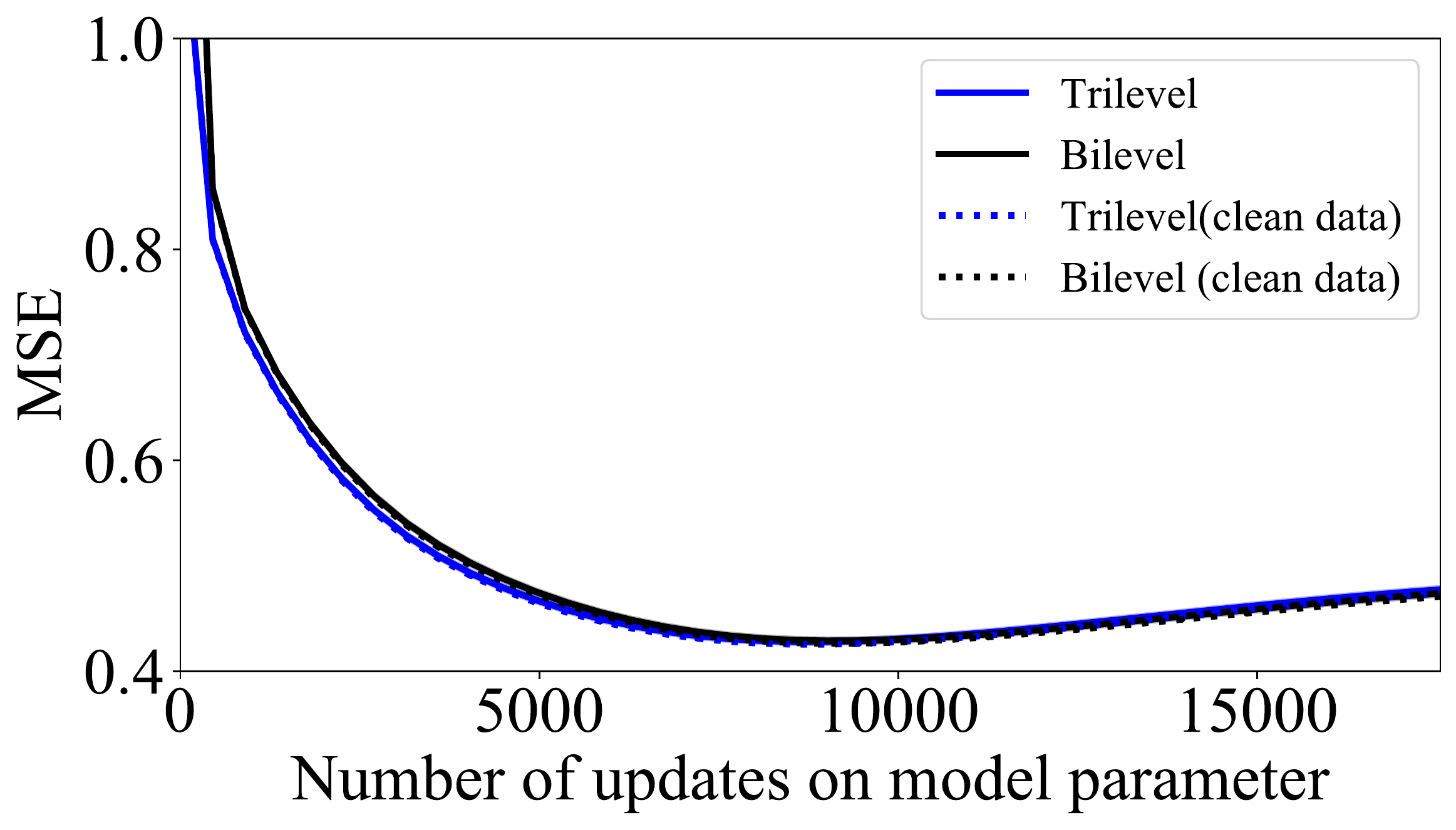}}
	\caption{MSE of test data with Gaussian noise with a standard deviation of $\sigma$ on Boston house-prices dataset.}
	\label{fig: suppl iter-boston}
\end{figure}

We show the transition of the MSE of test data with Gaussian noise
in Figures \ref{fig: suppl iter-diabetes}, \ref{fig: suppl iter-wine_red}, \ref{fig: suppl iter-wine_white},
and \ref{fig: suppl iter-boston}.
The solid line and colored belt respectively indicate the mean and the standard deviation
over 500 times of generation of Gaussian noise.
The dashed line indicates the mean squared error (MSE) without noise as a baseline.
In the results of the diabetes dataset (Figure \ref{fig: suppl iter-diabetes}),
the trilevel model provided a more robust parameter than the bilevel model,
because the MSE rises less with the large standard deviation of noise on test data.

\begin{figure}
	\centering
	\subfigure[Diabetes dataset]{
	\includegraphics[width = 0.49\textwidth]{fig_diabetes/early_stop_diabete.pdf}}
	\subfigure[Red wine quality dataset]{
	\includegraphics[width = 0.49\textwidth]{fig_wine_red/early_stop_wine_red.pdf}}
	\subfigure[White wine quality dataset]{
	\includegraphics[width = 0.49\textwidth]{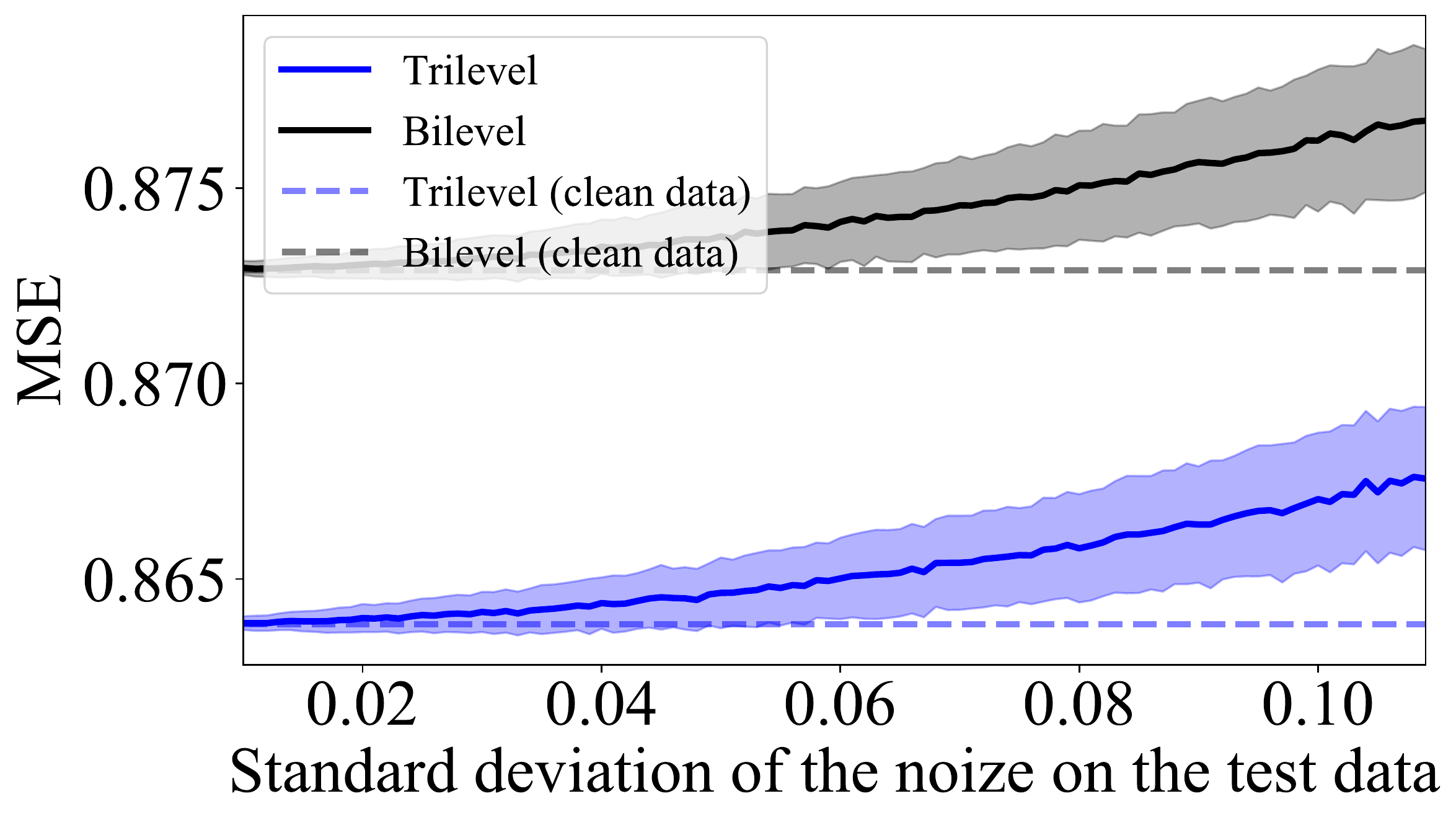}}
	\subfigure[Boston house-prices dataset]{
	\includegraphics[width = 0.49\textwidth]{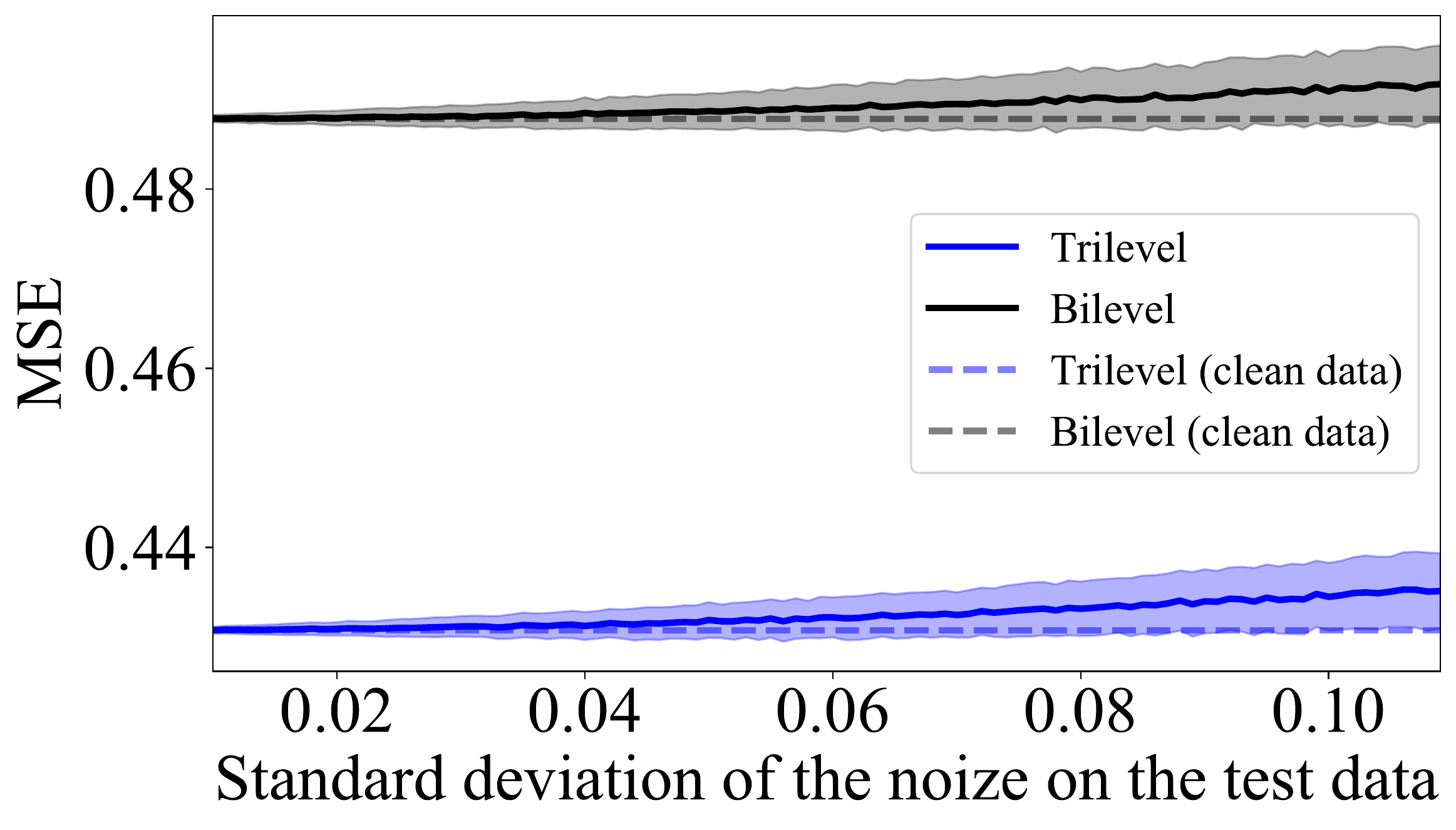}}
	\caption{MSE of test data with Gaussian noise, using early-stopped parameters for prediction.} 
	\label{fig: suppl MSE early stop}
\end{figure}

Next, we compared the quality of the resulting model parameters by the trilevel and bilevel models.
We set an early-stopping condition on learning parameters:
after 1000 times of updates on the model parameter,
if one time of update on hyperparameter $\lambda$ did not improve test error by $\epsilon$,
terminate the iteration and return the parameters at that time.
By using the early-stopped parameters obtained by this stopping condition,
we show the relationship of test error and standard deviation of the noise on the test data
in Figure \ref{fig: suppl MSE early stop} and Table \ref{tbl: suppl MSE}.
For the diabetes dataset, the growth of MSE of the trilevel model was slower than that of the bilevel model.
For wine quality and Boston house-prices datasets,
the MSE of the trilevel model was consistently lower than that of the bilevel model.
Therefore, the trilevel model provides more robust parameters than the bilevel model in these settings of problems.

\begin{table}
\centering
\caption{MSE of test data with Gaussian noise with standard deviation of 0.08, using early-stopped parameters for prediction.
	The better values are shown in bold face.}
\label{tbl: suppl MSE}
\begin{tabular}{rr@{${} \pm {}$}lr@{${} \pm {}$}lr@{${} \pm {}$}lr@{${} \pm {}$}l}
\toprule
&	\multicolumn{2}{c}{diabetes}&			\multicolumn{2}{c}{Boston}&
	\multicolumn{2}{c}{wine (red)}&			\multicolumn{2}{c}{wine (white)}\\
\midrule
\multicolumn{1}{c}{Trilevel}&
	\textbf{0.8601}&	\textbf{0.0479}&	\textbf{0.4333}&	\textbf{0.0032}&
	\textbf{0.7223}&	\textbf{0.0019}&	\textbf{0.8659}&	\textbf{0.0013}\\
\multicolumn{1}{c}{Bilevel}&
	1.0573&				0.0720&				0.4899&				0.0033&
	0.7277&				0.0019&				0.8750&				0.0014\\
\bottomrule
\end{tabular}
\end{table}